# A Particle Method for 1-D Compressible Fluid Flow

**Iasson Karafyllis[*] and Markos Papageorgiou[**],[a]**

[*]Dept. of Mathematics, National Technical University of Athens,
Zografou Campus, 15780, Athens, Greece.

[**]Dynamic Systems and Simulation Laboratory, Technical University of Crete,
Chania, 73100, Greece.

[a]Faculty of Maritime and Transportation, Ningbo University, Ningbo, China.

**Abstract**

This paper proposes a novel particle scheme that provides convergent approximations of a weak solution of the Navier-Stokes equations for the 1-D flow of a viscous compressible fluid. Moreover, it is shown that all differential inequalities that hold for the fluid model are preserved by the particle method: mass is conserved, mechanical energy is decaying, and a modified mechanical energy functional is also decaying. The proposed particle method can be used both as a numerical method and as a method of proving existence of solutions for compressible fluid models.

## 1. Introduction

The study of the flow of viscous compressible fluids is a major research topic that has attracted the attention of many researchers (see for instance [7, 10, 21, 22, 23, 24, 25, 27, 28, 29]). The Navier-Stokes equations that are used for the description of the flow of a viscous fluid have been studied in detail for existence/uniqueness issues, as well as for the dynamical behavior of the solutions and the numerical approximation of the solutions (see for instance [7] for a description of finite-volume and finite-elements methods).

Among the numerical methods used for the approximation of the solutions of the Navier-Stokes equations we should mention the so-called "particle methods", which is also known as the method of "Smoothed Particle Hydrodynamics" (see [8, 26, 33] and references therein) and has been used extensively for many difficult flow problems. Recently, particle methods have been used for other Partial Differential Equations (PDEs) than the Navier-Stokes equations (see [3, 6]). It can be said that particle methods provide the link between the macroscopic fluid model and the microscopic ODE model that describes the interaction between the fluid particles. This is important from a mathematical viewpoint (Hilbert's 6[th] problem), even when the forces acting on the particles are fictitious. Moreover, particle methods have been used recently for the proof of existence of entropy solutions for traffic flow problems (see [4, 11]). The use of particle methods in traffic flow problems is important due to a real-world application: the design and use of cruise controllers for automated vehicles. Cruise controllers can manipulate the accelerations of the vehicles and implement forces that may be considered as fictitious forces from a physical point of view. When the vehicles are considered as (self-driven) particles of a fluid ("the traffic fluid") the corresponding macroscopic PDE models are very similar to the equations of a viscous compressible fluid (see [14, 32]). It then becomes clear that the questions of convergence and accuracy are crucial.

On the other hand, feedback controllers have been proposed for 1-D compressible fluid models (and the related 1-D Saint-Venant models for incompressible fluids), which are based on proving differential inequalities for (Control) Lyapunov functionals (not only mass and mechanical energy, but also for other functionals; see [15, 16, 17, 18]). However, there is no guarantee that the differential inequalities hold (in a discretized form) for the numerical approximation. This is a well-known problem (see [9, 13, 30] for the analysis of the problem in a finite-dimensional setting), which is closely related to the issue of numerical stability for the applied scheme. To our knowledge, the issue has not been studied for particle methods.

The present work provides rigorous answers to all aforementioned issues. Adopting a methodology that is inspired from [4, 11, 12], we construct a novel particle scheme. The particle method is shown to provide approximations that converge to a weak solution of the Navier-Stokes equations for the 1-D flow of a viscous compressible fluid. The results of the paper cover the case where the viscosity depends on the density of the fluid. Moreover, we show that all differential inequalities that hold for the PDE model are preserved by the particle method: mass is conserved, mechanical energy is decaying, and a modified mechanical energy functional used recently in [15, 16, 17, 18] is also decaying. We see the proposed particle method both as a numerical method (because Lemma 4.3 and Lemma 4.4 below provide sharp estimates of the numerical accuracy of the scheme) and as a method of proving existence of solutions for compressible fluid models. In particular, it should be noticed that, when the particle method is considered as a method of proving existence of solutions for compressible fluid models then the method has a similar structure with the so-called "method of lines" (see [10, 24]). However, there are also crucial differences between the proposed particle method and the method of lines: in the method of lines the states (density and velocity) are approximated at fixed grid points in Lagrangian coordinates while in the particle method shown in this work the states (position and velocity of the particles) are evaluated at time-varying positions (the position of the particles) in Eulerian coordinates.

The structure of the paper is as follows. Section 2 is devoted to the presentation of the problem. Section 3 contains the statements of the assumptions and the main results of the paper. The proofs of all results are provided in Section 4. Finally, Section 5 gives the concluding remarks of the present work.

**Notation.** Throughout this paper, we adopt the following notation.

* $\mathbb{R}_+ = [0, +\infty)$ denotes the set of non-negative real numbers.

* Let $X$ be a given normed linear space with norm $\|\cdot\|_X$, let $B \subseteq X$ be a subset of $X$ and let $I \subseteq \mathbb{R}$ be a non-empty interval. By $C^0(I;B)$ we denote the class of continuous mappings $f: I \to B$, i.e., the class of mappings $f: I \to B$ for which the following property holds: for every $t_0 \in I$ and for every $\varepsilon > 0$ there exists $\delta > 0$ such that $\|f(t) - f(t_0)\|_X < \varepsilon$ for all $t \in I \cap (t_0 - \delta, t_0 + \delta)$. By $C_c^0(I;X)$ we denote the class of continuous mappings $f: I \to X$ with compact support. By $C_c^1(I)$, where $I \subseteq \mathbb{R}$ is an open interval, we denote the class of continuously differentiable functions $f: I \to \mathbb{R}$ with compact support. For any $\gamma \in (0,1]$, by $C^{0,\gamma}(I;B)$ we denote the class of Hölder continuous mappings $f: I \to B$ with exponent

$\gamma \in (0,1]$, i.e., the class of mappings $f : I \to B$ for which
$$\sup\left\{ \frac{\|f(t)-f(t_0)\|_X}{|t-t_0|^\gamma} : t, t_0 \in I, t \neq t_0 \right\} < +\infty.$$

* Let $S \subseteq \mathbb{R}^n$ be an open set and let $A \subseteq \mathbb{R}^n$ be a set that satisfies $S \subseteq A \subseteq cl(S)$. By $C^0(A;\Omega)$, we denote the class of continuous functions on $A$, which take values in $\Omega \subseteq \mathbb{R}^m$. By $C^k(A;\Omega)$, where $k \geq 1$ is an integer, we denote the class of functions on $A \subseteq \mathfrak{R}^n$, which take values in $\Omega \subseteq \mathfrak{R}^m$ and have continuous derivatives of order $k$. In other words, the functions of class $C^k(A;\Omega)$ are the functions which have continuous derivatives of order $k$ in $S = \text{int}(A)$ that can be continued continuously to all points in $\partial S \cap A$. When $\Omega = \mathbb{R}$ then we write $C^0(A)$ or $C^k(A)$. When $I \subseteq \mathbb{R}$ is an interval and $\eta \in C^1(I)$ is a function of a single variable, $\eta'(\rho)$ denotes the derivative with respect to $\rho \in I$.

* Let $I \subseteq \mathbb{R}$ be an interval and let $a < b$ be given constants. Let the function $u : I \times (a,b) \to \mathbb{R}$ be given. We use the notation $u[t]$ to denote the profile at certain $t \in I$, i.e., $(u[t])(x) = u(t,x)$ for all $x \in (a,b)$.

* We denote by $K_\infty$ the class of continuous, increasing functions $a : \mathbb{R}_+ \to \mathbb{R}_+$ with $a(\mathbb{R}_+) = \mathbb{R}_+$.

* Let $\Omega \subseteq \mathbb{R}^n$ be a given open set given. For $p \in [1,+\infty)$, $L^p(\Omega)$ is the set of equivalence classes of Lebesgue measurable functions $u : \Omega \to \mathbb{R}$ with $\|u\|_p := \left( \int_\Omega |u(x)|^p \, dx \right)^{1/p} < +\infty$. $L^\infty(\Omega)$ is the set of equivalence classes of Lebesgue measurable functions $u : \Omega \to \mathbb{R}$ with $\|u\|_\infty := \operatorname*{ess\,sup}_{x \in \Omega}(|u(x)|) < +\infty$.

* Let $X$ be a given Banach space with norm $\|\cdot\|_X$ and let $I \subseteq \mathbb{R}$ be a given open interval. A function $f : I \to X$ is measurable if there exists a set $E \subset I$ of measure zero and a sequence $\{f_n \in C_c^0(I;X) : n \geq 1\}$ such that $\lim_{n \to +\infty}(f_n(t)) = f(t)$ for all $t \in I \setminus E$. For $p \in [1,+\infty)$, $L^p(I;X)$ is the set of equivalence classes of measurable functions $f : I \to X$ with $\|f\|_p := \left( \int_I \|f(t)\|_X^p \, dx \right)^{1/p} < +\infty$. $L^\infty(I;X)$ is the set of equivalence classes of measurable functions $f : I \to X$ with $\|f\|_\infty := \operatorname*{ess\,sup}_{t \in I}(\|f(t)\|_X) < +\infty$.

* Let $a < b$ be given constants. For $p \in [1,+\infty]$, $W^{1,p}(a,b)$ denotes the Sobolev space of functions in $L^p(a,b)$ with weak derivative in $L^p(a,b)$. We set $H^1(a,b) = W^{1,2}(a,b)$. The closure of $C_c^1((a,b))$ in $H^1(a,b)$ is denoted by $H_0^1(a,b)$. The dual space of $H_0^1(a,b)$ is denoted by $H^{-1}(a,b)$.

## 2. Description of the Problem

Let $L > 0$ be a given constant. We consider the following model

$$\rho_t + (\rho v)_x = 0, \text{ for } t > 0, \ x \in [0, L] \tag{2.1}$$

$$(\rho v)_t + \left(\rho v^2 + P(\rho)\right)_x = \left(\mu(\rho) v_x\right)_x, \text{ for } t > 0, \ x \in (0, L) \tag{2.2}$$

$$v(t, 0) = v(t, L) = 0, \text{ for } t \geq 0 \tag{2.3}$$

where $\rho(t, x) > 0$, $v(t, x) \in \mathbb{R}$, and the functions $P, \mu : (0, +\infty) \to (0, +\infty)$ are $C^2((0, +\infty))$ functions with $P'(\rho) > 0$ for all $\rho > 0$.

Mathematical model (2.1), (2.2), (2.3) may describe:

- the 1-D isentropic (or polytropic) motion of a compressible Newtonian fluid contained between two walls placed at $x = 0$ and $x = L$; in this case $\rho(t, x) > 0$ and $v(t, x) \in \mathbb{R}$ are the fluid density and fluid velocity, respectively, at time $t \geq 0$ and position $x \in [0, L]$, the dynamic viscosity of the fluid is $\mu(\rho)$ and the pressure is given by the equation $P(\rho) = k\rho^\gamma$ with $k > 0, \gamma > 1$,
- the motion of an incompressible Newtonian liquid within a tank of length $L > 0$ (case of viscous Saint-Venant model; see [16, 17, 18, 20, 31]); in this case, $\rho(t, x) > 0$ and $v(t, x) \in \mathbb{R}$ are the liquid level and the liquid velocity, respectively, at time $t \geq 0$ and position $x \in [0, L]$, while $P(\rho) := \frac{1}{2} g \rho^2$ and $\mu(\rho) := \mu \rho$, where $g, \mu > 0$ (constants), are the acceleration of gravity and the kinematic viscosity of the liquid.

Using (2.1) and (2.3), we can prove that for every classical solution of (2.1)-(2.3) it holds that $\frac{d}{dt}\left(\int_0^L \rho(t, x) dx\right) = 0$ for all $t > 0$. Therefore, the total mass of the fluid $m > 0$ is constant. Thus, without loss of generality, we assume that every solution of (2.1)-(2.3) satisfies the equation

$$\int_0^L \rho(t, x) dx \equiv m \tag{2.4}$$

where $m > 0$ is a given constant. It should be emphasized that for obvious physical reasons, $\rho(t, x)$ must be positive for all times, i.e., we must have:

$$\min_{x \in [0, L]} (\rho(t, x)) > 0, \text{ for } t \geq 0. \tag{2.5}$$

System (2.1)-(2.4) allows a unique equilibrium point, namely the point

$$\rho(x) \equiv \rho^*, \ v(x) \equiv 0, \text{ for } x \in [0, L] \tag{2.6}$$

where $\rho^* = m / L$.

Define the functions

$$k(\rho) := \int_{\rho^*}^{\rho} \frac{\mu(\tau)}{\tau} d\tau \qquad (2.7)$$

$$Q(\rho) := \rho \int_{\rho^*}^{\rho} \frac{P(\tau)}{\tau^2} d\tau - \frac{P(\rho^*)}{\rho^*} \rho + P(\rho^*). \qquad (2.8)$$

Since $Q''(\rho) = \rho^{-1} P'(\rho) > 0$ for all $\rho > 0$, $Q(\rho^*) = Q'(\rho^*) = 0$, it follows that $Q(\rho) > 0$ for all $\rho > 0$, $\rho \neq \rho^*$.

Given $\rho \in C^0([0, L]; (0, +\infty)) \cap H^1((0, L))$, $v \in L^2((0, L))$ with $\int_0^L \rho(x) dx = m$, we define the following functionals:

$$E(\rho, v) := \frac{1}{2} \int_0^L \rho(x) v^2(x) dx + U(\rho) \qquad (2.9)$$

$$W(\rho, v) := \frac{1}{2} \int_0^L \frac{1}{\rho(x)} \left( \rho(x) v(x) + \frac{\partial}{\partial x} (k(\rho(x))) \right)^2 dx + U(\rho) \qquad (2.10)$$

$$U(\rho) := \int_0^L Q(\rho(x)) dx. \qquad (2.11)$$

We notice that:
- the functional $U$ is the potential energy of the fluid;
- the functional $E$ is the total mechanical energy of the two-piston system. Indeed, notice that $E$ is the sum of the potential energy ($U$) and the kinetic energy of the fluid ($\frac{1}{2} \int_0^L \rho(x) v^2(x) dx$);
- the functional $W$ is a modified mechanical energy and has been constructed based on a specific transformation that has been used extensively in the literature of isentropic, compressible fluid flow (see [19, 27, 28, 29]). More specifically, the transformation $w = \rho v + (k(\rho))_x$ has the effect of making the viscosity term disappear from the momentum equation, which then becomes $w_t + (vw)_x = -(P(\rho))_x$. The modified mechanical energy $W$ has been used recently for the construction of stabilizing feedback laws in control problems (see [15, 16, 17, 18]).

The functionals $E, W$ are non-increasing along the classical solutions of (2.1)-(2.4). More specifically, we have the following lemma. Its proof is straightforward and is omitted.

**Lemma 1.1:** *For every classical solution of the PDE system (2.1)-(2.4) the following equations hold for all $t > 0$:*

$$\frac{d}{dt} E(\rho[t], v[t]) = -\int_0^L \mu(\rho(t, x)) v_x^2(t, x) dx \qquad (2.12)$$

$$\frac{d}{dt} W(\rho[t], v[t]) = -\int_0^L \rho^{-2}(t, x) P'(\rho(t, x)) \mu(\rho(t, x)) \rho_x^2(t, x) dx \qquad (2.13)$$

*where $E, W$ are defined by (2.9), (2.10), respectively.*

In what follows, we assume the following property for the pressure function $P$.

**(A)** $\lim\limits_{\rho \to +\infty} \left( \int\limits_{\rho^*}^{\rho} s^{-2} P(s) ds \right) = +\infty$ and $\inf\limits_{\rho > 0} \left( \int\limits_{\rho^*}^{\rho} s^{-2} P(s) ds \right) > -\infty$.

Assumption (A) holds automatically for the viscous Saint-Venant model as well as for the isentropic (or polytropic) motion of a compressible Newtonian fluid since every pressure function of the form $P(\rho) = k\rho^\gamma$ with $k > 0, \gamma > 1$ satisfies Assumption (A).

We next provide the precise notion of a weak solution for the problem (2.1)-(2.4).

**Definition 1.2:** *Let* $\rho_0 \in W^{1,\infty}((0,L)) \cap C^0([0,L];(0,+\infty))$ *with* $\int\limits_0^L \rho_0(x)dx = m$, $v_0 \in H_0^1((0,L))$ *be given. A weak solution of the initial-boundary value problem (2.1)-(2.4) with*

$$\rho[0] = \rho_0 \quad , \quad v[0] = v_0 \tag{2.14}$$

*on the interval* $[0,T]$, *where* $T > 0$, *is a pair of functions*

$$\rho \in L^\infty((0,T); H^1((0,L))) \cap C^{0,1/2}([0,T]; L^2((0,L))),$$
$$v \in L^\infty((0,T); H_0^1((0,L))) \cap C^{0,1/4}([0,T]; L^2((0,L)))$$
$$\text{with } \rho_x \in L^\infty((0,T) \times (0,L)), \; \rho_t \in L^\infty((0,T); L^2((0,L))),$$
$$v_t \in L^2((0,T); H^{-1}((0,L)))$$
$$\text{and } \rho[t], v[t] \in L^\infty((0,L)) \text{ for all } t \in [0,T]$$

*for which there exist constants* $\rho_{\max} > \rho_{\min} > 0$, $v_{\max} > 0$, *such that (2.4) holds for all* $t \in [0,T]$ *and the following additional conditions hold:*

$$\int\limits_0^L \varphi(0,x)\rho_0(x)dx + \int\limits_0^T \int\limits_0^L \rho(t,x)(\varphi_t(t,x) + v(t,x)\varphi_x(t,x))dxdt = 0$$

for every $\varphi \in C^1([0,T] \times [0,L])$ with $\varphi(T,x) = 0$ for all $x \in [0,L]$ \hfill (2.15)

$$\int\limits_0^L \varphi(0,x)\rho_0(x)v_0(x)dx + \int\limits_0^T \int\limits_0^L \varphi_t(t,x)\rho(t,x)v(t,x)dxdt$$
$$+ \int\limits_0^T \int\limits_0^L \varphi_x(t,x)\left(\rho(t,x)v^2(t,x) + P(\rho(t,x)) - \mu(\rho(t,x))v_x(t,x)\right)dxdt = 0$$

for every $\varphi \in C^2([0,T] \times [0,L])$ with $\varphi(T,x) = 0$ for all $x \in [0,L]$ and
$$\varphi(t,0) = \varphi(t,L) = \varphi_x(t,L) = 0 \text{ for all } t \in [0,T] \tag{2.16}$$

$$E(\rho[t], v[t]) \leq E(\rho[t_0], v[t_0]), \text{ for all } T \geq t \geq t_0 \geq 0 \tag{2.17}$$

$$\frac{1}{h}\int_{t}^{t+h}W(\rho[s],v[s])ds \leq mL\|v_0'\|_2^2 + 2m^5\bar{M}^2\frac{\|\rho_0'\|_\infty^2}{\rho_{\min}^6} + m\Phi\left(\frac{m}{\|\rho_0\|_\infty}\right),$$

for all $t \in [0,T)$ and $h > 0$ with $t+h \leq T$ and $\bar{M} = \max\left\{K'(s): \frac{m}{\|\rho_0\|_\infty} \leq s \leq \frac{m}{\rho_{\min}}\right\}$ (2.18)

$$\rho_{\min} \leq \rho(t,x) \leq \rho_{\max}, \text{ for } (t,x) \in (0,T)\times(0,L) \text{ a.e.} \quad (2.19)$$

$$\rho_{\min} \leq \rho(t,x) \leq \rho_{\max}, \text{ for all } t \in [0,T] \text{ and for } x \in (0,L) \text{ a.e.} \quad (2.20)$$

$$|v(t,x)| \leq v_{\max}, \text{ for } (t,x) \in (0,T)\times(0,L) \text{ a.e.} \quad (2.21)$$

$$|v(t,x)| \leq v_{\max}, \text{ for all } t \in [0,T] \text{ and for } x \in (0,L) \text{ a.e.} \quad (2.22)$$

**Remarks on Definition 1.2:** **(i)** The main difference of the notion of weak solution introduced in the present work and other notions used in the literature (see [7, 10, 21, 25]) is condition (2.18). Indeed, while some works are using the mechanical energy functional $E$ for the definition of a weak solution of (2.1)-(2.4) (see for instance [7]), we could not find any other work that uses both the mechanical energy functional $E$ and the modified mechanical energy functional $W$. However, as noticed above, the use of both functionals is important for control purposes (see [15, 16, 17, 18]).
**(ii)** Another difference of the notion of weak solution introduced in the present work and other notions used in the literature is the fact that we require the density to be positive (recall (2.19), (2.20)). Other works in the literature (see [7, 10, 25]) simply require the density to be non-negative. The difference is crucial: we do not allow vacuum to be formed while in other works in the literature vacuum is allowed. Therefore, the assumptions in the main results that are presented below (namely, Assumption (A) and inequality (3.30) below) can be considered to be sufficient conditions for the exclusion of vacuum. This is a useful point when one faces the question of vacuum: vacuum can only appear when the existence of a weak solution is verified using relevant results in the literature that allow density to be zero and when the assumptions of the present work are violated.
**(iii)** Another difference of the notion of weak solution introduced in the present work and other notions used in the literature is the regularity guaranteed for the solution and the regularity properties required for the initial condition. However, it should be noted that the notion of weak solution of Definition 1.2 requires less regularity than the notion of strong solution used in [22].

## 3. A Novel Particle Method

<u>3.I. Equations</u>

We consider $n \geq 2$ particles of fluid spaced at $x_i(t) \in [0,L]$, $i = 1,...,n$, each with mass $\frac{m}{n}$. We also consider a "ghost" particle at $x_0(t) \equiv L$ and we assume that $x_n(t) = 0 < x_{n-1}(t) < ... < x_1(t) < x_0(t) = L$. The ghost particle has no mass and does not move, i.e.,

$$\dot{x}_0(t) = v_0(t) \equiv 0 \quad (3.1)$$

while all other particles move according to the following equations for $i = 1,...,n-1$:

$$\dot{x}_i(t) = v_i(t)$$
$$\dot{v}_i(t) = n\Phi'(n(x_{i-1}(t) - x_i(t))) - n\Phi'(n(x_i(t) - x_{i+1}(t)))$$
$$+ n^2 K'(n(x_{i-1}(t) - x_i(t)))(v_{i-1}(t) - v_i(t)) + n^2 K'(n(x_i(t) - x_{i+1}(t)))(v_{i+1}(t) - v_i(t))$$
(3.2)

and the equation
$$\dot{x}_n(t) = v_n(t) \equiv 0.$$
(3.3)

The functions $\Phi, K : (0, +\infty) \to \mathbb{R}$ are selected below in an appropriate way. The state space for the dynamical system (3.2) is the open set

$$\Omega_n := \left\{ (x_1, ..., x_{n-1}, v_1, ..., v_{n-1}) \in (0, L)^{n-1} \times \mathbb{R}^{n-1} : x_{n-1} < ... < x_1 \right\}.$$
(3.4)

If we follow the methodology in [14] and relate finite-dimensional system (3.1), (3.2), (3.3) with the infinite-dimensional system (2.1)-(2.4) and if we follow the approximation

$$\rho(x_i) \approx \frac{m}{n(x_{i-1} - x_i)}, \quad f\left(\frac{m}{\rho(x_{i+1})}\right) \approx f\left(\frac{m}{\rho(x_i)}\right) + f'\left(\frac{m}{\rho(x_i)}\right)\frac{m^2}{n\rho^3(x_i)}\rho_x(x_i)$$

then the functions $\Phi, K : (0, +\infty) \to \mathbb{R}$ must satisfy the following relations for all $x > 0$

$$\Phi''(x) = x^{-2} P'\left(\frac{m}{x}\right) > 0$$
$$\Phi'(x) = -\frac{1}{m} P\left(\frac{m}{x}\right) < 0$$
(3.5)
$$K'(x) = \frac{1}{mx} \mu\left(\frac{m}{x}\right) > 0$$

as well as the equation
$$\Phi(L) = 0.$$
(3.6)

Consequently, by virtue of (2.7) and (3.5), (3.6), we get for all $x > 0$:

$$\Phi(x) = \int_{\rho^*}^{m/x} s^{-2} P(s) \, ds$$
$$K(x) = -\frac{1}{m} k\left(\frac{m}{x}\right)$$
(3.7)

The potential function $\Phi$ is related to the function $Q$ defined by (2.8) that gives the potential energy of the fluid $U(\rho) = \int_0^L Q(\rho(x))dx$ by means of the following formula that holds for all $\rho > 0$:

$$Q(\rho) = \rho \Phi\left(\frac{m}{\rho}\right) - P(\rho^*)\left(\frac{\rho}{\rho^*} - 1\right).$$
(3.8)

Assumption (A) guarantees the following property for the potential function $\Phi$ (recall (3.7)):

$$\lim_{x \to 0^+} (\Phi(x)) = +\infty \text{ and } \Phi:(0,+\infty) \to \mathbb{R} \text{ is bounded from below, i.e., } \inf_{x>0} (\Phi(x)) > -\infty. \quad (3.9)$$

We define the following functions for $n \geq 2$:

$$E_n(t) = \frac{m}{2n} \sum_{i=1}^{n} v_i^2(t) + \frac{m}{n} \sum_{i=1}^{n} \Phi\left(n\left(x_{i-1}(t) - x_i(t)\right)\right) \quad (3.10)$$

$$W_n(t) = \frac{m}{2n} \sum_{i=1}^{n-1} w_i^2(t) + \frac{m}{n} \sum_{i=1}^{n} \Phi\left(n\left(x_{i-1}(t) - x_i(t)\right)\right) \quad (3.11)$$

$$Z_n(t) = \frac{1}{2} \sum_{i=1}^{n} \frac{(v_{i-1}(t) - v_i(t))^2}{x_{i-1}(t) - x_i(t)} \quad (3.12)$$

$$H_n(t) = \max_{i=1,\ldots,n-1} \left( n \left| K\left(n(x_{i-1}(t) - x_i(t))\right) - K\left(n(x_i(t) - x_{i+1}(t))\right) \right| \right) \quad (3.13)$$

where

$$w_i(t) = v_i(t) - nK(n(x_{i-1}(t) - x_i(t))) + nK(n(x_i(t) - x_{i+1}(t))), \text{ for } i = 1,\ldots,n-1. \quad (3.14)$$

We consider the parameterized family of finite-dimensional systems (3.1), (3.2), (3.3) defined on $\Omega_n$ parameterized by $n \geq 2$ with initial conditions $(x_1(0),\ldots,x_{n-1}(0), v_1(0),\ldots,v_{n-1}(0))$ in $\Omega_n$, for which there exist constants $\bar{E}, \bar{W}, \bar{Z}, \bar{A} > 0$ such that the following conditions hold:

$$\frac{m}{2n} \sum_{i=1}^{n} v_i^2(0) + \frac{m}{n} \sum_{i=1}^{n} \Phi\left(n\left(x_{i-1}(0) - x_i(0)\right)\right) \leq \bar{E},$$

$$\frac{m}{2n} \sum_{i=1}^{n-1} w_i^2(0) + \frac{m}{n} \sum_{i=1}^{n} \Phi\left(n\left(x_{i-1}(0) - x_i(0)\right)\right) \leq \bar{W},$$

$$\frac{1}{2} \sum_{i=1}^{n} \frac{(v_{i-1}(0) - v_i(0))^2}{x_{i-1}(0) - x_i(0)} \leq \bar{Z},$$

$$\max_{i=1,\ldots,n-1} \left( n \left| K\left(n(x_{i-1}(0) - x_i(0))\right) - K\left(n(x_i(0) - x_{i+1}(0))\right) \right| \right) \leq \bar{A}. \quad (3.15)$$

Define the increasing function $F:(0,+\infty) \to \mathbb{R}$ by means of the formula:

$$F(\rho) = \begin{cases} \max\left(\sqrt{\frac{F_1(\rho)}{2\sqrt{2}}}, \frac{F_2(\rho)}{2\sqrt{2L}}, \frac{k(\rho)}{\sqrt{2m}}\right) & \text{if } \rho \geq \rho^* \\ \min\left(-\sqrt{-\frac{F_1(\rho)}{2\sqrt{2}}}, \frac{F_2(\rho)}{2\sqrt{2L}}, \frac{k(\rho)}{\sqrt{2m}}\right) & \text{if } 0 < \rho < \rho^* \end{cases} \quad (3.16)$$

where

$$F_1(\rho) = \int_{\rho^*}^{\rho} s^{-3/2} \mu(s) \sqrt{Q(s)} ds \quad (3.17)$$

$$F_2(\rho) = \int_{\rho^*}^{\rho} s^{-3/2} \mu(s) ds.  \tag{3.18}$$

The following theorem provides results for the parameterized family of finite-dimensional systems (3.1), (3.2), (3.3).

**Theorem 3.1:** *Suppose that Assumption (A) holds. Let $\bar{E}, \bar{W}, \bar{Z}, \bar{A} > 0$ be given constants with*

$$\sqrt{\bar{W}} + \sqrt{\bar{E}} < \min\left(\lim_{\rho \to +\infty}(F(\rho)), -\lim_{\rho \to 0^+}(F(\rho))\right). \tag{3.19}$$

*Then there exist constants $\kappa, R_i > 0$, $i = 1,...,6$, (independent of $n \geq 2$ but dependent on $\bar{E}, \bar{W}, \bar{Z}, \bar{A} > 0$) such that every solution of (3.1), (3.2), (3.3) with initial condition that satisfies (3.15) also satisfies the following estimates for all $n \geq 2$ and $t \geq 0$:*

$$E_n(t) \leq \bar{E} \tag{3.20}$$

$$W_n(t) \leq \bar{W} \tag{3.21}$$

$$n \sum_{i=1}^{n-1} \left( K(n(x_{i-1}(t) - x_i(t))) - K(n(x_i(t) - x_{i+1}(t))) \right)^2 \leq R_1 \tag{3.22}$$

$$a \leq n(x_{i-1}(t) - x_i(t)) \leq b, \text{ for } i = 1,...,n \tag{3.23}$$

$$Z_n(t) \leq R_2 \tag{3.24}$$

$$H_n(t) \leq R_3 \exp(\kappa t) \tag{3.25}$$

$$\max_{i=0,...,n}\left(|v_i(t)|\right) \leq R_4 \tag{3.26}$$

$$n \int_0^t \sum_{i=1}^{n-1} \left( \frac{v_i(s) - v_{i+1}(s)}{x_i(s) - x_{i+1}(s)} - \frac{v_{i-1}(s) - v_i(s)}{x_{i-1}(s) - x_i(s)} \right)^2 ds \leq R_5 + R_6 t \tag{3.27}$$

*where* $a := \dfrac{m}{F^{-1}\left(\sqrt{\dfrac{2}{m}}\left(\sqrt{\bar{E}} + \sqrt{\bar{W}}\right)\right)} \in (0, L]$ *and* $b := \dfrac{m}{F^{-1}\left(-\sqrt{\dfrac{2}{m}}\left(\sqrt{\bar{E}} + \sqrt{\bar{W}}\right)\right)} \geq L$.

**Remarks on Theorem 3.1:**
**(i)** Theorem 3.1 can be interpreted as a stability result for the numerical scheme (3.1), (3.2), (3.3). Indeed, notice that under the assumptions of Theorem 3.1, the solutions of (3.1), (3.2), (3.3) are bounded and the approximations of the solutions of the initial-boundary value problem (2.1)-(2.4) with (2.14) are bounded in various spaces.
**(ii)** When $\lim_{\rho \to +\infty}(F(\rho)) = \lim_{\rho \to 0^+}(-F(\rho)) = +\infty$ then condition (3.19) is automatically satisfied and arbitrary initial conditions for the finite-dimensional system (3.1), (3.2), (3.3) can be allowed. For an ideal gas under constant entropy, it holds that $\lim_{\rho \to +\infty}(F(\rho)) = \lim_{\rho \to 0^+}(-F(\rho)) = +\infty$, and consequently, arbitrary initial conditions for the finite-dimensional system (3.1), (3.2), (3.3) can be allowed. Indeed, for an ideal gas under constant entropy we have $P(\rho) = c\rho^\gamma$, $P = f\rho T$ and $\mu = \beta\sqrt{T}$, where $T > 0$ is the temperature and $c, f, \beta, \gamma > 0$ are constants with $\gamma \in (1,2)$ (see [10, 15]). Therefore, we have $\mu(\rho) = A\rho^{\frac{\gamma-1}{2}}$ with $A = \beta\sqrt{f^{-1}c}$ and

$$Q(\rho) := \frac{c}{\gamma-1}\left(\rho^\gamma - \gamma(\rho^*)^{\gamma-1}\rho + (\gamma-1)(\rho^*)^\gamma\right), \qquad k(\rho) := \frac{2A}{\gamma-1}\left(\rho^{\frac{\gamma-1}{2}} - (\rho^*)^{\frac{\gamma-1}{2}}\right),$$

$F_1(\rho) = A\int_{\rho^*}^{\rho} s^{\frac{\gamma-4}{2}}\sqrt{Q(s)}ds$. Notice that the fact that $Q(\rho) \geq \frac{c}{2}(\rho^*)^\gamma$ for $0 < \rho < \frac{\gamma-1}{2\gamma}\rho^*$ implies that

$$F_1(\rho) \leq -\frac{2A}{2-\gamma}\sqrt{\frac{c}{2}}(\rho^*)^{\frac{\gamma}{2}}\left[\rho^{\frac{\gamma-2}{2}} - \left(\frac{\gamma-1}{2\gamma}\rho^*\right)^{\frac{\gamma-2}{2}}\right] \text{ for } 0 < \rho < \frac{\gamma-1}{2\gamma}\rho^*.$$ Consequently, the fact that $\gamma \in (1,2)$

and definition (3.16) guarantees that $\lim_{\rho \to +\infty}(F(\rho)) = \lim_{\rho \to 0^+}(-F(\rho)) = +\infty$ for an ideal gas under constant entropy.

**(iii)** However, for the case of the viscous Saint-Venant problem (where $\mu(\rho) = \mu\rho$, $k(\rho) = \mu(\rho - \rho^*)$ for a constant $\mu > 0$), condition (3.19) gives

$$\sqrt{\overline{W}} + \sqrt{\overline{E}} < \sqrt{\frac{\mu\rho^*}{3}}\sqrt{\rho^*g}\max\left(1, \sqrt{\frac{3\mu}{2L\sqrt{\rho^*g}}}\right),$$ and therefore Theorem 3.1 does not hold for

arbitrary initial conditions for the finite-dimensional system (3.1), (3.2), (3.3). This result is in accordance with the fact that the viscous Saint-Venant problem is well-posed for initial data close to the equilibrium (see [20, 31]) as well as with our everyday experience that waves can create areas where the fluid is absent (and thus $\rho = 0$).

**(iv)** Working as above, it is possible to show that $\lim_{\rho \to +\infty}(F(\rho)) = \lim_{\rho \to 0^+}(-F(\rho)) = +\infty$ are valid for any gas that satisfies the inequalities $P(\rho) \geq c\rho^\gamma$ and $\mu(\rho) \geq A\rho^\eta$, where $c, A, \gamma > 0$, $\eta \in \Re$ are constants with $\gamma \in (1,2)$ and $\eta \in [0, 1/2]$. Therefore, condition $\lim_{\rho \to +\infty}(F(\rho)) = \lim_{\rho \to 0^+}(-F(\rho)) = +\infty$ is valid for a wide class of gases.

**(v)** The function $F_1(\rho)$ defined by (3.17) was recently used in [15, 16, 17, 18] for control problems related to model (2.1)-(2.4) in order to obtain estimates of the upper and lower bounds of density or liquid level.

3.II. The choice of the initial condition

Given $n \geq 2$, $v \in H_0^1((0,L))$, $\rho \in W^{1,\infty}((0,L)) \cap C^0([0,L];(0,+\infty))$ with $\int_0^L \rho(x)dx = m$, we can find a unique set of points $x_n = 0 < x_{n-1} < ... < x_1 < x_0 = L$ with $\int_{x_i}^{x_{i-1}} \rho(x)dx = \frac{m}{n}$ for $i = 1,...,n$. We also define $v_i = v(x_i)$ for $i = 0, 1, ..., n$. In this way, we may define the mapping

$$X \ni (\rho, v) \to P_n(\rho, v) = (x_1, ..., x_{n-1}, v_1, ..., v_{n-1}) \in \Omega_n \tag{3.28}$$

with
$$X = S \times H_0^1((0,L)) \text{ and}$$
$$S = \left\{ \rho \in W^{1,\infty}((0,L)) \cap C^0([0,L]; (0,+\infty)) : \int_0^L \rho(x)dx = m \right\}. \quad (3.29)$$

For the numerical scheme (3.1), (3.2), (3.3) and the initial-boundary value problem (2.1)-(2.4) with (2.14) we consider as initial condition the point $(x_1(0),...,x_{n-1}(0), v_1(0),...,v_{n-1}(0)) = P_n(\rho_0, v_0) \in \Omega_n$. For this particular selection, we have the following result.

**Theorem 3.2:** *Let an initial condition $(\rho_0, v_0) \in X$ for the initial-boundary value problem (2.1)-(2.4), (2.14) be given and suppose that*

$$\sqrt{mL\|v_0'\|_2^2 + 2m^5 \bar{M}^2 \frac{\|\rho_0'\|_\infty^2}{\rho_{\min}^6} + m\Phi\left(\frac{m}{\|\rho_0\|_\infty}\right)} + \sqrt{\frac{mL}{2}\|v_0'\|_2^2 + m\Phi\left(\frac{m}{\|\rho_0\|_\infty}\right)} \quad (3.30)$$
$$< \min\left( \lim_{\rho \to +\infty}(F(\rho)), -\lim_{\rho \to 0^+}(F(\rho)) \right)$$

*where $\bar{M} = \max\left\{ K'(s) : \frac{m}{\|\rho_0\|_\infty} \leq s \leq \frac{m}{\rho_{\min}} \right\}$ and $\rho_{\min} = \min_{0 \leq x \leq L}(\rho_0(x))$. Then there exist constants $\bar{E}, \bar{W}, \bar{Z}, \bar{A} > 0$ such that conditions (3.15), (3.19) are valid for all $n \geq 2$ with $(x_1(0),...,x_{n-1}(0), v_1(0),...,v_{n-1}(0)) = P_n(\rho_0, v_0) \in \Omega_n$.*

Again, it should be noted that if $\lim_{\rho \to +\infty}(F(\rho)) = \lim_{\rho \to 0^+}(-F(\rho)) = +\infty$ then condition (3.30) is automatically satisfied.

3.III. Properties of the numerical scheme

Let a solution of (3.1), (3.2), (3.3) defined for all $t \geq 0$ be given. Define for all $t \geq 0$:

$$\rho_i(t) := \frac{m}{n(x_{i-1}(t) - x_i(t))}, \text{ for } i = 1,...,n \quad (3.31)$$

$$\rho_0(t) = \rho_1(t) = \frac{m}{n(L - x_1(t))} \quad (3.32)$$

$$v^{(n)}(t,x) = v_i(t) + \frac{v_{i-1}(t) - v_i(t)}{x_{i-1}(t) - x_i(t)}(x - x_i(t)) \quad , \quad x \in [x_i(t), x_{i-1}(t))$$

$$\rho^{(n)}(t,x) = \rho_i(t) + \frac{\rho_{i-1}(t) - \rho_i(t)}{x_{i-1}(t) - x_i(t)}(x - x_i(t)) \quad , \quad x \in [x_i(t), x_{i-1}(t))$$
$$\text{for } i = 1,...,n \quad (3.33)$$

$$v^{(n)}(t,L) = 0 \quad , \quad \rho^{(n)}(t,L) = \rho_0(t). \quad (3.34)$$

The functions $v^{(n)}, \rho^{(n)} : \mathbb{R}_+ \times [0,L] \to \mathbb{R}$ are $C^1$ on the set $\bigcup_{t \geq 0} \{(t,x) : x \in [0,L], x \neq x_i(t), i = 0,...,n\}$ (are $C^1$ a.e. on $\mathbb{R}_+ \times [0,L]$) and are used as approximations of the solution of the initial-boundary value problem (2.1)-(2.4) with (2.14). The functions $\dot{\rho}^{(n)}, \dot{v}^{(n)}$ defined for each $n \geq 2$ and $t \geq 0$ by means of the formulas

$$\dot{\rho}^{(n)}(t,x) = \dot{\rho}_i(t) + \frac{\dot{\rho}_{i-1}(t) - \dot{\rho}_i(t)}{x_{i-1}(t) - x_i(t)}(x - x_i(t)) - \frac{\rho_{i-1}(t) - \rho_i(t)}{x_{i-1}(t) - x_i(t)} v^{(n)}(t,x)$$

for $x \in [x_i(t), x_{i-1}(t))$, $i = 1,...,n$ (3.35)

$$\dot{v}^{(n)}(t,x) = \dot{v}_i(t) + \frac{\dot{v}_{i-1}(t) - \dot{v}_i(t)}{x_{i-1}(t) - x_i(t)}(x - x_i(t))$$

$$- \frac{(v_{i-1}(t) - v_i(t))^2}{(x_{i-1}(t) - x_i(t))^2}(x - x_i(t)) - \frac{v_{i-1}(t) - v_i(t)}{x_{i-1}(t) - x_i(t)} v_i(t)$$

for $x \in [x_i(t), x_{i-1}(t))$, $i = 1,...,n$ (3.36)

are the weak time derivatives of the functions $\rho^{(n)}, v^{(n)}$ defined by (3.33).

The following theorem summarizes our convergence results.

**Theorem 3.3:** *Suppose that Assumption (A) holds. Let $(\rho_0, v_0) \in X$ that satisfies (3.30) be given. Consider the unique solution of (3.1), (3.2), (3.3) with $n \geq 2$ and initial condition $(x_1(0),...,x_{n-1}(0), v_1(0),...,v_{n-1}(0)) = P_n(\rho_0, v_0) \in \Omega_n$. Define $v^{(n)}, \rho^{(n)} : \mathbb{R}_+ \times [0,L] \to \mathbb{R}$ by means of (3.31)-(3.34). Then there exists a subsequence of $\{(\rho^{(n)}, v^{(n)}) : n \geq 2\}$ that converges to a weak solution $(\rho, v)$ of the initial-boundary value problem (2.1)-(2.4) with (2.14) in the following sense:*

$$\rho^{(n)} \to \rho \text{ and } v^{(n)} \to v \text{ in } C^0\left([0,T]; L^2((0,L))\right) \text{ strongly}$$

$$v^{(n)} \to v \text{ in } L^\infty\left((0,T); H_0^1((0,L))\right) \text{ weak star}$$

$$\rho^{(n)} \to \rho \text{ in } L^\infty\left((0,T); H^1((0,L))\right) \text{ weak star}$$

$$\dot{\rho}^{(n)} \to \rho_t \text{ in } L^\infty\left((0,T); L^2((0,L))\right) \text{ weak star}$$

$$\dot{v}^{(n)} \to v_t \text{ in } L^2\left((0,T); H^{-1}((0,L))\right) \text{ weakly}$$

$$\rho_x^{(n)} \to \rho_x \text{ in } L^\infty\left((0,T) \times (0,L)\right) \text{ weak star.}$$

*Moreover, if in addition*

$$\rho_x^{(n)} \to \rho_x \text{ in } C^0\left([0,T]; L^2((0,L))\right) \text{ strongly}$$

*then the weak solution satisfies the following additional estimate:*

$$W(\rho[t], v[t]) \leq W(\rho[t_0], v[t_0]), \text{ for all } T \geq t \geq t_0 \geq 0. \quad (3.37)$$

Theorem 3.3 is a convergence result for the particle scheme (3.1), (3.2), (3.3). Indeed, a subsequence of the numerical solutions generated by the particle scheme (3.1), (3.2), (3.3) is guaranteed to converge to a weak solution of the initial-boundary value problem (2.1)-(2.4) with (2.14). However, Theorem 3.3 can also be considered to be an existence result for the initial-boundary value problem (2.1)-(2.4) with (2.14) because it guarantees that there exists a weak solution of the initial-boundary value problem (2.1)-(2.4) with (2.14). Indeed, we have the following result.

**Corollary 3.4:** *Suppose that Assumption (A) holds. Let $(\rho_0, v_0) \in X$ that satisfies (3.30) be given. Then there exists weak solution $(\rho, v)$ of the initial-boundary value problem (2.1)-(2.4) with (2.14)*

Again, it should be noted that if $\lim_{\rho \to +\infty}(F(\rho)) = \lim_{\rho \to 0^+}(-F(\rho)) = +\infty$ then condition (3.30) is automatically satisfied. In this case, Theorem 3.3 and Corollary 3.4 provide global existence results (with no restriction on the size of the initial conditions).

## 4. Proofs of Main Results

The following facts hold for the parameterized family of finite-dimensional systems (3.1), (3.2), (3.3) defined on $\Omega_n$ (recall (3.4)).

<u>Fact 1:</u> *For every solution of (3.1), (3.2), (3.3) the following equation holds*

$$\dot{E}_n(t) = -nm \sum_{i=1}^{n} K'(n(x_{i-1}(t) - x_i(t)))(v_{i-1}(t) - v_i(t))^2 \leq 0 \qquad (4.1)$$

*as long as the solution of (3.1), (3.2), (3.3) is defined.*

**Remark:** Inequality (4.1) is the discretized version of inequality (2.12) and shows that the discretized mechanical energy $E_n$ is decaying.

**Proof:** Using (3.1), (3.2), (3.3) and (3.10), we get:

$$\dot{E}_n(t) = m \sum_{i=1}^{n-1} v_i(t) \left( \Phi'(n(x_{i-1}(t) - x_i(t))) - \Phi'(n(x_i(t) - x_{i+1}(t))) \right)$$

$$+ m \sum_{i=1}^{n} \Phi'\left(n(x_{i-1}(t) - x_i(t))\right)(v_{i-1}(t) - v_i(t))$$

$$+ mn \sum_{i=1}^{n-1} v_i(t) \left( K'(n(x_{i-1}(t) - x_i(t)))(v_{i-1}(t) - v_i(t)) + K'(n(x_i(t) - x_{i+1}(t)))(v_{i+1}(t) - v_i(t)) \right)$$

$$= mn \sum_{i=1}^{n-1} v_i(t) \left( K'(n(x_{i-1}(t) - x_i(t)))(v_{i-1}(t) - v_i(t)) + K'(n(x_i(t) - x_{i+1}(t)))(v_{i+1}(t) - v_i(t)) \right)$$

$$= -mn \sum_{i=1}^{n} K'(n(x_{i-1}(t) - x_i(t)))(v_{i-1}(t) - v_i(t))^2$$

The fact that $\sum_{i=1}^{n} K'(n(x_{i-1}(t) - x_i(t)))(v_{i-1}(t) - v_i(t))^2 \geq 0$ follows from definitions (3.5) (which implies that $K'(x) > 0$, for all $x > 0$). The proof is complete. ◁

**Fact 2:** *Suppose that Assumption (A) is valid. Then every solution of (3.1), (3.2), (3.3) is defined for all $t \geq 0$.*

**Proof:** This fact is a direct consequence of Assumption (A), (3.9) and Fact 1, i.e., the fact that

$$\frac{m}{2n}\sum_{i=1}^{n} v_i^2(t) + \frac{m}{n}\sum_{i=1}^{n} \Phi\big(n(x_{i-1}(t) - x_i(t))\big)$$

$$\leq E_n(0) = \frac{m}{2n}\sum_{i=1}^{n} v_i^2(0) + \frac{m}{n}\sum_{i=1}^{n} \Phi\big(n(x_{i-1}(0) - x_i(0))\big)$$

as long as the solution of (3.1), (3.2), (3.3) is defined, and the facts that $\lim_{x \to 0^+}(\Phi(x)) = +\infty$, $\inf_{x>0}(\Phi(x)) > -\infty$ which shows that the solution of (3.1), (3.2), (3.3) remains bounded and cannot approach the boundary of the open set $\Omega_n = \left\{ (x_1,...,x_{n-1}, v_1,...,v_{n-1}) \in (0,L)^{n-1} \times \mathbb{R}^{n-1} : x_{n-1} < ... < x_1 \right\}$. The proof is complete. ◁

**Fact 3:** *Suppose that Assumption (A) is valid. Then for every solution of (3.1), (3.2), (3.3) the following equation holds for all $t \geq 0$:*

$$\dot{W}_n(t) =$$
$$-mn\sum_{i=1}^{n-1} \big(\Phi'(n(x_{i-1}(t)-x_i(t))) - \Phi'(n(x_i(t)-x_{i+1}(t)))\big)\big(K(n(x_{i-1}(t)-x_i(t))) - K(n(x_i(t)-x_{i+1}(t)))\big) \leq 0$$

(4.2)

**Remark:** Inequality (4.2) is the discretized version of inequality (2.13) and shows that the discretized mechanical energy $W_n$ is decaying.

**Proof:** Using (3.1), (3.2), (3.3) and definition (3.14) we get the following equations for $i = 1,...,n-1$:

$$\dot{x}_i(t) = w_i(t) + nK(n(x_{i-1}(t) - x_i(t))) - nK(n(x_i(t) - x_{i+1}(t)))$$
$$\dot{w}_i(t) = n\Phi'(n(x_{i-1}(t) - x_i(t))) - n\Phi'(n(x_i(t) - x_{i+1}(t)))$$

(4.3)

Using definition (3.11) and equations (4.3) we get:

$$\dot{W}_n(t) = m\sum_{i=1}^{n-1} w_i \big(\Phi'(n(x_{i-1}(t)-x_i(t))) - \Phi'(n(x_i(t)-x_{i+1}(t)))\big)$$

$$+m\sum_{i=2}^{n-1} \Phi'\big(n(x_{i-1}(t)-x_i(t))\big)\big(w_{i-1}(t) - w_i(t)\big)$$

$$+mn\sum_{i=2}^{n-1} \Phi'\big(n(x_{i-1}(t)-x_i(t))\big)\big(K(n(x_{i-2}(t)-x_{i-1}(t))) - 2K(n(x_{i-1}(t)-x_i(t))) + K(n(x_i(t)-x_{i+1}(t)))\big)$$

$$-m\Phi'\big(n(L-x_1(t))\big)\big(w_1(t) + nK(n(L-x_1(t))) - nK(n(x_1(t)-x_2(t)))\big)$$

$$+m\Phi'\big(nx_{n-1}(t)\big)\big(w_{n-1}(t) + nK(n(x_{n-2}(t)-x_{n-1}(t))) - nK(nx_{n-1}(t))\big)$$

$$= -mn\sum_{i=1}^{n-1} \big(\Phi'\big(n(x_{i-1}(t)-x_i(t))\big) - \Phi'\big(n(x_i(t)-x_{i+1}(t))\big)\big)\big(K(n(x_{i-1}(t)-x_i(t))) - K(n(x_i(t)-x_{i+1}(t)))\big)$$

The fact that

$$\sum_{i=1}^{n-1}\left(\Phi'(n(x_{i-1}(t)-x_i(t)))-\Phi'(n(x_i(t)-x_{i+1}(t)))\right)\left(K(n(x_{i-1}(t)-x_i(t)))-K(n(x_i(t)-x_{i+1}(t)))\right)\geq 0$$

follows from definitions (3.5) (which implies that both $K(x)$ and $\Phi'(x)$ are increasing functions). The proof is complete. ◁

The following fact is trivial and is proof is omitted.

<u>Fact 4:</u> *Let the set $\{u_i \in \mathbb{R} : i=1,...,n\}$ be given. Then the following inequality holds:*

$$\max_{i=1,...,n}(u_i) - \min_{i=1,...,n}(u_i) \leq \sum_{i=1}^{n-1}|u_{i+1}-u_i| \qquad (4.4)$$

*Moreover, if $\max_{i=1,...,n}(u_i) \geq 0 \geq \min_{i=1,...,n}(u_i)$, then the following inequality holds:*

$$\max_{i=1,...,n}(|u_i|) \leq \sum_{i=1}^{n-1}|u_{i+1}-u_i| \qquad (4.5)$$

The following fact plays a crucial role in what follows.

<u>Fact 5:</u> *The following inequality holds for all $t \geq 0$:*

$$n\sum_{i=1}^{n-1}\left(K(n(x_{i-1}(t)-x_i(t)))-K(n(x_i(t)-x_{i+1}(t)))\right)^2 \leq \frac{2}{m}\left(\sqrt{W_n(t)}+\sqrt{E_n(t)}\right)^2 \qquad (4.6)$$

*Define the increasing function $F:(0,+\infty) \to \mathbb{R}$ by means of (3.16). If the following inequality holds*

$$\sqrt{W_n(t)}+\sqrt{E_n(t)} < \min\left(\lim_{\rho\to+\infty}(F(\rho)), -\lim_{\rho\to 0^+}(F(\rho))\right) \qquad (4.7)$$

*then the following inequalities hold for all $i=1,...,n$:*

$$0 < F^{-1}\left(-\left(\sqrt{W_n(t)}+\sqrt{E_n(t)}\right)\right) \leq \frac{m}{n(x_{i-1}(t)-x_i(t))} \leq F^{-1}\left(\sqrt{W_n(t)}+\sqrt{E_n(t)}\right) < +\infty \qquad (4.8)$$

**Proof:** Since $\sum_{i=1}^{n}(x_{i-1}(t)-x_i(t))\left(\frac{m}{n(x_{i-1}(t)-x_i(t))}-\rho^*\right)=0$ (a consequence of the fact that $\rho^* = m/L$), we get

$$\min_{i=1,...,n}\left(\frac{m}{n(x_{i-1}(t)-x_i(t))}\right) \leq \rho^* \leq \max_{i=1,...,n}\left(\frac{m}{n(x_{i-1}(t)-x_i(t))}\right) \qquad (4.9)$$

Using (4.9) and the facts that $F_1, F_2, k$ are increasing (recall (2.7), (3.17), (3.18)), we get:

$$\min_{i=1,...,n}\left(F_j\left(\frac{m}{n(x_{i-1}(t)-x_i(t))}\right)\right) \leq 0 \leq \max_{i=1,...,n}\left(F_j\left(\frac{m}{n(x_{i-1}(t)-x_i(t))}\right)\right), \text{ for } j=1,2 \qquad (4.10)$$

$$\min_{i=1,\ldots,n}\left(k\left(\frac{m}{n(x_{i-1}(t)-x_i(t))}\right)\right) \leq 0 \leq \max_{i=1,\ldots,n}\left(k\left(\frac{m}{n(x_{i-1}(t)-x_i(t))}\right)\right) \quad (4.11)$$

Fact 4 in conjunction with (4.10) and (4.11) implies the following inequalities

$$\max_{i=1,\ldots,n}\left(\left|F_1\left(\frac{m}{n(x_{i-1}(t)-x_i(t))}\right)\right|\right) \leq \sum_{i=1}^{n-1}\left|F_1\left(\frac{m}{n(x_i(t)-x_{i+1}(t))}\right) - F_1\left(\frac{m}{n(x_{i-1}(t)-x_i(t))}\right)\right|$$

$$\max_{i=1,\ldots,n}\left(\left|F_2\left(\frac{m}{n(x_{i-1}(t)-x_i(t))}\right)\right|\right) \leq \sum_{i=1}^{n-1}\left|F_2\left(\frac{m}{n(x_i(t)-x_{i+1}(t))}\right) - F_2\left(\frac{m}{n(x_{i-1}(t)-x_i(t))}\right)\right| \quad (4.12)$$

$$\max_{i=1,\ldots,n}\left(\left|k\left(\frac{m}{n(x_{i-1}(t)-x_i(t))}\right)\right|\right) \leq \sum_{i=1}^{n-1}\left|k\left(\frac{m}{n(x_i(t)-x_{i+1}(t))}\right) - k\left(\frac{m}{n(x_{i-1}(t)-x_i(t))}\right)\right|$$

Define the following functions for all $u \in K\left((0,+\infty)\right)$:

$$\tilde{F}_j(u) = F_j\left(\frac{m}{K^{-1}(u)}\right), \text{ for } j=1,2 \quad (4.13)$$

Notice that definitions (4.13) imply that

$$F_j\left(\frac{m}{x}\right) = \tilde{F}_j(K(x)), \text{ for } j=1,2 \text{ and } x>0 \quad (4.14)$$

Using (4.12), (4.14) and the fact that $K(x) = -\frac{1}{m}k\left(\frac{m}{x}\right)$ for $x>0$ (recall (3.5)), we get:

$$\max_{i=1,\ldots,n}\left(\left|F_1\left(\frac{m}{n(x_{i-1}(t)-x_i(t))}\right)\right|\right) \leq \sum_{i=1}^{n-1}\left|\tilde{F}_1\left(K\left(n(x_i(t)-x_{i+1}(t))\right)\right) - \tilde{F}_1\left(K\left(n(x_{i-1}(t)-x_i(t))\right)\right)\right|$$

$$\max_{i=1,\ldots,n}\left(\left|F_2\left(\frac{m}{n(x_{i-1}(t)-x_i(t))}\right)\right|\right) \leq \sum_{i=1}^{n-1}\left|\tilde{F}_2\left(K\left(n(x_i(t)-x_{i+1}(t))\right)\right) - \tilde{F}_2\left(K\left(n(x_{i-1}(t)-x_i(t))\right)\right)\right| \quad (4.15)$$

$$\max_{i=1,\ldots,n}\left(\left|k\left(\frac{m}{n(x_{i-1}(t)-x_i(t))}\right)\right|\right) \leq m\sum_{i=1}^{n-1}\left|K\left(n(x_i(t)-x_{i+1}(t))\right) - K\left(n(x_{i-1}(t)-x_i(t))\right)\right|$$

Definitions (3.17), (3.18), (4.13) and (3.5) imply that

$$\tilde{F}_1'(u) = -\sqrt{mK^{-1}(u)Q\left(\frac{m}{K^{-1}(u)}\right)}, \text{ for all } u \in K\left((0,+\infty)\right) \quad (4.16)$$

$$\tilde{F}_2'(u) = -\sqrt{mK^{-1}(u)}, \text{ for all } u \in K\left((0,+\infty)\right) \quad (4.17)$$

Equations (4.16) and (4.17) give us the following estimates for all $i=1,\ldots,n-1$:

$$\left|\tilde{F}_1(y_{i+1}) - \tilde{F}_1(y_i)\right|$$
$$\leq \sqrt{m}\left|y_{i+1} - y_i\right| \max\left\{\left|\sqrt{K^{-1}(u)Q\left(\frac{m}{K^{-1}(u)}\right)}\right| : \min(y_{i+1}, y_i) \leq u \leq \max(y_{i+1}, y_i)\right\} \quad (4.18)$$

$$\left|\tilde{F}_2(y_{i+1}) - \tilde{F}_2(y_i)\right|$$
$$\leq \sqrt{m}\left|y_{i+1} - y_i\right| \max\left\{\left|\sqrt{K^{-1}(u)}\right| : \min(y_{i+1}, y_i) \leq u \leq \max(y_{i+1}, y_i)\right\} \quad (4.19)$$

with $y_i = K(n(x_{i-1}(t) - x_i(t)))$. Since $K^{-1}$ is increasing (recall (3.5)) we get from (4.19):

$$\left|\tilde{F}_2(y_{i+1}) - \tilde{F}_2(y_i)\right| \leq \sqrt{nm \max(x_i(t) - x_{i+1}(t), x_{i-1}(t) - x_i(t))}\left|y_{i+1} - y_i\right| \quad (4.20)$$

The function $g(s) = sQ\left(\dfrac{m}{s}\right)$ for $s > 0$ satisfies $g'(s) = P(\rho^*) - P(m/s)$ (recall (2.8)). Since $P$ is increasing, it follows that $g$ is convex. Consequently, we get from (4.18):

$$\left|\tilde{F}_1(y_{i+1}) - \tilde{F}_1(y_i)\right|$$
$$\leq \sqrt{m}\left|y_{i+1} - y_i\right|\sqrt{\max\left\{g(s) : K^{-1}(\min(y_{i+1}, y_i)) \leq s \leq K^{-1}(\max(y_{i+1}, y_i))\right\}}$$
$$= \sqrt{m}\left|y_{i+1} - y_i\right|\sqrt{\max\left(g(K^{-1}(\min(y_{i+1}, y_i))), g(K^{-1}(\max(y_{i+1}, y_i)))\right)} \quad (4.21)$$
$$\leq \sqrt{m}\left|y_{i+1} - y_i\right|\sqrt{g(K^{-1}(\min(y_{i+1}, y_i))) + g(K^{-1}(\max(y_{i+1}, y_i)))}$$
$$= \sqrt{nm}\left|y_{i+1} - y_i\right|\sqrt{(x_i(t) - x_{i+1}(t))Q\left(\frac{m}{n(x_i(t) - x_{i+1}(t))}\right) + (x_{i-1}(t) - x_i(t))Q\left(\frac{m}{n(x_{i-1}(t) - x_i(t))}\right)}$$

Using (4.15), (4.20) and (4.21), we get:

$$\max_{i=1,\ldots,n}\left(\left|F_1\left(\frac{m}{n(x_{i-1}(t) - x_i(t))}\right)\right|\right)$$
$$\leq \sqrt{nm}\sum_{i=1}^{n-1}\left|K(n(x_i(t) - x_{i+1}(t))) - K(n(x_{i-1}(t) - x_i(t)))\right|\sqrt{(x_i(t) - x_{i+1}(t))Q\left(\frac{m}{n(x_i(t) - x_{i+1}(t))}\right)} \quad (4.22)$$
$$+ \sqrt{nm}\sum_{i=1}^{n-1}\left|K(n(x_i(t) - x_{i+1}(t))) - K(n(x_{i-1}(t) - x_i(t)))\right|\sqrt{(x_{i-1}(t) - x_i(t))Q\left(\frac{m}{n(x_{i-1}(t) - x_i(t))}\right)}$$

$$\max_{i=1,\ldots,n}\left(\left|F_2\left(\frac{m}{n(x_{i-1}(t) - x_i(t))}\right)\right|\right)$$
$$\leq \sqrt{nm}\sum_{i=1}^{n-1}\left|K(n(x_i(t) - x_{i+1}(t))) - K(n(x_{i-1}(t) - x_i(t)))\right|\sqrt{x_i(t) - x_{i+1}(t)} \quad (4.23)$$
$$+ \sqrt{nm}\sum_{i=1}^{n-1}\left|K(n(x_i(t) - x_{i+1}(t))) - K(n(x_{i-1}(t) - x_i(t)))\right|\sqrt{x_{i-1}(t) - x_i(t)}$$

Applying the Cauchy-Schwarz inequality, we obtain from (4.15), (4.22) and (4.23) and the fact that $\sum_{i=1}^{n}(x_{i-1}(t)-x_i(t))=L$:

$$\max_{i=1,\ldots,n}\left(\left|k\left(\frac{m}{n(x_{i-1}(t)-x_i(t))}\right)\right|\right) \leq m\left(n\sum_{i=1}^{n-1}\left|K\left(n(x_i(t)-x_{i+1}(t))\right)-K\left(n(x_{i-1}(t)-x_i(t))\right)\right|^2\right)^{1/2} \quad (4.24)$$

$$\max_{i=1,\ldots,n}\left(\left|F_1\left(\frac{m}{n(x_{i-1}(t)-x_i(t))}\right)\right|\right)$$
$$\leq 2\sqrt{m}\left(n\sum_{i=1}^{n-1}\left|K\left(n(x_i(t)-x_{i+1}(t))\right)-K\left(n(x_{i-1}(t)-x_i(t))\right)\right|^2\right)^{1/2}\left(\sum_{i=1}^{n}(x_{i-1}(t)-x_i(t))Q\left(\frac{m}{n(x_{i-1}(t)-x_i(t))}\right)\right)^{1/2}$$
(4.25)

$$\max_{i=1,\ldots,n}\left(\left|F_2\left(\frac{m}{n(x_{i-1}(t)-x_i(t))}\right)\right|\right) \leq 2\sqrt{Lm}\left(n\sum_{i=1}^{n-1}\left|K\left(n(x_i(t)-x_{i+1}(t))\right)-K\left(n(x_{i-1}(t)-x_i(t))\right)\right|^2\right)^{1/2} \quad (4.26)$$

Using (3.8), (3.10), (3.11) and the facts that $\rho^* = m/L$, $\sum_{i=1}^{n}(x_{i-1}(t)-x_i(t))=L$, we get:

$$\sum_{i=1}^{n}(x_{i-1}(t)-x_i(t))Q\left(\frac{m}{n(x_{i-1}(t)-x_i(t))}\right) = \frac{m}{n}\sum_{i=1}^{n}\Phi\left(n(x_{i-1}(t)-x_i(t))\right) \leq \min\left(E_n(t),W_n(t)\right) \quad (4.27)$$

Using the inequality

$$\frac{\delta+1}{\delta}\left(v_i(t)-nK(n(x_{i-1}(t)-x_i(t)))+nK(n(x_i(t)-x_{i+1}(t)))\right)^2 + (\delta+1)v_i^2(t)$$
$$\geq n^2\left(K(n(x_{i-1}(t)-x_i(t)))-K(n(x_i(t)-x_{i+1}(t)))\right)^2$$

that holds for all $\delta > 0$, we obtain the inequality

$$\frac{2(\delta+1)}{\delta m}\left(\frac{m}{2n}\sum_{i=1}^{n-1}\left(v_i(t)-nK(n(x_{i-1}(t)-x_i(t)))+nK(n(x_i(t)-x_{i+1}(t)))\right)^2 + \frac{\delta m}{2n}\sum_{i=1}^{n}v_i^2(t)\right)$$
$$\geq n\sum_{i=1}^{n-1}\left(K(n(x_{i-1}(t)-x_i(t)))-K(n(x_i(t)-x_{i+1}(t)))\right)^2$$
(4.28)

Using definitions (3.10), (3.11), (3.14) and the fact that $\sum_{i=1}^{n}\Phi\left(n(x_{i-1}(t)-x_i(t))\right)\geq 0$ (a consequence of the facts that $\rho^* = m/L$, $\sum_{i=1}^{n}(x_{i-1}(t)-x_i(t))=L$, $Q(\rho)\geq 0$ for all $\rho > 0$ and (3.8) which implies that $\Phi(x) = \frac{x}{m}Q(m/x) + P(\rho^*)\left(\frac{1}{\rho^*}-\frac{x}{m}\right)$ for all $x > 0$), we get from (4.28):

$$\frac{2}{m}\left(\frac{(\delta+1)}{\delta}W_n(t)+(\delta+1)E_n(t)\right) \geq n\sum_{i=1}^{n-1}\left(K(n(x_{i-1}(t)-x_i(t)))-K(n(x_i(t)-x_{i+1}(t)))\right)^2 \quad (4.29)$$

Selecting $\delta = \sqrt{W_n(t)/E_n(t)}$ when $E_n(t) > 0$ we get (4.6). Estimate (4.29) shows that inequality (4.6) also holds for the case where $E_n(t) = 0$.

Combining (4.6), (4.24), (4.25), (4.26) and (4.27) we get the following estimates for all $i = 1,...,n$:

$$-\sqrt{2m}\left(\sqrt{W_n(t)} + \sqrt{E_n(t)}\right) \leq k\left(\frac{m}{n(x_{i-1}(t) - x_i(t))}\right) \leq \sqrt{2m}\left(\sqrt{W_n(t)} + \sqrt{E_n(t)}\right)$$

$$-2\sqrt{2L}\left(\sqrt{W_n(t)} + \sqrt{E_n(t)}\right) \leq F_2\left(\frac{m}{n(x_{i-1}(t) - x_i(t))}\right) \leq 2\sqrt{2L}\left(\sqrt{W_n(t)} + \sqrt{E_n(t)}\right) \quad (4.30)$$

$$-2\sqrt{2}\left(\sqrt{W_n(t)} + \sqrt{E_n(t)}\right)^2 \leq F_1\left(\frac{m}{n(x_{i-1}(t) - x_i(t))}\right) \leq 2\sqrt{2}\left(\sqrt{W_n(t)} + \sqrt{E_n(t)}\right)^2$$

Combining (4.30), (3.16), (3.17), (3.18) and using (4.9), we obtain estimates (4.8). The proof is complete. ◁

We next show the following lemma.

**Lemma 4.1:** *For every pair of constants $a, b > 0$ with $0 < a \leq L \leq b$ there exist constants $\sigma, \omega, r > 0$ (independent of $n \geq 2$) and a function $B \in K_\infty$ (independent of $n \geq 2$) for which the following property holds: If*

$$a \leq n(x_{i-1}(t) - x_i(t)) \leq b, \text{ for all } i = 1,...,n \quad (4.31)$$

*then*

$$\dot{Z}_n(t) \leq -\sigma Z_n(t) + \omega Z_n^2(t) + B(W_n(t) + E_n(t)) \quad (4.32)$$

$$\dot{Z}_n(t) \leq -\frac{rn}{6}\sum_{i=1}^{n-1}\left(\frac{v_i(t) - v_{i+1}(t)}{x_i(t) - x_{i+1}(t)} - \frac{v_{i-1}(t) - v_i(t)}{x_{i-1}(t) - x_i(t)}\right)^2 + \omega Z_n^2(t) + B(W_n(t) + E_n(t)) \quad (4.33)$$

**Proof:** Equations (3.1), (3.2), (3.3) and definition (3.12) give:

$$\dot{Z}_n(t) = -\frac{1}{2}\sum_{i=1}^{n}\frac{(v_{i-1}(t) - v_i(t))^3}{(x_{i-1}(t) - x_i(t))^2} + n\sum_{i=2}^{n-1}\frac{(v_{i-1}(t) - v_i(t))}{x_{i-1}(t) - x_i(t)}\left(\Phi'(n(x_{i-2}(t) - x_{i-1}(t))) - \Phi'(n(x_{i-1}(t) - x_i(t)))\right)$$

$$+ n\sum_{i=2}^{n-1}\frac{(v_{i-1}(t) - v_i(t))}{x_{i-1}(t) - x_i(t)}\left(-\Phi'(n(x_{i-1}(t) - x_i(t))) + \Phi'(n(x_i(t) - x_{i+1}(t)))\right)$$

$$+ n^2\sum_{i=2}^{n-1}\frac{(v_{i-1}(t) - v_i(t))}{x_{i-1}(t) - x_i(t)}\left(K'(n(x_{i-2}(t) - x_{i-1}(t)))(v_{i-2}(t) - v_{i-1}(t)) - K'(n(x_{i-1}(t) - x_i(t)))(v_{i-1}(t) - v_i(t))\right)$$

$$+ n^2\sum_{i=2}^{n-1}\frac{(v_{i-1}(t) - v_i(t))}{x_{i-1}(t) - x_i(t)}\left(-K'(n(x_{i-1}(t) - x_i(t)))(v_{i-1}(t) - v_i(t)) + K'(n(x_i(t) - x_{i+1}(t)))(v_i(t) - v_{i+1}(t))\right)$$

$$+ \frac{nv_1(t)}{L - x_1(t)}\left(\Phi'(n(L - x_1(t))) - \Phi'(n(x_1(t) - x_2(t)))\right) + \frac{nv_{n-1}(t)}{x_{n-1}(t)}\left(\Phi'(n(x_{n-2}(t) - x_{n-1}(t))) - \Phi'(nx_{n-1}(t))\right)$$

$$+ \frac{n^2 v_1(t)}{L - x_1(t)}\left(-K'(n(L - x_1(t)))v_1(t) + K'(n(x_1(t) - x_2(t)))(v_2(t) - v_1(t))\right)$$

$$+ \frac{n^2 v_{n-1}(t)}{x_{n-1}(t)}\left(K'(n(x_{n-2}(t) - x_{n-1}(t)))(v_{n-2}(t) - v_{n-1}(t)) - K'(nx_{n-1}(t))v_{n-1}(t)\right)$$

(4.34)

Rearranging terms in (4.34) we obtain:

$$\dot{Z}_n(t) = -\frac{1}{2}\sum_{i=1}^{n}\frac{\left(v_{i-1}(t)-v_i(t)\right)^3}{\left(x_{i-1}(t)-x_i(t)\right)^2}$$
$$+n\sum_{i=1}^{n-1}\left(\frac{v_i(t)-v_{i+1}(t)}{x_i(t)-x_{i+1}(t)}-\frac{v_{i-1}(t)-v_i(t)}{x_{i-1}(t)-x_i(t)}\right)\left(\Phi'(n(x_{i-1}(t)-x_i(t)))-\Phi'(n(x_i(t)-x_{i+1}(t)))\right)$$
$$+n^2\sum_{i=1}^{n-1}\left(\frac{v_i(t)-v_{i+1}(t)}{x_i(t)-x_{i+1}(t)}-\frac{v_{i-1}(t)-v_i(t)}{x_{i-1}(t)-x_i(t)}\right)K'(n(x_{i-1}(t)-x_i(t)))\left(v_{i-1}(t)-v_i(t)\right)$$
$$-n^2\sum_{i=1}^{n-1}\left(\frac{v_i(t)-v_{i+1}(t)}{x_i(t)-x_{i+1}(t)}-\frac{v_{i-1}(t)-v_i(t)}{x_{i-1}(t)-x_i(t)}\right)K'(n(x_i(t)-x_{i+1}(t)))\left(v_i(t)-v_{i+1}(t)\right)$$
(4.35)

Next define the function $G:(0,+\infty)\to(0,+\infty)$ by means of the formula:

$$G(x) = xK'(x), \text{ for all } x > 0 \tag{4.36}$$

Using definition (4.36), rearranging terms in (4.35) and using equations (3.1), (3.3), we get:

$$\dot{Z}_n(t) = -\frac{1}{2}\sum_{i=1}^{n-1}v_i(t)\left(\frac{v_i(t)-v_{i+1}(t)}{x_i(t)-x_{i+1}(t)}-\frac{v_{i-1}(t)-v_i(t)}{x_{i-1}(t)-x_i(t)}\right)\left(\frac{v_i(t)-v_{i+1}(t)}{x_i(t)-x_{i+1}(t)}+\frac{v_{i-1}(t)-v_i(t)}{x_{i-1}(t)-x_i(t)}\right)$$
$$+n\sum_{i=1}^{n-1}\left(\frac{v_i(t)-v_{i+1}(t)}{x_i(t)-x_{i+1}(t)}-\frac{v_{i-1}(t)-v_i(t)}{x_{i-1}(t)-x_i(t)}\right)\left(\Phi'(n(x_{i-1}(t)-x_i(t)))-\Phi'(n(x_i(t)-x_{i+1}(t)))\right)$$
$$-n\sum_{i=1}^{n-1}G(n(x_{i-1}(t)-x_i(t)))\left(\frac{v_i(t)-v_{i+1}(t)}{x_i(t)-x_{i+1}(t)}-\frac{v_{i-1}(t)-v_i(t)}{x_{i-1}(t)-x_i(t)}\right)^2$$
$$+n\sum_{i=1}^{n-1}\left(\frac{v_i(t)-v_{i+1}(t)}{x_i(t)-x_{i+1}(t)}-\frac{v_{i-1}(t)-v_i(t)}{x_{i-1}(t)-x_i(t)}\right)\frac{v_i(t)-v_{i+1}(t)}{x_i(t)-x_{i+1}(t)}\left(G(n(x_{i-1}(t)-x_i(t)))-G(n(x_i(t)-x_{i+1}(t)))\right)$$
(4.37)

Define
$$r = \min_{a\le x\le b}\left(G(x)\right) > 0 \tag{4.38}$$

Using (4.37) and definition (4.38) we get:

$$\dot{Z}_n(t) \le -rn\sum_{i=1}^{n-1}\left(\frac{v_i(t)-v_{i+1}(t)}{x_i(t)-x_{i+1}(t)}-\frac{v_{i-1}(t)-v_i(t)}{x_{i-1}(t)-x_i(t)}\right)^2$$
$$+\frac{1}{2}\sum_{i=1}^{n-1}|v_i(t)|\left|\frac{v_i(t)-v_{i+1}(t)}{x_i(t)-x_{i+1}(t)}-\frac{v_{i-1}(t)-v_i(t)}{x_{i-1}(t)-x_i(t)}\right|\left|\frac{v_i(t)-v_{i+1}(t)}{x_i(t)-x_{i+1}(t)}+\frac{v_{i-1}(t)-v_i(t)}{x_{i-1}(t)-x_i(t)}\right|$$
$$+n\sum_{i=1}^{n-1}\left|\frac{v_i(t)-v_{i+1}(t)}{x_i(t)-x_{i+1}(t)}-\frac{v_{i-1}(t)-v_i(t)}{x_{i-1}(t)-x_i(t)}\right|\left|\Phi'(n(x_{i-1}(t)-x_i(t)))-\Phi'(n(x_i(t)-x_{i+1}(t)))\right|$$
$$+n\sum_{i=1}^{n-1}\left|\frac{v_i(t)-v_{i+1}(t)}{x_i(t)-x_{i+1}(t)}-\frac{v_{i-1}(t)-v_i(t)}{x_{i-1}(t)-x_i(t)}\right|\left|\frac{v_i(t)-v_{i+1}(t)}{x_i(t)-x_{i+1}(t)}\right|\left|G(n(x_{i-1}(t)-x_i(t)))-G(n(x_i(t)-x_{i+1}(t)))\right|$$
(4.39)

Using the inequalities

$$\left|v_i(t)\right|\left|\frac{v_i(t)-v_{i+1}(t)}{x_i(t)-x_{i+1}(t)}-\frac{v_{i-1}(t)-v_i(t)}{x_{i-1}(t)-x_i(t)}\right|\left|\frac{v_i(t)-v_{i+1}(t)}{x_i(t)-x_{i+1}(t)}+\frac{v_{i-1}(t)-v_i(t)}{x_{i-1}(t)-x_i(t)}\right|$$

$$\leq \frac{2nr}{3}\left(\frac{v_i(t)-v_{i+1}(t)}{x_i(t)-x_{i+1}(t)}-\frac{v_{i-1}(t)-v_i(t)}{x_{i-1}(t)-x_i(t)}\right)^2 + \frac{3}{8nr}\left|v_i(t)\right|^2\left|\frac{v_i(t)-v_{i+1}(t)}{x_i(t)-x_{i+1}(t)}+\frac{v_{i-1}(t)-v_i(t)}{x_{i-1}(t)-x_i(t)}\right|^2$$

$$\left|\frac{v_i(t)-v_{i+1}(t)}{x_i(t)-x_{i+1}(t)}-\frac{v_{i-1}(t)-v_i(t)}{x_{i-1}(t)-x_i(t)}\right|\left|\Phi'(n(x_{i-1}(t)-x_i(t)))-\Phi'(n(x_i(t)-x_{i+1}(t)))\right|$$

$$\leq \frac{r}{6}\left(\frac{v_i(t)-v_{i+1}(t)}{x_i(t)-x_{i+1}(t)}-\frac{v_{i-1}(t)-v_i(t)}{x_{i-1}(t)-x_i(t)}\right)^2 + \frac{3}{2r}\left|\Phi'(n(x_{i-1}(t)-x_i(t)))-\Phi'(n(x_i(t)-x_{i+1}(t)))\right|^2$$

$$\left|\frac{v_i(t)-v_{i+1}(t)}{x_i(t)-x_{i+1}(t)}-\frac{v_{i-1}(t)-v_i(t)}{x_{i-1}(t)-x_i(t)}\right|\left|\frac{v_i(t)-v_{i+1}(t)}{x_i(t)-x_{i+1}(t)}\right|\left|G(n(x_{i-1}(t)-x_i(t)))-G(n(x_i(t)-x_{i+1}(t)))\right|$$

$$\leq \frac{r}{6}\left(\frac{v_i(t)-v_{i+1}(t)}{x_i(t)-x_{i+1}(t)}-\frac{v_{i-1}(t)-v_i(t)}{x_{i-1}(t)-x_i(t)}\right)^2$$

$$+\frac{3}{2r}\left|\frac{v_i(t)-v_{i+1}(t)}{x_i(t)-x_{i+1}(t)}\right|^2 \left|G(n(x_{i-1}(t)-x_i(t)))-G(n(x_i(t)-x_{i+1}(t)))\right|^2$$

we obtain from (4.39):

$$\dot{Z}_n(t) \leq -\frac{nr}{3}\sum_{i=1}^{n-1}\left(\frac{v_i(t)-v_{i+1}(t)}{x_i(t)-x_{i+1}(t)}-\frac{v_{i-1}(t)-v_i(t)}{x_{i-1}(t)-x_i(t)}\right)^2$$

$$+\frac{3}{16nr}\sum_{i=1}^{n-1}\left|v_i(t)\right|^2\left|\frac{v_i(t)-v_{i+1}(t)}{x_i(t)-x_{i+1}(t)}+\frac{v_{i-1}(t)-v_i(t)}{x_{i-1}(t)-x_i(t)}\right|^2 \qquad (4.40)$$

$$+\frac{3n}{2r}\sum_{i=1}^{n-1}\left|\Phi'(n(x_{i-1}(t)-x_i(t)))-\Phi'(n(x_i(t)-x_{i+1}(t)))\right|^2$$

$$+\frac{3n}{2r}\sum_{i=1}^{n-1}\left|\frac{v_i(t)-v_{i+1}(t)}{x_i(t)-x_{i+1}(t)}\right|^2\left|G(n(x_{i-1}(t)-x_i(t)))-G(n(x_i(t)-x_{i+1}(t)))\right|^2$$

Using the fact that $K'(x) > 0$ for all $x > 0$, we conclude that there exists a constant $M > 0$ such that the following inequalities hold for all $x, y \in [a,b]$:

$$\begin{aligned}\left|\Phi'(x)-\Phi'(y)\right| &\leq M\left|K(x)-K(y)\right| \\ \left|G(x)-G(y)\right| &\leq M\left|K(x)-K(y)\right|\end{aligned} \qquad (4.41)$$

Combining (4.40), (4.41), (4.6) and (4.31) we obtain:

$$\dot{Z}_n(t) \leq -\frac{nr}{3}\sum_{i=1}^{n-1}\left(\frac{v_i(t)-v_{i+1}(t)}{x_i(t)-x_{i+1}(t)} - \frac{v_{i-1}(t)-v_i(t)}{x_{i-1}(t)-x_i(t)}\right)^2$$

$$+\frac{3}{16nr}\sum_{i=1}^{n-1}|v_i(t)|^2\left|\frac{v_i(t)-v_{i+1}(t)}{x_i(t)-x_{i+1}(t)} + \frac{v_{i-1}(t)-v_i(t)}{x_{i-1}(t)-x_i(t)}\right|^2 + \frac{6M^2}{rm}(W_n(t)+E_n(t)) \quad (4.42)$$

$$+\frac{3M^2 n}{2r}\sum_{i=1}^{n-1}\left|\frac{v_i(t)-v_{i+1}(t)}{x_i(t)-x_{i+1}(t)}\right|^2 |K(n(x_{i-1}(t)-x_i(t)))-K(n(x_i(t)-x_{i+1}(t)))|^2$$

Using (3.1), (3.3), (4.5), definition (3.12), the fact that $\sum_{i=1}^{n}(x_{i-1}(t)-x_i(t)) = L$ and the Cauchy-Schwarz inequality, we get:

$$\max_{i=0,1,\ldots,n}(|v_i(t)|) \leq \sum_{i=1}^{n}|v_{i-1}(t)-v_i(t)| \leq \left(L\sum_{i=1}^{n}\frac{|v_{i-1}(t)-v_i(t)|^2}{x_{i-1}(t)-x_i(t)}\right)^{1/2} = \sqrt{2LZ_n(t)} \quad (4.43)$$

Consequently, we obtain from (4.42), (4.43):

$$\dot{Z}_n(t) \leq -\frac{nr}{3}\sum_{i=1}^{n-1}\left(\frac{v_i(t)-v_{i+1}(t)}{x_i(t)-x_{i+1}(t)} - \frac{v_{i-1}(t)-v_i(t)}{x_{i-1}(t)-x_i(t)}\right)^2$$

$$+\frac{3LZ_n(t)}{8nr}\sum_{i=1}^{n-1}\left|\frac{v_i(t)-v_{i+1}(t)}{x_i(t)-x_{i+1}(t)} + \frac{v_{i-1}(t)-v_i(t)}{x_{i-1}(t)-x_i(t)}\right|^2 + \frac{6M^2}{rm}(W_n(t)+E_n(t)) \quad (4.44)$$

$$+\frac{3M^2 n}{2r}\sum_{i=1}^{n-1}\left|\frac{v_i(t)-v_{i+1}(t)}{x_i(t)-x_{i+1}(t)}\right|^2 |K(n(x_{i-1}(t)-x_i(t)))-K(n(x_i(t)-x_{i+1}(t)))|^2$$

Using the inequality $\left|\frac{v_i(t)-v_{i+1}(t)}{x_i(t)-x_{i+1}(t)} + \frac{v_{i-1}(t)-v_i(t)}{x_{i-1}(t)-x_i(t)}\right|^2 \leq 2\left|\frac{v_i(t)-v_{i+1}(t)}{x_i(t)-x_{i+1}(t)}\right|^2 + 2\left|\frac{v_{i-1}(t)-v_i(t)}{x_{i-1}(t)-x_i(t)}\right|^2$ which gives $\sum_{i=1}^{n-1}\left|\frac{v_i(t)-v_{i+1}(t)}{x_i(t)-x_{i+1}(t)} + \frac{v_{i-1}(t)-v_i(t)}{x_{i-1}(t)-x_i(t)}\right|^2 \leq 4\sum_{i=1}^{n}\left|\frac{v_{i-1}(t)-v_i(t)}{x_{i-1}(t)-x_i(t)}\right|^2$, we get from (4.44) and (4.6):

$$\dot{Z}_n(t) \leq -\frac{nr}{3}\sum_{i=1}^{n-1}\left(\frac{v_i(t)-v_{i+1}(t)}{x_i(t)-x_{i+1}(t)} - \frac{v_{i-1}(t)-v_i(t)}{x_{i-1}(t)-x_i(t)}\right)^2 + \frac{3LZ_n(t)}{2nr}\sum_{i=1}^{n}\left|\frac{v_{i-1}(t)-v_i(t)}{x_{i-1}(t)-x_i(t)}\right|^2$$

$$+\frac{6M^2}{rm}(W_n(t)+E_n(t)) + \frac{6M^2}{rm}\max_{i=1,\ldots,n-1}\left(\left|\frac{v_i(t)-v_{i+1}(t)}{x_i(t)-x_{i+1}(t)}\right|^2\right)(W_n(t)+E_n(t)) \quad (4.45)$$

Using definition (3.12) and (4.31) we get:

$$n\min_{i=1,\ldots,n}\left(\left(\frac{v_{i-1}(t)-v_i(t)}{x_{i-1}(t)-x_i(t)}\right)^2\right) \leq \sum_{i=1}^{n}\left|\frac{v_{i-1}(t)-v_i(t)}{x_{i-1}(t)-x_i(t)}\right|^2 \leq \frac{2n}{a}Z_n(t) \quad (4.46)$$

Consequently, we obtain from (4.45) and (4.46):

$$\dot{Z}_n(t) \leq -\frac{nr}{3} \sum_{i=1}^{n-1} \left( \frac{v_i(t) - v_{i+1}(t)}{x_i(t) - x_{i+1}(t)} - \frac{v_{i-1}(t) - v_i(t)}{x_{i-1}(t) - x_i(t)} \right)^2 + \frac{6M^2}{rm}(W_n(t) + E_n(t))$$
$$+ \frac{3L}{ar} Z_n^2(t) + \frac{6M^2}{rm} \max_{i=1,\ldots,n} \left( \left| \frac{v_i(t) - v_{i-1}(t)}{x_i(t) - x_{i-1}(t)} \right|^2 \right)(W_n(t) + E_n(t))$$
(4.47)

Using (4.4), the Cauchy-Schwarz inequality, the inequality $2xy \leq \frac{1}{\varepsilon n} x^2 + \varepsilon n y^2$ that holds for every $\varepsilon > 0, x, y \in \mathbb{R}$, definition (3.12) and (4.31), we get:

$$\max_{i=1,\ldots,n} \left( \left( \frac{v_{i-1}(t) - v_i(t)}{x_{i-1}(t) - x_i(t)} \right)^2 \right) - \min_{i=1,\ldots,n} \left( \left( \frac{v_{i-1}(t) - v_i(t)}{x_{i-1}(t) - x_i(t)} \right)^2 \right)$$
$$\leq \sum_{i=1}^{n-1} \left| \left( \frac{v_i(t) - v_{i+1}(t)}{x_i(t) - x_{i+1}(t)} \right)^2 - \left( \frac{v_{i-1}(t) - v_i(t)}{x_{i-1}(t) - x_i(t)} \right)^2 \right|$$
$$= \sum_{i=1}^{n-1} \left| \frac{v_i(t) - v_{i+1}(t)}{x_i(t) - x_{i+1}(t)} - \frac{v_{i-1}(t) - v_i(t)}{x_{i-1}(t) - x_i(t)} \right| \left| \frac{v_i(t) - v_{i+1}(t)}{x_i(t) - x_{i+1}(t)} + \frac{v_{i-1}(t) - v_i(t)}{x_{i-1}(t) - x_i(t)} \right|$$
$$\leq 2 \left( \sum_{i=1}^{n} \left| \frac{v_{i-1}(t) - v_i(t)}{x_{i-1}(t) - x_i(t)} \right|^2 \right)^{1/2} \left( \sum_{i=1}^{n-1} \left| \frac{v_i(t) - v_{i+1}(t)}{x_i(t) - x_{i+1}(t)} - \frac{v_{i-1}(t) - v_i(t)}{x_{i-1}(t) - x_i(t)} \right|^2 \right)^{1/2}$$
$$\leq \frac{1}{\varepsilon n} \sum_{i=1}^{n} \left| \frac{v_{i-1}(t) - v_i(t)}{x_{i-1}(t) - x_i(t)} \right|^2 + \varepsilon n \sum_{i=1}^{n-1} \left| \frac{v_i(t) - v_{i+1}(t)}{x_i(t) - x_{i+1}(t)} - \frac{v_{i-1}(t) - v_i(t)}{x_{i-1}(t) - x_i(t)} \right|^2$$
$$\leq \frac{2}{\varepsilon a} Z_n(t) + \varepsilon n \sum_{i=1}^{n-1} \left| \frac{v_i(t) - v_{i+1}(t)}{x_i(t) - x_{i+1}(t)} - \frac{v_{i-1}(t) - v_i(t)}{x_{i-1}(t) - x_i(t)} \right|^2$$
(4.48)

Estimate (4.48) in conjunction with (4.46) gives for every $\varepsilon > 0$:

$$\max_{i=1,\ldots,n} \left( \left( \frac{v_{i-1}(t) - v_i(t)}{x_{i-1}(t) - x_i(t)} \right)^2 \right) \leq \frac{2}{a}\left(1 + \frac{1}{\varepsilon}\right) Z_n(t) + \varepsilon n \sum_{i=1}^{n-1} \left| \frac{v_i(t) - v_{i+1}(t)}{x_i(t) - x_{i+1}(t)} - \frac{v_{i-1}(t) - v_i(t)}{x_{i-1}(t) - x_i(t)} \right|^2$$
(4.49)

Combining (4.47) and (4.49) we get for every $\varepsilon > 0$:

$$\dot{Z}_n(t) \leq -\left( \frac{r}{3} - \varepsilon \frac{6M^2}{rm}(W_n(t) + E_n(t)) \right) n \sum_{i=1}^{n-1} \left( \frac{v_i(t) - v_{i+1}(t)}{x_i(t) - x_{i+1}(t)} - \frac{v_{i-1}(t) - v_i(t)}{x_{i-1}(t) - x_i(t)} \right)^2$$
$$+ \frac{3L}{ar} Z_n^2(t) + \frac{6M^2}{rm}(W_n(t) + E_n(t)) + \left(1 + \frac{1}{\varepsilon}\right) \frac{12M^2}{arm}(W_n(t) + E_n(t)) Z_n(t)$$
(4.50)

Selecting $\varepsilon = \frac{r^2 m}{36M^2(W_n(t) + E_n(t))}$ when $W_n(t) + E_n(t) > 0$, we get from (4.50):

$$\dot{Z}_n(t) \leq -\frac{rn}{6}\sum_{i=1}^{n-1}\left(\frac{v_i(t)-v_{i+1}(t)}{x_i(t)-x_{i+1}(t)} - \frac{v_{i-1}(t)-v_i(t)}{x_{i-1}(t)-x_i(t)}\right)^2 + \frac{6M^2}{rm}(W_n(t)+E_n(t))$$
$$+\frac{3L}{ar}Z_n^2(t) + \left(36M^2(W_n(t)+E_n(t))+r^2m\right)\frac{12M^2}{ar^3m^2}(W_n(t)+E_n(t))Z_n(t) \qquad (4.51)$$

Inequality (4.51) holds when $W_n(t)+E_n(t)=0$ as well. Therefore, inequality (4.51) holds in any case. Using the inequality

$$\left(36M^2(W_n(t)+E_n(t))+r^2m\right)\frac{12M^2}{ar^3m^2}(W_n(t)+E_n(t))Z_n(t)$$
$$\leq Z_n^2(t) + \left(36M^2(W_n(t)+E_n(t))+r^2m\right)^2\left(\frac{6M^2}{ar^3m^2}\right)^2(W_n(t)+E_n(t))^2$$

it follows from (4.51) that the following estimate holds:

$$\dot{Z}_n(t) \leq -\frac{rn}{6}\sum_{i=1}^{n-1}\left(\frac{v_i(t)-v_{i+1}(t)}{x_i(t)-x_{i+1}(t)} - \frac{v_{i-1}(t)-v_i(t)}{x_{i-1}(t)-x_i(t)}\right)^2 + \left(1+\frac{3L}{ar}\right)Z_n^2(t) + B(W_n(t)+E_n(t)) \qquad (4.52)$$

where $B(s) := \frac{6M^2}{rm}s + (36M^2s+r^2m)^2\left(\frac{6M^2}{ar^3m^2}\right)^2 s^2$, for all $s>0$. Inequality (4.33) is a direct consequence of (4.52).

The fact that $\sum_{i=1}^{n}(x_{i-1}(t)-x_i(t))\frac{v_{i-1}(t)-v_i(t)}{x_{i-1}(t)-x_i(t)} = v_0(t)-v_n(t)=0$ implies that $\min_{i=1,\ldots,n}\left(\frac{v_{i-1}(t)-v_i(t)}{x_{i-1}(t)-x_i(t)}\right) \leq 0 \leq \max_{i=1,\ldots,n}\left(\frac{v_{i-1}(t)-v_i(t)}{x_{i-1}(t)-x_i(t)}\right)$. Consequently, we obtain from Fact 4, (4.5) and the Cauchy-Schwarz inequality:

$$\max_{i=1,\ldots,n}\left(\left|\frac{v_{i-1}(t)-v_i(t)}{x_{i-1}(t)-x_i(t)}\right|\right) \leq \sum_{i=1}^{n-1}\left|\frac{v_i(t)-v_{i+1}(t)}{x_i(t)-x_{i+1}(t)} - \frac{v_{i-1}(t)-v_i(t)}{x_{i-1}(t)-x_i(t)}\right|$$
$$\leq \left(n\sum_{i=1}^{n-1}\left|\frac{v_i(t)-v_{i+1}(t)}{x_i(t)-x_{i+1}(t)} - \frac{v_{i-1}(t)-v_i(t)}{x_{i-1}(t)-x_i(t)}\right|^2\right)^{1/2} \qquad (4.53)$$

Definition (3.12), (4.53), the Cauchy-Schwarz inequality and the fact that $\sum_{i=1}^{n}(x_{i-1}(t)-x_i(t))=L$ give the estimate

$$2Z_n(t) \leq \max_{i=1,\ldots,n}\left(\left|\frac{v_{i-1}(t)-v_i(t)}{x_{i-1}(t)-x_i(t)}\right|\right)\sum_{i=1}^{n}|v_{i-1}(t)-v_i(t)|$$
$$\leq \left(n\sum_{i=1}^{n-1}\left|\frac{v_i(t)-v_{i+1}(t)}{x_i(t)-x_{i+1}(t)} - \frac{v_{i-1}(t)-v_i(t)}{x_{i-1}(t)-x_i(t)}\right|^2\right)^{1/2}\left(L\sum_{i=1}^{n}\frac{|v_{i-1}(t)-v_i(t)|^2}{x_{i-1}(t)-x_i(t)}\right)^{1/2}$$
$$= \left(n\sum_{i=1}^{n-1}\left|\frac{v_i(t)-v_{i+1}(t)}{x_i(t)-x_{i+1}(t)} - \frac{v_{i-1}(t)-v_i(t)}{x_{i-1}(t)-x_i(t)}\right|^2\right)^{1/2}\sqrt{2LZ_n(t)}$$

which directly implies the estimate

$$2Z_n(t) \le Ln \sum_{i=1}^{n-1} \left| \frac{v_i(t) - v_{i+1}(t)}{x_i(t) - x_{i+1}(t)} - \frac{v_{i-1}(t) - v_i(t)}{x_{i-1}(t) - x_i(t)} \right|^2 \quad (4.54)$$

Consequently, we obtain from (4.52) and (4.54) the estimate:

$$\dot{Z}_n(t) \le -\frac{r}{3L} Z_n(t) + \left(1 + \frac{3L}{ar}\right) Z_n^2(t) + B\left(W_n(t) + E_n(t)\right) \quad (4.55)$$

Estimate (4.55) directly implies estimate (4.32). The proof is complete. ◁

We are now ready to provide the proof of Theorem 3.1.

**Proof of Theorem 3.1:** Estimates (3.20), (3.21), (3.22), (3.23) are consequences of (4.1), (4.2), (4.6), (4.7), (4.8), (3.15) and definitions (3.10), (3.11), (3.14).

We next show estimates (3.24), (3.26). Notice that (4.1) and definitions (3.12), (4.36), (4.38) imply the following estimate

$$\dot{E}_n(t) = -m \sum_{i=1}^{n} G(n(x_{i-1}(t) - x_i(t))) \frac{(v_{i-1}(t) - v_i(t))^2}{x_{i-1}(t) - x_i(t)} \le -2mrZ_n(t)$$

from which we get for all $t \ge 0$:

$$E_n(t) + 2mr \int_0^t Z_n(s)ds \le E_n(0) \quad (4.56)$$

Lemma 4.1 implies that there exist constants $\sigma, \omega, r > 0$ (independent of $n \ge 2$ and $r > 0$ defined by (4.38)) and a function $B \in K_\infty$ (independent of $n \ge 2$) such that the differential inequalities (4.32), (4.33) hold for all $t \ge 0$. Using (4.32) and defining $y(t) = \exp\left(\sigma t - \omega \int_0^t Z_n(s)ds\right) Z_n(t)$, we get for all $t \ge 0$:

$$\dot{y}(t) \le \exp\left(\sigma t - \omega \int_0^t Z_n(s)ds\right) B\left(W_n(t) + E_n(t)\right)$$

which directly implies the following estimate

$$Z_n(t) \le \exp\left(-\sigma t + \omega \int_0^t Z_n(s)ds\right) Z_n(0)$$

$$+ \int_0^t B(W_n(\tau) + E_n(\tau)) \exp\left(-\sigma(t - \tau) + \omega \int_\tau^t Z_n(s)ds\right) d\tau$$

The above estimate combined with the fact that $W_n(t) + E_n(t) \le W_n(0) + E_n(0)$ for all $t \ge 0$ (a consequence of (4.1) and (4.2)) and (4.56) gives:

$$Z_n(t) \le \exp\left(-\sigma t + \omega \int_0^t Z_n(s)ds\right) Z_n(0)$$

$$+ B(W_n(0) + E_n(0)) \int_0^t \exp\left(-\sigma(t-\tau) + \omega \int_0^t Z_n(s)ds\right) d\tau$$

$$\le \exp\left(\frac{\omega E_n(0)}{2mr}\right) Z_n(0) + B(W_n(0) + E_n(0)) \exp\left(\frac{\omega E_n(0)}{2mr}\right) \int_0^t \exp(-\sigma(t-\tau)) d\tau$$

$$\le \exp\left(\frac{\omega E_n(0)}{2mr}\right) \left(Z_n(0) + \frac{1}{\sigma} B(W_n(0) + E_n(0))\right)$$

It is clear that the above estimate and (3.15) implies (3.24) with $R_2 = \exp\left(\frac{\omega \bar{E}}{2mr}\right) \left(\bar{Z} + \frac{1}{\sigma} B(\bar{W} + \bar{E})\right)$. Moreover, estimates (4.43) and (3.24) imply estimate (3.26) with $R_4 = \sqrt{2LR_2}$.

We next show estimate (3.27). Integrating (4.33) we get for all $t \ge 0$:

$$Z_n(t) + \frac{rn}{6} \int_0^t \sum_{i=1}^{n-1} \left(\frac{v_i(s) - v_{i+1}(s)}{x_i(s) - x_{i+1}(s)} - \frac{v_{i-1}(s) - v_i(s)}{x_{i-1}(s) - x_i(s)}\right)^2 ds$$

$$\le Z_n(0) + \omega \int_0^t Z_n^2(s)ds + \int_0^t B(W_n(s) + E_n(s))ds \quad (4.57)$$

$$\le R_2 + \omega R_2^2 t + B(\bar{W} + \bar{E})t$$

Estimate (4.57) combined with (3.20), (3.21), (3.24) gives estimate (3.27) with $R_5 = \frac{6}{r} R_2$ and $R_6 = \frac{6}{r}(\omega R_2^2 + B(\bar{W} + \bar{E}))$.

Finally, we show estimate (3.25). Equations (3.1), (3.2), (3.3) imply the following equation for all $t \ge 0$ and $i = 1, ..., n-1$:

$$\frac{d}{dt}\left(K(n(x_{i-1}(t) - x_i(t))) - K(n(x_i(t) - x_{i+1}(t)))\right)$$
$$= \frac{1}{n}\dot{v}_i(t) + \Phi'(n(x_i(t) - x_{i+1}(t))) - \Phi'(n(x_{i-1}(t) - x_i(t))) \quad (4.58)$$

Integrating (4.58) we get for all $t \ge 0$ and $i = 1, ..., n-1$:

$$n\left(K(n(x_{i-1}(t) - x_i(t))) - K(n(x_i(t) - x_{i+1}(t)))\right)$$
$$= n\left(K(n(x_{i-1}(0) - x_i(0))) - K(n(x_i(0) - x_{i+1}(0)))\right) + v_i(t) - v_i(0) \quad (4.59)$$
$$+ n \int_0^t \left(\Phi'(n(x_i(s) - x_{i+1}(s))) - \Phi'(n(x_{i-1}(s) - x_i(s)))\right) ds$$

Using the fact that $K'(x) > 0$ we conclude that there exists a constant $M > 0$ such that the following inequalities hold for all $x, y \in [a, b]$:

$$|\Phi'(x) - \Phi'(y)| \leq M |K(x) - K(y)| \tag{4.60}$$

Combining (4.59) and (3.15), (3.26), (4.60) we obtain for all $t \geq 0$ and $i = 1, \ldots, n-1$:

$$\begin{aligned} & n \left| K\left(n(x_{i-1}(t) - x_i(t))\right) - K\left(n(x_i(t) - x_{i+1}(t))\right) \right| \\ & \leq \bar{A} + 2R_4 + M \int_0^t n \left| K\left(n(x_i(s) - x_{i+1}(s))\right) - K\left(n(x_{i-1}(s) - x_i(s))\right) \right| ds \end{aligned} \tag{4.61}$$

Applying the Gronwall-Bellman Lemma in (4.61) we obtain for all $t \geq 0$ and $i = 1, \ldots, n-1$:

$$n \left| K\left(n(x_{i-1}(t) - x_i(t))\right) - K\left(n(x_i(t) - x_{i+1}(t))\right) \right| \leq \left(\bar{A} + 2R_4\right) \exp(Mt)$$

The above estimate and definition (3.13) imply estimate (3.25) with $\kappa = M$ and $R_3 = \bar{A} + 2R_4$. The proof is complete. ◁

Next, we provide the proof of Theorem 3.2.

**Proof of Theorem 3.2:** Let $(\rho_0, v_0) \in X$ be given and suppose that (3.30) holds. We set $(x_1(0), \ldots, x_{n-1}(0), v_1(0), \ldots, v_{n-1}(0)) = P_n(\rho_0, v_0) \in \Omega_n$, where $x_n(0) = 0 < x_{n-1}(0) < \ldots < x_1(0) < x_0(0) = L$ is the unique set of points with $\int_{x_i(0)}^{x_{i-1}(0)} \rho_0(x) dx = \frac{m}{n}$, $v_i(0) = v_0(x_i(0))$ for $i = 0, 1, \ldots, n$. Using the Cauchy-Schwarz inequality, the fact that $\frac{m}{\|\rho_0\|_\infty} \leq n(x_{i-1}(0) - x_i(0)) \leq \frac{m}{\rho_{\min}}$ for $i = 1, \ldots, n$ with $\rho_{\min} = \min_{0 \leq x \leq L}(\rho_0(x))$, the fact that $\Phi: (0, +\infty) \to \mathbb{R}$ is a decreasing function and the fact that $v_0(0) = 0$, we get:

$$\begin{aligned} & \frac{m}{2n} \sum_{i=1}^n v_i^2(0) + \frac{m}{n} \sum_{i=1}^n \Phi\left(n(x_{i-1}(0) - x_i(0))\right) \\ & \leq \frac{m}{2n} \sum_{i=1}^n \left( \int_0^{x_i(0)} v_0'(s) ds \right)^2 + \frac{m}{n} \sum_{i=1}^n \Phi\left(\frac{m}{\|\rho_0\|_\infty}\right) \leq \frac{m}{2n} \sum_{i=1}^n x_i(0) \|v_0'\|_2^2 + m\Phi\left(\frac{m}{\|\rho_0\|_\infty}\right) \\ & \leq \frac{mL}{2} \|v_0'\|_2^2 + m\Phi\left(\frac{m}{\|\rho_0\|_\infty}\right) \end{aligned} \tag{4.62}$$

Since $\rho_0$ is continuous, for every $i = 1, \ldots, n$ there exists $\xi_i \in [x_i(0), x_{i-1}(0)]$ with $n(x_{i-1}(0) - x_i(0)) = \frac{m}{\rho(\xi_i)}$. Setting $\bar{M} = \max\left\{ K'(s) : \frac{m}{\|\rho_0\|_\infty} \leq s \leq \frac{m}{\rho_{\min}} \right\}$ and using the fact that $K'(x) > 0$ for all $x > 0$ and the fact that $\frac{m}{\|\rho_0\|_\infty} \leq n(x_{i-1}(0) - x_i(0)) \leq \frac{m}{\rho_{min}}$ for $i = 1,...,n$ with $\rho_{min} = \min_{0 \leq x \leq L}(\rho_0(x))$, we get for all $i = 1,...,n-1$:

$$n\left|K(n(x_{i-1}(0) - x_i(0))) - K(n(x_i(0) - x_{i+1}(0)))\right| = n\left|K\left(\frac{m}{\rho(\xi_i)}\right) - K\left(\frac{m}{\rho(\xi_{i+1})}\right)\right|$$
$$\leq mn\bar{M}\left|\frac{1}{\rho(\xi_i)} - \frac{1}{\rho(\xi_{i+1})}\right| \leq mn\bar{M}\frac{\|\rho_0'\|_\infty}{\rho_{min}^2}|\xi_i - \xi_{i+1}| \quad (4.63)$$
$$\leq mn\bar{M}\frac{\|\rho_0'\|_\infty}{\rho_{min}^2}(x_{i-1}(0) - x_{i+1}(0)) \leq 2m^2\bar{M}\frac{\|\rho_0'\|_\infty}{\rho_{min}^3}$$

It follows from (4.63) that

$$\max_{i=1,...,n-1}\left(n\left|K(n(x_{i-1}(0) - x_i(0))) - K(n(x_i(0) - x_{i+1}(0)))\right|\right) \leq 2m^2\bar{M}\frac{\|\rho_0'\|_\infty}{\rho_{min}^3} \quad (4.64)$$

Using the Cauchy-Schwarz inequality, we get:

$$\frac{1}{2}\sum_{i=1}^n \frac{(v_{i-1}(0) - v_i(0))^2}{x_{i-1}(0) - x_i(0)} = \frac{1}{2}\sum_{i=1}^n \frac{\left(\int_{x_i(0)}^{x_{i-1}(0)} v_0'(s)ds\right)^2}{x_{i-1}(0) - x_i(0)} \leq \frac{1}{2}\sum_{i=1}^n \int_{x_i(0)}^{x_{i-1}(0)} |v_0'(s)|^2 ds = \frac{1}{2}\|v_0'\|_2^2 \quad (4.65)$$

Furthermore, (4.63) and the fact that $\|v_0\|_\infty \leq \sqrt{L}\|v_0'\|_2$ gives for all $i = 1,...,n-1$:

$$|v_i(0) - nK(n(x_{i-1}(0) - x_i(0))) + nK(n(x_i(0) - x_{i+1}(0)))|$$
$$\leq |v_0(x_i(0))| + n|K(n(x_{i-1}(0) - x_i(0))) - K(n(x_i(0) - x_{i+1}(0)))| \quad (4.66)$$
$$\leq \sqrt{L}\|v_0'\|_2 + 2m^2\bar{M}\frac{\|\rho_0'\|_\infty}{\rho_{min}^3}$$

Inequality (4.66) implies that

$$|v_i(0) - nK(n(x_{i-1}(0) - x_i(0))) + nK(n(x_i(0) - x_{i+1}(0)))|^2 \leq 2L\|v_0'\|_2^2 + 4m^4\bar{M}^2\frac{\|\rho_0'\|_\infty^2}{\rho_{min}^6}$$

for all $i = 1,...,n-1$. Consequently, using the fact that $\frac{m}{\|\rho_0\|_\infty} \leq n(x_{i-1}(0) - x_i(0)) \leq \frac{m}{\rho_{min}}$ for $i = 1,...,n$ with $\rho_{min} = \min_{0 \leq x \leq L}(\rho_0(x))$, the fact that $\Phi : (0, +\infty) \to \mathbb{R}$ is a decreasing function, we get:

$$\frac{m}{2n}\sum_{i=1}^{n-1}\left(v_i(0)-nK(n(x_{i-1}(0)-x_i(0)))+nK(n(x_i(0)-x_{i+1}(0)))\right)^2+\frac{m}{n}\sum_{i=1}^{n}\Phi\left(n\left(x_{i-1}(0)-x_i(0)\right)\right) \quad (4.67)$$

$$\leq mL\|v_0'\|_2^2+2m^5\bar{M}^2\frac{\|\rho_0'\|_\infty^2}{\rho_{\min}^6}+m\Phi\left(\frac{m}{\|\rho_0\|_\infty}\right)$$

We conclude from (4.62), (4.64), (4.65) and (4.67) that if we are given $(\rho_0,v_0)\in X$ then, no matter what $n\geq 2$ is, the vector $(x_1(0),...,x_{n-1}(0),v_1(0),...,v_{n-1}(0))=P_n(\rho_0,v_0)\in\Omega_n$ satisfies inequalities

(3.15) with $\bar{E}=\dfrac{mL}{2}\|v_0'\|_2^2+m\Phi\left(\dfrac{m}{\|\rho_0\|_\infty}\right)$, $\bar{A}=2m^2\bar{M}\dfrac{\|\rho_0'\|_\infty}{\rho_{\min}^3}$, $\bar{Z}=\dfrac{1}{2}\|v_0'\|_2^2$,

$\bar{W}=mL\|v_0'\|_2^2+2m^5\bar{M}^2\dfrac{\|\rho_0'\|_\infty^2}{\rho_{\min}^6}+m\Phi\left(\dfrac{m}{\|\rho_0\|_\infty}\right)$, $\rho_{\min}=\min\limits_{0\leq x\leq L}(\rho_0(x))$ and

$\bar{M}=\max\left\{K'(s):\dfrac{m}{\|\rho_0\|_\infty}\leq s\leq\dfrac{m}{\rho_{\min}}\right\}$. Therefore, by virtue of (3.30) we get (3.19). The proof is complete. ◁

In order to prove Theorem 3.3 we first need to show some auxiliary technical results.

**Lemma 4.2:** *Suppose that Assumption (A) holds. Let $(\rho_0,v_0)\in X$ that satisfies (3.30) be given. Consider the unique solution of (3.1), (3.2), (3.3) with initial condition $(x_1(0),...,x_{n-1}(0),v_1(0),...,v_{n-1}(0))=P_n(\rho_0,v_0)\in\Omega_n$ and define $v^{(n)},\rho^{(n)}:\mathbb{R}_+\times[0,L]\to\mathbb{R}$ by means of (3.31)-(3.34). Then for every $T>0$ it holds that*

$$\rho^{(n)}\in L^\infty\left((0,T);H^1((0,L))\right),\ v^{(n)}\in L^\infty\left((0,T);H_0^1((0,L))\right)$$

*Moreover, it holds that $\rho_x^{(n)}\in L^\infty((0,T)\times(0,L))$.*

**Proof:** The functions $v^{(n)},\rho^{(n)}:\mathbb{R}_+\times[0,L]\to\mathbb{R}$ are $C^1$ on the set $\bigcup\limits_{t\geq 0}\{(t,x):x\in[0,L],x\neq x_i(t),i=0,...,n\}$ (are $C^1$ a.e. on $\mathbb{R}_+\times[0,L]$). Furthermore, for each $t\geq 0$ the functions $v^{(n)}[t],\rho^{(n)}[t]:[0,L]\to\mathbb{R}$ are of class $W^{1,\infty}((0,L))$ with $v^{(n)}[t]\in H_0^1((0,L))$. Moreover, inequalities (3.23), (3.26) and definitions (3.31)-(3.34) guarantee the following estimates for all $t\geq 0$:

$$\frac{m}{b}\leq\min_{0\leq x\leq L}\left(\rho^{(n)}(t,x)\right)\leq\max_{0\leq x\leq L}\left(\rho^{(n)}(t,x)\right)=\left\|\rho^{(n)}[t]\right\|_\infty\leq\frac{m}{a} \quad (4.68)$$

$$\left\|v^{(n)}[t]\right\|_\infty\leq R_4 \quad (4.69)$$

Using (3.22), (3.23), (3.25), (3.13), (3.31), (3.32) and the fact that $K'(x)>0$ for all $x>0$, we get for all $t\geq 0$:

$$R_1 \geq n \sum_{i=1}^{n-1} \left( K(n(x_{i-1}(t) - x_i(t))) - K(n(x_i(t) - x_{i+1}(t))) \right)^2$$

$$\geq n \left( \min_{a \leq s \leq b} (K'(s)) \right)^2 \sum_{i=1}^{n-1} \left( \frac{m}{\rho_i(t)} - \frac{m}{\rho_{i+1}(t)} \right)^2 = nm^2 \left( \min_{a \leq s \leq b} (K'(s)) \right)^2 \sum_{i=1}^{n-1} \left( \frac{\rho_{i+1}(t) - \rho_i(t)}{\rho_i(t) \rho_{i+1}(t)} \right)^2$$

$$\geq n \left( \min_{a \leq s \leq b} (K'(s)) \right)^2 \frac{a^4}{m^2} \sum_{i=1}^{n-1} (\rho_{i+1}(t) - \rho_i(t))^2$$

$$R_3 \exp(\kappa t) \geq \max_{i=1,\ldots,n-1} \left( n \left| K(n(x_{i-1}(t) - x_i(t))) - K(n(x_i(t) - x_{i+1}(t))) \right| \right)$$

$$\geq \left( \min_{a \leq s \leq b} (K'(s)) \right) \max_{i=1,\ldots,n-1} \left( n \left| \frac{m}{\rho_i(t)} - \frac{m}{\rho_{i+1}(t)} \right| \right) = m \left( \min_{a \leq s \leq b} (K'(s)) \right) \max_{i=1,\ldots,n-1} \left( n \left| \frac{\rho_{i+1}(t) - \rho_i(t)}{\rho_i(t) \rho_{i+1}(t)} \right| \right)$$

$$\geq \left( \min_{a \leq s \leq b} (K'(s)) \right) \frac{a^2}{m} \max_{i=1,\ldots,n-1} \left( n |\rho_{i+1}(t) - \rho_i(t)| \right)$$

We conclude from the above inequalities and (3.32) that there exist constants $R_7, R_8 > 0$ (independent of $n \geq 2$ but dependent on $\rho_0, v_0$) such that the following estimates hold for all $t \geq 0$:

$$n \sum_{i=0}^{n-1} (\rho_{i+1}(t) - \rho_i(t))^2 \leq R_7 \tag{4.70}$$

$$\max_{i=0,\ldots,n-1} \left( n |\rho_{i+1}(t) - \rho_i(t)| \right) \leq R_8 \exp(\kappa t) \tag{4.71}$$

Using (3.12), (3.23) and definitions (3.33), it follows that for every $t \geq 0$ we get:

$$\|v_x^{(n)}[t]\|_2^2 = \int_0^L \left( v_x^{(n)}(t, x) \right)^2 dx = \sum_{i=1}^n \int_{x_i(t)}^{x_{i-1}(t)} \left( v_x^{(n)}(t, x) \right)^2 dx = \sum_{i=1}^n \frac{(v_{i-1}(t) - v_i(t))^2}{x_{i-1}(t) - x_i(t)} = Z_n(t) \tag{4.72}$$

$$\|\rho_x^{(n)}[t]\|_2^2 = \int_0^L \left( \rho_x^{(n)}(t, x) \right)^2 dx = \sum_{i=1}^n \int_{x_i(t)}^{x_{i-1}(t)} \left( \rho_x^{(n)}(t, x) \right)^2 dx$$
$$= \sum_{i=1}^n \frac{(\rho_{i-1}(t) - \rho_i(t))^2}{x_{i-1}(t) - x_i(t)} \leq \frac{n}{a} \sum_{i=1}^n (\rho_{i-1}(t) - \rho_i(t))^2 \tag{4.73}$$

$$\|\rho_x^{(n)}[t]\|_\infty = \max_{i=1,\ldots,n} \left( \frac{|\rho_{i-1}(t) - \rho_i(t)|}{x_{i-1}(t) - x_i(t)} \right) \leq \frac{n}{a} \max_{i=1,\ldots,n} \left( |\rho_{i-1}(t) - \rho_i(t)| \right) \tag{4.74}$$

We conclude from (3.24), (4.70)-(4.74) that there exists a constant $R_9 > 0$ (independent of $n \geq 2$ but dependent on $\rho_0, v_0$) such that the following estimates hold for all $t \geq 0$:

$$\|v_x^{(n)}[t]\|_2^2 \leq R_9 \tag{4.75}$$

$$\|\rho_x^{(n)}[t]\|_2^2 \leq R_9 \tag{4.76}$$

$$\left\|\rho_x^{(n)}[t]\right\|_\infty \leq R_9 \exp(\kappa t) \qquad (4.77)$$

Next consider a continuous linear functional $F$ on $H^1((0,L))$ which by virtue of Remark 21 on page 220 in [1] is given by

$$F(u) = \int_0^L f_0(x)u(x)dx + \int_0^L f_1(x)u'(x)dx, \text{ for all } u \in H^1((0,L))$$

for some $f_0, f_1 \in L^2((0,L))$. Using definition (3.33) we get

$$F(\rho^{(n)}[t]) = \sum_{i=1}^n \rho_i(t) \int_{x_i(t)}^{x_{i-1}(t)} f_0(x)dx + \sum_{i=1}^n \frac{\rho_{i-1}(t) - \rho_i(t)}{x_{i-1}(t) - x_i(t)} \int_{x_i(t)}^{x_{i-1}(t)} \left(f_0(x)(x - x_i(t)) + f_1(x)\right)dx$$

from which we can conclude that the mapping $t \to F(\rho^{(n)}[t])$ is continuous for every continuous linear functional $F$ on $H^1((0,L))$. Hence the mapping $t \to \rho^{(n)}[t] \in H^1((0,L))$ is weakly continuous and Corollary 1.4.8 on page 6 in [2] implies that the mapping $t \to \rho^{(n)}[t] \in H^1((0,L))$ is measurable. Therefore, we conclude from (4.68), (4.69), (4.75), (4.76) that for every $T > 0$ we have $\rho^{(n)} \in L^\infty((0,T); H^1((0,L)))$. Similarly, we can prove (using definitions (3.33) and (3.34)) that for every $T > 0$ we have $v^{(n)} \in L^\infty((0,T); H_0^1((0,L)))$. Moreover, from (4.77) we conclude that for every $T > 0$ we have $\rho_x^{(n)} \in L^\infty((0,T) \times (0,L))$. The proof is complete. ◁

**Lemma 4.3:** *Suppose that Assumption (A) holds. Let $(\rho_0, v_0) \in X$ that satisfies (3.30) be given. Consider the unique solution of (3.1), (3.2), (3.3) with initial condition $(x_1(0),...,x_{n-1}(0), v_1(0),...,v_{n-1}(0)) = P_n(\rho_0, v_0) \in \Omega_n$ and define $v^{(n)}, \rho^{(n)} : \mathbb{R}_+ \times [0,L] \to \mathbb{R}$ by means of (3.31)-(3.34). Then for every $T > 0$ and for every $\varphi \in C^1([0,T] \times [0,L])$ with $\varphi(T,x) = 0$ for all $x \in [0,L]$ there exists a constant $R_{10} > 0$ (independent of $n \geq 2$ but dependent on $\rho_0, v_0, T, \varphi$) such that the following estimate holds for all $n \geq 2$:*

$$\left| \int_0^L \varphi(0,x)\rho_0(x)dx + \int_0^T \int_0^L \rho^{(n)}(t,x)\left(\varphi_t(t,x) + v^{(n)}(t,x)\varphi_x(t,x)\right)dxdt \right| \leq \frac{R_{10}}{n} \qquad (4.78)$$

**Remark:** Estimate (4.78) shows how accurately the approximate solution $v^{(n)}, \rho^{(n)}$ satisfies the continuity equation. Notice that the accuracy provided by the particle scheme is of order $1/n$.

**Proof:** Let $T > 0$ (arbitrary) be given and let $\varphi \in C^1([0,T] \times [0,L])$ with $\varphi(T,x) = 0$ for all $x \in [0,L]$ (but otherwise arbitrary) be given. We get

$$\int_0^L \varphi(0,x)\rho_0(x)dx + \int_0^T\int_0^L \rho^{(n)}(t,x)\big(\varphi_t(t,x)+v^{(n)}(t,x)\varphi_x(t,x)\big)dxdt$$

$$= \sum_{i=1}^n \int_{x_i(0)}^{x_{i-1}(0)} \varphi(0,x)\rho_0(x)dx + \int_0^T\left(\sum_{i=1}^n \int_{x_i(t)}^{x_{i-1}(t)} \rho^{(n)}(t,x)\big(\varphi_t(t,x)+v^{(n)}(t,x)\varphi_x(t,x)\big)dx\right)dt$$

$$= \sum_{i=1}^n \int_{x_i(0)}^{x_{i-1}(0)} \varphi(0,x)\rho_0(x)dx + \int_0^T\left(\sum_{i=1}^n \rho_i(t)\int_{x_i(t)}^{x_{i-1}(t)} \varphi_t(t,x)dx\right)dt$$

$$+ \int_0^T\left(\sum_{i=1}^n \rho_i(t)\int_{x_i(t)}^{x_{i-1}(t)} v^{(n)}(t,x)\varphi_x(t,x)dx\right)dt + I$$

(4.79)

where

$$I = \int_0^T\left(\sum_{i=1}^n \int_{x_i(t)}^{x_{i-1}(t)} \big(\rho^{(n)}(t,x)-\rho_i(t)\big)\big(\varphi_t(t,x)+v^{(n)}(t,x)\varphi_x(t,x)\big)dx\right)dt \quad (4.80)$$

Using (4.69) and definition (4.80) we get

$$|I| \leq (R_4+1)B\int_0^T\left(\sum_{i=1}^n \int_{x_i(t)}^{x_{i-1}(t)} |\rho^{(n)}(t,x)-\rho_i(t)|dx\right)dt \quad (4.81)$$

where

$$B := \max_{(t,x)\in[0,T]\times[0,L]} \big(|\varphi(t,x)|+|\varphi_t(t,x)|+|\varphi_x(t,x)|\big) \quad (4.82)$$

Combining (4.81) with definition (3.33) we obtain:

$$|I| \leq \frac{B}{2}(R_4+1)\int_0^T\left(\sum_{i=1}^n |\rho_{i-1}(t)-\rho_i(t)|(x_{i-1}(t)-x_i(t))\right)dt \quad (4.83)$$

Using (4.70), (4.83) and the Cauchy-Schwarz inequality we get:

$$|I| \leq \frac{B}{2\sqrt{n}}(R_4+1)\int_0^T\left(n\sum_{i=1}^n |\rho_{i-1}(t)-\rho_i(t)|^2\right)^{1/2}\left(\sum_{i=1}^n (x_{i-1}(t)-x_i(t))^2\right)^{1/2}dt$$

$$\leq \frac{B\sqrt{R_7}}{2\sqrt{n}}(R_4+1)\int_0^T\left(\max_{i=1,\ldots,n}(x_{i-1}(t)-x_i(t))\sum_{i=1}^n (x_{i-1}(t)-x_i(t))\right)^{1/2}dt$$

(4.84)

Using (3.23), (4.84) and the fact that $\sum_{i=1}^n (x_{i-1}(t)-x_i(t)) = L$ we obtain:

$$|I| \leq \frac{BT\sqrt{bLR_7}}{2n}(R_4+1) \quad (4.85)$$

Using (3.33), (4.79), integration by parts and the fact that

$$\frac{d}{dt}\left(\sum_{i=1}^{n}\rho_i(t)\int_{x_i(t)}^{x_{i-1}(t)}\varphi(t,x)dx\right)=\sum_{i=1}^{n}\dot{\rho}_i(t)\int_{x_i(t)}^{x_{i-1}(t)}\varphi(t,x)dx$$

$$+\sum_{i=1}^{n}\rho_i(t)\int_{x_i(t)}^{x_{i-1}(t)}\varphi_t(t,x)dx+\sum_{i=1}^{n}\rho_i(t)\big(\varphi(t,x_{i-1}(t))v_{i-1}(t)-\varphi(t,x_i(t))v_i(t)\big)$$

which is a consequence of (3.1), (3.2), (3.3), we get:

$$\int_0^L\varphi(0,x)\rho_0(x)dx+\int_0^T\int_0^L\rho^{(n)}(t,x)\big(\varphi_t(t,x)+v^{(n)}(t,x)\varphi_x(t,x)\big)dxdt$$

$$=\int_0^T\left(\sum_{i=1}^{n}\rho_i(t)\big(v^{(n)}(t,x_{i-1}(t))\varphi(t,x_{i-1}(t))-v^{(n)}(t,x_i(t))\varphi(t,x_i(t))\big)\right)dt$$

$$-\int_0^T\left(\sum_{i=1}^{n}\rho_i(t)\int_{x_i(t)}^{x_{i-1}(t)}v_x^{(n)}(t,x)\varphi(t,x)dx\right)dt+I+\int_0^T\frac{d}{dt}\left(\sum_{i=1}^{n}\rho_i(t)\int_{x_i(t)}^{x_{i-1}(t)}\varphi(t,x)dx\right)dt \quad (4.86)$$

$$+\sum_{i=1}^{n}\int_{x_i(0)}^{x_{i-1}(0)}\varphi(0,x)\rho_0(x)dx-\int_0^T\left(\sum_{i=1}^{n}\dot{\rho}_i(t)\int_{x_i(t)}^{x_{i-1}(t)}\varphi(t,x)dx\right)dt$$

$$-\int_0^T\left(\sum_{i=1}^{n}\rho_i(t)\big(\varphi(t,x_{i-1}(t))v_{i-1}(t)-\varphi(t,x_i(t))v_i(t)\big)\right)dt$$

Definition (3.33) in conjunction with (4.86) gives:

$$\int_0^L\varphi(0,x)\rho_0(x)dx+\int_0^T\int_0^L\rho^{(n)}(t,x)\big(\varphi_t(t,x)+v^{(n)}(t,x)\varphi_x(t,x)\big)dxdt$$

$$=I+\sum_{i=1}^{n}\rho_i(T)\int_{x_i(T)}^{x_{i-1}(T)}\varphi(T,x)dx+J \quad (4.87)$$

$$-\int_0^T\left(\sum_{i=1}^{n}\left(\dot{\rho}_i(t)+\rho_i(t)\frac{v_{i-1}(t)-v_i(t)}{x_{i-1}(t)-x_i(t)}\right)\int_{x_i(t)}^{x_{i-1}(t)}\varphi(t,x)dx\right)dt$$

where

$$J=\sum_{i=1}^{n}\int_{x_i(0)}^{x_{i-1}(0)}\varphi(0,x)\big(\rho_0(x)-\rho_i(0)\big)dx \quad (4.88)$$

Using (3.31) and the fact that $\int_{x_i(0)}^{x_{i-1}(0)}\rho_0(x)dx=\frac{m}{n}$ for $i=1,\ldots,n$, we conclude that $\rho_i(0)=\frac{1}{x_{i-1}(0)-x_i(0)}\int_{x_i(0)}^{x_{i-1}(0)}\rho_0(x)dx$ for $i=1,\ldots,n$. Since $\rho_0$ is continuous, for each $i=1,\ldots,n$ there exists $\xi_i\in[x_i(0),x_{i-1}(0)]$ such that $\rho_i(0)=\rho_0(\xi_i)$ and consequently, $|\rho_0(x)-\rho_i(0)|\leq\|\rho_0'\|_\infty(x_{i-1}(0)-x_i(0))$ for all $x\in[x_i(0),x_{i-1}(0)]$. Therefore, using (4.88), (3.23), (4.82) and the fact that $\sum_{i=1}^{n}(x_{i-1}(t)-x_i(t))=L$ we get:

$$|J| \leq B\|\rho_0'\|_\infty \sum_{i=1}^{n}(x_{i-1}(0) - x_i(0))^2 \leq \frac{bBL}{n}\|\rho_0'\|_\infty \quad (4.89)$$

The fact that $\varphi(T,x) = 0$ for all $x \in [0,L]$ and equations (3.1), (3.2), (3.3) in conjunction with definition (3.31) (which imply that $\dot{\rho}_i(t) = -\rho_i(t)\frac{v_{i-1}(t) - v_i(t)}{x_{i-1}(t) - x_i(t)}$ for $i = 1,...,n$) allow us to obtain from (4.87):

$$\left| \int_0^L \varphi(0,x)\rho_0(x)dx + \int_0^T\int_0^L \rho^{(n)}(t,x)\left(\varphi_t(t,x) + v^{(n)}(t,x)\varphi_x(t,x)\right)dxdt \right|$$
$$= |I + J| \leq |I| + |J|$$

The above inequality in conjunction with (4.85) and (4.89) imply estimate (4.78) with $R_{10} = bBL\|\rho_0'\|_\infty + \frac{BT\sqrt{bLR_7}}{2}(R_4 + 1)$. The proof is complete. ◁

**Lemma 4.4:** *Suppose that Assumption (A) holds. Let $(\rho_0, v_0) \in X$ that satisfies (3.30) be given. Consider the unique solution of (3.1), (3.2), (3.3) with initial condition $(x_1(0),...,x_{n-1}(0), v_1(0),...,v_{n-1}(0)) = P_n(\rho_0, v_0) \in \Omega_n$ and define $v^{(n)}, \rho^{(n)} : \mathbb{R}_+ \times [0,L] \to \mathbb{R}$ by means of (3.31)-(3.34). Then for every $T > 0$ and for every $\varphi \in C^2([0,T] \times [0,L])$ with $\varphi(T,x) = 0$ for all $x \in [0,L]$ and $\varphi(t,0) = \varphi(t,L) = \varphi_x(t,L) = 0$ for all $t \in [0,T]$ there exists a constant $R_{11} > 0$ (independent of $n \geq 2$ but dependent on $\rho_0, v_0, T, \varphi$) such that the following estimate holds for all $n \geq 2$:*

$$\left| \int_0^L \varphi(0,x)\rho_0(x)v_0(x)dx + \int_0^T\int_0^L \left(\varphi_t(t,x)\rho^{(n)}(t,x)v^{(n)}(t,x) + \varphi_x(t,x)q(t,x)\right)dxdt \right| \leq \frac{R_{11}}{n} \quad (4.90)$$

where $q(t,x) = \rho^{(n)}(t,x)\left(v^{(n)}(t,x)\right)^2 + P(\rho^{(n)}(t,x)) - \mu(\rho^{(n)}(t,x))v_x^{(n)}(t,x)$.

**Remark:** Estimate (4.90) shows how accurately the approximate solution $v^{(n)}, \rho^{(n)}$ satisfies the momentum equation. Notice that the accuracy provided by the particle scheme is of order $1/n$.

**Proof:** Let $T > 0$ (arbitrary) be given and let $\varphi \in C^1([0,T] \times [0,L])$ with $\varphi(T,x) = 0$ for all $x \in [0,L]$ and $\varphi(t,0) = \varphi(t,L) = \varphi_x(t,L) = 0$ for all $t \in [0,T]$ (but otherwise arbitrary) be given. We get

$$\int_0^L \varphi(0,x)\rho_0(x)v_0(x)dx + \int_0^T\int_0^L \left(\varphi_t(t,x)\rho^{(n)}(t,x)v^{(n)}(t,x) + \varphi_x(t,x)q(t,x)\right)dxdt$$
$$= \sum_{i=1}^{n} \int_{x_i(0)}^{x_{i-1}(0)} \varphi(0,x)\rho_0(x)v_0(x)dx$$
$$+ \int_0^T \left( \sum_{i=1}^{n} \int_{x_i(t)}^{x_{i-1}(t)} \left(\varphi_t(t,x)\rho_i(t)v_{i-1}(t) + \varphi_x(t,x)q(t,x)\right)dx \right)dt \quad (4.91)$$
$$+ \int_0^T \left( \sum_{i=1}^{n} \int_{x_i(t)}^{x_{i-1}(t)} \varphi_t(t,x)\left(\rho^{(n)}(t,x)v^{(n)}(t,x) - \rho_i(t)v_{i-1}(t)\right)dx \right)dt$$

Using (3.1), (3.2), (3.3), (3.5) and definitions (3.31) we get for all $i = 2,...,n$ and $t \geq 0$:

$$\frac{d}{dt}\left(\rho_i(t)v_{i-1}(t)\right) = -\rho_i(t)v_{i-1}(t)\frac{v_{i-1}(t)-v_i(t)}{x_{i-1}(t)-x_i(t)} + \frac{n}{m}\rho_i(t)\left(P(\rho_i(t))-P(\rho_{i-1}(t))\right)$$
$$+\frac{n}{m}\rho_i(t)\left(\mu(\rho_{i-1}(t))\frac{v_{i-2}(t)-v_{i-1}(t)}{x_{i-2}(t)-x_{i-1}(t)} - \mu(\rho_i(t))\frac{v_{i-1}(t)-v_i(t)}{x_{i-1}(t)-x_i(t)}\right) \quad (4.92)$$

Multiplying (4.92) by $\varphi(t,x)$, integrating and using (3.33) we get for all $i = 2,...,n$ and $t \geq 0$:

$$\int_{x_i(t)}^{x_{i-1}(t)}\frac{d}{dt}\left(\rho_i(t)v_{i-1}(t)\right)\varphi(t,x)dx = -\int_{x_i(t)}^{x_{i-1}(t)}\rho_i(t)v_{i-1}(t)v_x^{(n)}(t,x)\varphi(t,x)dx$$
$$+\frac{n}{m}\int_{x_i(t)}^{x_{i-1}(t)}\rho_i(t)\left(P(\rho_i(t))-P(\rho_{i-1}(t))\right)\varphi(t,x)dx \quad (4.93)$$
$$+\frac{n}{m}\int_{x_i(t)}^{x_{i-1}(t)}\rho_i(t)\left(\mu(\rho_{i-1}(t))\frac{v_{i-2}(t)-v_{i-1}(t)}{x_{i-2}(t)-x_{i-1}(t)} - \mu(\rho_i(t))\frac{v_{i-1}(t)-v_i(t)}{x_{i-1}(t)-x_i(t)}\right)\varphi(t,x)dx$$

Using (3.1), (3.2), (3.3) we also obtain for all $i = 2,...,n$ and $t \geq 0$

$$\frac{d}{dt}\int_{x_i(t)}^{x_{i-1}(t)}\rho_i(t)v_{i-1}(t)\varphi(t,x)dx = \int_{x_i(t)}^{x_{i-1}(t)}\frac{d}{dt}\left(\rho_i(t)v_{i-1}(t)\right)\varphi(t,x)dx$$
$$+\int_{x_i(t)}^{x_{i-1}(t)}\rho_i(t)v_{i-1}(t)\varphi_t(t,x)dx + \rho_i(t)v_{i-1}(t)\left(\varphi(t,x_{i-1}(t))v_{i-1}(t)-\varphi(t,x_i(t))v_i(t)\right)$$

which combined with (4.93) gives for all $i = 2,...,n$ and $t \geq 0$:

$$\frac{d}{dt}\int_{x_i(t)}^{x_{i-1}(t)}\rho_i(t)v_{i-1}(t)\varphi(t,x)dx = \int_{x_i(t)}^{x_{i-1}(t)}\rho_i(t)v_{i-1}(t)\varphi_t(t,x)dx$$
$$+\rho_i(t)v_{i-1}(t)\left(\varphi(t,x_{i-1}(t))v_{i-1}(t)-\varphi(t,x_i(t))v_i(t) - \int_{x_i(t)}^{x_{i-1}(t)}v_x^{(n)}(t,x)\varphi(t,x)dx\right)$$
$$-\frac{n}{m}\int_{x_i(t)}^{x_{i-1}(t)}\rho_i(t)P'(\rho^{(n)}(t,x))\left(\rho_{i-1}(t)-\rho_i(t)\right)\varphi(t,x)dx + \frac{n}{m}\int_{x_i(t)}^{x_{i-1}(t)}\rho_i(t)g_i(t,x)\varphi(t,x)dx$$
$$+\frac{n}{m}\int_{x_i(t)}^{x_{i-1}(t)}\rho_i(t)\left(\mu(\rho_{i-1}(t))\frac{v_{i-2}(t)-v_{i-1}(t)}{x_{i-2}(t)-x_{i-1}(t)} - \mu(\rho_i(t))\frac{v_{i-1}(t)-v_i(t)}{x_{i-1}(t)-x_i(t)}\right)\varphi(t,x)dx$$

(4.94)

where

$$g_i(t,x) = P(\rho_i(t)) - P(\rho_{i-1}(t)) - P'(\rho^{(n)}(t,x))\left(\rho_i(t)-\rho_{i-1}(t)\right) \quad (4.95)$$

Using (3.31), (3.33) we obtain from (4.94) for all $i = 2,...,n$ and $t \geq 0$:

$$\frac{d}{dt}\int_{x_i(t)}^{x_{i-1}(t)} \rho_i(t)v_{i-1}(t)\varphi(t,x)dx = \int_{x_i(t)}^{x_{i-1}(t)} \rho_i(t)v_{i-1}(t)\left(\varphi_t(t,x)+v^{(n)}(t,x)\varphi_x(t,x)\right)dx$$

$$-\int_{x_i(t)}^{x_{i-1}(t)}\left(P'\left(\rho^{(n)}(t,x)\right)-\mu'\left(\rho^{(n)}(t,x)\right)\frac{v_{i-1}(t)-v_i(t)}{x_{i-1}(t)-x_i(t)}\right)\rho_x^{(n)}(t,x)\varphi(t,x)dx$$

$$+\frac{\mu(\rho_{i-1}(t))}{x_{i-1}(t)-x_i(t)}\left(\frac{v_{i-2}(t)-v_{i-1}(t)}{x_{i-2}(t)-x_{i-1}(t)}-\frac{v_{i-1}(t)-v_i(t)}{x_{i-1}(t)-x_i(t)}\right)\int_{x_i(t)}^{x_{i-1}(t)}\varphi(t,x)dx$$

$$+\frac{1}{x_{i-1}(t)-x_i(t)}\int_{x_i(t)}^{x_{i-1}(t)} g_i(t,x)\varphi(t,x)dx + \frac{v_{i-1}(t)-v_i(t)}{\left(x_{i-1}(t)-x_i(t)\right)^2}\int_{x_i(t)}^{x_{i-1}(t)} f_i(t,x)\varphi(t,x)dx$$

(4.96)

where

$$f_i(t,x)=\mu(\rho_{i-1}(t))-\mu(\rho_i(t))-\mu'\left(\rho^{(n)}(t,x)\right)(\rho_{i-1}(t)-\rho_i(t)) \quad (4.97)$$

Using (3.31), (3.33), integrating by parts and using the fact that $q(t,x)=\rho^{(n)}(t,x)\left(v^{(n)}(t,x)\right)^2+P(\rho^{(n)}(t,x))-\mu(\rho^{(n)}(t,x))v_x^{(n)}(t,x)$, we obtain from (4.96) for all $i=2,\ldots,n$ and $t\geq 0$:

$$\frac{d}{dt}\int_{x_i(t)}^{x_{i-1}(t)} \rho_i(t)v_{i-1}(t)\varphi(t,x)dx = \int_{x_i(t)}^{x_{i-1}(t)}\left(\rho_i(t)v_{i-1}(t)\varphi_t(t,x)+q(t,x)\varphi_x(t,x)\right)dx$$

$$-P(\rho_{i-1}(t))\varphi(t,x_{i-1}(t))+\mu(\rho_{i-1}(t))\frac{v_{i-2}(t)-v_{i-1}(t)}{x_{i-2}(t)-x_{i-1}(t)}\varphi(t,x_{i-1}(t))$$

$$+P(\rho_i(t))\varphi(t,x_i(t))-\mu(\rho_i(t))\frac{v_{i-1}(t)-v_i(t)}{x_{i-1}(t)-x_i(t)}\varphi(t,x_i(t))$$

$$+\left(\frac{v_{i-2}(t)-v_{i-1}(t)}{x_{i-2}(t)-x_{i-1}(t)}-\frac{v_{i-1}(t)-v_i(t)}{x_{i-1}(t)-x_i(t)}\right)\frac{\mu(\rho_{i-1}(t))}{x_{i-1}(t)-x_i(t)}\int_{x_i(t)}^{x_{i-1}(t)}\left(\varphi(t,x)-\varphi(t,x_{i-1}(t))\right)dx$$

$$+\frac{1}{x_{i-1}(t)-x_i(t)}\int_{x_i(t)}^{x_{i-1}(t)} g_i(t,x)\varphi(t,x)dx + \frac{v_{i-1}(t)-v_i(t)}{\left(x_{i-1}(t)-x_i(t)\right)^2}\int_{x_i(t)}^{x_{i-1}(t)} f_i(t,x)\varphi(t,x)dx$$

$$+\int_{x_i(t)}^{x_{i-1}(t)}\left(\rho_i(t)v_{i-1}(t)-\rho^{(n)}(t,x)v^{(n)}(t,x)\right)v^{(n)}(t,x)\varphi_x(t,x)dx$$

(4.98)

Summing up all equations (4.98) for $i=2,\ldots,n$, using (3.1), (3.3) the facts that $x_0(t)\equiv L$, $x_n(t)\equiv 0$, $\varphi(T,x)=0$ for all $x\in[0,L]$, $\varphi(t,0)=\varphi(t,L)=0$ for all $t\in[0,T]$ and integrating we get:

$$-\sum_{i=1}^{n}\int_{x_i(0)}^{x_{i-1}(0)}\rho_i(0)v_{i-1}(0)\varphi(0,x)dx = \int_0^T\left(\sum_{i=1}^{n}\int_{x_i(t)}^{x_{i-1}(t)}\left(\rho_i(t)v_{i-1}(t)\varphi_t(t,x)+q(t,x)\varphi_x(t,x)\right)dx\right)dt$$

$$-\int_0^T\left(\int_{x_1(t)}^{L}q(t,x)\varphi_x(t,x)dx\right)dt - \int_0^T P(\rho_1(t))\varphi(t,x_1(t))dt - \int_0^T\frac{\mu(\rho_1(t))v_1(t)}{L-x_1(t)}\varphi(t,x_1(t))dt$$

$$+\int_0^T\left(\sum_{i=2}^{n}\left(\frac{v_{i-2}(t)-v_{i-1}(t)}{x_{i-2}(t)-x_{i-1}(t)}-\frac{v_{i-1}(t)-v_i(t)}{x_{i-1}(t)-x_i(t)}\right)\frac{\mu(\rho_{i-1}(t))}{x_{i-1}(t)-x_i(t)}\int_{x_i(t)}^{x_{i-1}(t)}\left(\varphi(t,x)-\varphi(t,x_{i-1}(t))\right)dx\right)dt$$

$$+\int_0^T\left(\sum_{i=2}^{n}\left(\frac{1}{x_{i-1}(t)-x_i(t)}\int_{x_i(t)}^{x_{i-1}(t)}g_i(t,x)\varphi(t,x)dx+\frac{v_{i-1}(t)-v_i(t)}{(x_{i-1}(t)-x_i(t))^2}\int_{x_i(t)}^{x_{i-1}(t)}f_i(t,x)\varphi(t,x)dx\right)\right)dt$$

$$+\int_0^T\left(\sum_{i=2}^{n}\int_{x_i(t)}^{x_{i-1}(t)}\left(\rho_i(t)v_{i-1}(t)-\rho^{(n)}(t,x)v^{(n)}(t,x)\right)v^{(n)}(t,x)\varphi_x(t,x)dx\right)dt$$

The above equation combined with (4.91) gives the following equation:

$$\int_0^L\varphi(0,x)\rho_0(x)v_0(x)dx + \int_0^T\int_0^L\left(\varphi_t(t,x)\rho^{(n)}(t,x)v^{(n)}(t,x)+\varphi_x(t,x)q(t,x)\right)dxdt \quad (4.99)$$
$$= I_1 - I_2 + I_3 - I_4 + I_5$$

where

$$I_1 = \int_0^T\left(\sum_{i=2}^{n}\int_{x_i(t)}^{x_{i-1}(t)}\left(\rho^{(n)}(t,x)v^{(n)}(t,x)-\rho_i(t)v_{i-1}(t)\right)v^{(n)}(t,x)\varphi_x(t,x)dx\right)dt$$
$$+\int_0^T\left(\sum_{i=1}^{n}\int_{x_i(t)}^{x_{i-1}(t)}\varphi_t(t,x)\left(\rho^{(n)}(t,x)v^{(n)}(t,x)-\rho_i(t)v_{i-1}(t)\right)dx\right)dt \quad (4.100)$$

$$I_2 = \int_0^T\left(\sum_{i=2}^{n}\frac{1}{x_{i-1}(t)-x_i(t)}\int_{x_i(t)}^{x_{i-1}(t)}g_i(t,x)\varphi(t,x)dx\right)dt$$
$$+\int_0^T\left(\sum_{i=2}^{n}\frac{v_{i-1}(t)-v_i(t)}{(x_{i-1}(t)-x_i(t))^2}\int_{x_i(t)}^{x_{i-1}(t)}f_i(t,x)\varphi(t,x)dx\right)dt \quad (4.101)$$

$$I_3 = \int_0^T\left(\int_{x_1(t)}^{L}q(t,x)\varphi_x(t,x)dx\right)dt$$
$$+\int_0^T P(\rho_1(t))\varphi(t,x_1(t))dt + \int_0^T\frac{\mu(\rho_1(t))v_1(t)}{L-x_1(t)}\varphi(t,x_1(t))dt \quad (4.102)$$

$$I_4 =$$
$$\int_0^T\left(\sum_{i=2}^{n}\left(\frac{v_{i-2}(t)-v_{i-1}(t)}{x_{i-2}(t)-x_{i-1}(t)}-\frac{v_{i-1}(t)-v_i(t)}{x_{i-1}(t)-x_i(t)}\right)\frac{\mu(\rho_{i-1}(t))}{x_{i-1}(t)-x_i(t)}\int_{x_i(t)}^{x_{i-1}(t)}\left(\varphi(t,x)-\varphi(t,x_{i-1}(t))\right)dx\right)dt \quad (4.103)$$

$$I_5 = \sum_{i=1}^{n}\int_{x_i(0)}^{x_{i-1}(0)}\varphi(0,x)\left(\rho_0(x)v_0(x)-\rho_i(0)v_{i-1}(0)\right)dx \quad (4.104)$$

We next estimate each one of the terms $I_i$, $i=1,\ldots,5$. We start from the term $I_1$. Using (4.69) and (4.100) we get:

$$|I_1| \leq B(1+R_4) \int_0^T \left( \sum_{i=1}^n \int_{x_i(t)}^{x_{i-1}(t)} |v^{(n)}(t,x)| |\rho^{(n)}(t,x) - \rho_i(t,x)| dx \right) dt$$
$$+ B(1+R_4) \int_0^T \left( \sum_{i=1}^n \int_{x_i(t)}^{x_{i-1}(t)} \rho_i(t,x) |v^{(n)}(t,x) - v_{i-1}(t)| dx \right) dt \quad (4.105)$$

where

$$B := \max_{(t,x) \in [0,T] \times [0,L]} \left( |\varphi(t,x)| + |\varphi_t(t,x)| + |\varphi_x(t,x)| + |\varphi_{xx}(t,x)| \right) \quad (4.106)$$

Using (3.23), (3.31), (4.69) and (4.105) we get:

$$|I_1| \leq B(1+R_4) R_4 \int_0^T \left( \sum_{i=1}^n \int_{x_i(t)}^{x_{i-1}(t)} |\rho^{(n)}(t,x) - \rho_i(t,x)| dx \right) dt$$
$$+ B(1+R_4) \frac{m}{a} \int_0^T \left( \sum_{i=1}^n \int_{x_i(t)}^{x_{i-1}(t)} |v^{(n)}(t,x) - v_{i-1}(t)| dx \right) dt \quad (4.107)$$

Using (3.33), (4.107) and the Cauchy-Schwarz inequality we get:

$$|I_1| \leq \frac{1}{2} B(1+R_4) R_4 \int_0^T \left( \sum_{i=1}^n |\rho_{i-1}(t) - \rho_i(t)| (x_{i-1}(t) - x_i(t)) \right) dt$$
$$+ B(1+R_4) \frac{m}{2a} \int_0^T \left( \sum_{i=1}^n |v_{i-1}(t) - v_i(t)| (x_{i-1}(t) - x_i(t)) \right) dt$$
$$\leq \frac{1}{2} B(1+R_4) R_4 \int_0^T \left( \sum_{i=1}^n |\rho_{i-1}(t) - \rho_i(t)|^2 \right)^{1/2} \left( \sum_{i=1}^n (x_{i-1}(t) - x_i(t))^2 \right)^{1/2} dt \quad (4.108)$$
$$+ B(1+R_4) \frac{m}{2a} \int_0^T \left( \sum_{i=1}^n \frac{|v_{i-1}(t) - v_i(t)|^2}{x_{i-1}(t) - x_i(t)} \right)^{1/2} \left( \sum_{i=1}^n (x_{i-1}(t) - x_i(t))^3 \right)^{1/2} dt$$

Using (3.12), (3.24), (4.70) and (4.108) we get:

$$|I_1| \leq \frac{1}{2\sqrt{n}} B(1+R_4) R_4 \sqrt{R_7} \int_0^T \left( \max_{i=1,\ldots,n} (x_{i-1}(t) - x_i(t)) \sum_{i=1}^n (x_{i-1}(t) - x_i(t)) \right)^{1/2} dt$$
$$+ B(1+R_4) \frac{m\sqrt{R_2}}{a\sqrt{2}} \int_0^T \left( \max_{i=1,\ldots,n} (x_{i-1}(t) - x_i(t)) \right) \left( \sum_{i=1}^n (x_{i-1}(t) - x_i(t)) \right)^{1/2} dt \quad (4.109)$$

Using the fact that $\sum_{i=1}^n (x_{i-1}(t) - x_i(t)) = L$ and (3.23), (4.109) we get:

$$|I_1| \leq \frac{BT}{2n} (1+R_4) \sqrt{bL} \left( R_4 \sqrt{R_7} + \frac{m\sqrt{2R_2}}{a\sqrt{b}} \right) \quad (4.110)$$

We continue with the term $I_2$. Definitions (4.95), (4.97) in conjunction with (3.33), (4.68) imply that for every $t \geq 0$, $i = 2,\ldots,n$ and $x \in [x_i(t), x_{i-1}(t)]$ we get:

$$|g_i(t,x)| + |f_i(t,x)| \leq \Gamma (\rho_i(t) - \rho_{i-1}(t))^2 \tag{4.111}$$

where $\Gamma = \max\left\{|P''(s)|: \dfrac{m}{b} \leq s \leq \dfrac{m}{a}\right\} + \max\left\{|\mu''(s)|: \dfrac{m}{b} \leq s \leq \dfrac{m}{a}\right\}$. It follows from (4.101), (4.106) and (4.111) that:

$$\begin{aligned}|I_2| &\leq \Gamma B \int_0^T \left(\sum_{i=1}^n (\rho_{i-1}(t) - \rho_i(t))^2\right) dt \\ &+ \Gamma B \int_0^T \left(\sum_{i=1}^n \frac{(v_{i-1}(t) - v_i(t))(\rho_{i-1}(t) - \rho_i(t))^2}{x_{i-1}(t) - x_i(t)}\right) dt\end{aligned} \tag{4.112}$$

We conclude from (4.71), (4.74) and (4.112) that the following estimate holds:

$$\begin{aligned}|I_2| &\leq \Gamma B \frac{R_8^2}{a^2} \exp(2\kappa T) \int_0^T \left(\sum_{i=1}^n (x_{i-1}(t) - x_i(t))^2\right) dt \\ &+ \Gamma B \frac{R_8^2}{a^2} \exp(2\kappa T) \int_0^T \left(\sum_{i=1}^n (v_{i-1}(t) - v_i(t))(x_{i-1}(t) - x_i(t))\right) dt\end{aligned} \tag{4.113}$$

Working as in the case of the term $I_1$ we get from (4.113):

$$|I_2| \leq \Gamma B \frac{b R_8^2}{n a^2} \left(L + \sqrt{2 L R_2}\right) T \exp(2\kappa T) \tag{4.114}$$

We next continue with the term $I_3$. Since $\varphi(t,L) = \varphi_x(t,L) = 0$ for all $t \in [0,T]$ and since $x_0(t) \equiv L$, we get from (4.106):

$$|\varphi(t, x_1(t))| \leq \frac{B}{2}(L - x_1(t))^2, \quad |\varphi_x(t,x)| \leq B(L - x_1(t)) \text{ for all } t \in [0,T] \text{ and } x \in [x_1(t), 0]$$

Consequently, we obtain from definition (4.102) and (3.23), (3.31):

$$\begin{aligned}|I_3| &\leq B \int_0^T \left((L - x_1(t)) \int_{x_1(t)}^L |q(t,x)| dx\right) dt \\ &+ \frac{B}{2} \max\left\{P(s): \frac{m}{b} \leq s \leq \frac{m}{a}\right\} \int_0^T (L - x_1(t))^2 dt \\ &+ \frac{B}{2} \max\left\{\mu(s): \frac{m}{b} \leq s \leq \frac{m}{a}\right\} \int_0^T |v_1(t)|(L - x_1(t)) dt\end{aligned} \tag{4.115}$$

Since $q(t,x) = \rho^{(n)}(t,x)(v^{(n)}(t,x))^2 + P(\rho^{(n)}(t,x)) - \mu(\rho^{(n)}(t,x)) v_x^{(n)}(t,x)$, we get from (3.1), (3.33), (4.68), (4.69) for all $t \in [0,T]$ and $x \in [x_1(t), 0]$:

$$|q(t,x)| \leq \frac{m}{a} R_4^2 + \max\left\{P(s): \frac{m}{b} \leq s \leq \frac{m}{a}\right\}$$

$$+ \max\left\{\mu(s): \frac{m}{b} \leq s \leq \frac{m}{a}\right\} \frac{|v_1(t)|}{L - x_1(t)}$$

Combining the above estimate with (4.115) we get:

$$|I_3| \leq B\left(\frac{m}{a} R_4^2 + \frac{3}{2}\max\left\{P(s): \frac{m}{b} \leq s \leq \frac{m}{a}\right\}\right) \int_0^T (L - x_1(t))^2 \, dt$$

$$+ \frac{3B}{2} \max\left\{\mu(s): \frac{m}{b} \leq s \leq \frac{m}{a}\right\} \int_0^T |v_1(t)|(L - x_1(t)) \, dt \qquad (4.116)$$

Using (3.23), (3.26), the facts that $n \geq 2$, $x_0(t) \equiv L$ and (4.116) we get:

$$|I_3| \leq \frac{BbT}{n}\left(b\left(\frac{m}{a} R_4^2 + \frac{3}{2}\max\left\{P(s): \frac{m}{b} \leq s \leq \frac{m}{a}\right\}\right) + \frac{3}{2}\max\left\{\mu(s): \frac{m}{b} \leq s \leq \frac{m}{a}\right\} R_4\right) \qquad (4.117)$$

We next continue with the term $I_4$. Using (4.106), (3.23), (3.31) we get

$$\frac{\mu(\rho_{i-1}(t))}{x_{i-1}(t) - x_i(t)} \int_{x_i(t)}^{x_{i-1}(t)} |\varphi(t,x) - \varphi(t, x_{i-1}(t))| \, dx \leq \frac{B}{2}(x_{i-1}(t) - x_i(t)) \max\left\{\mu(s): \frac{m}{b} \leq s \leq \frac{m}{a}\right\}$$

which combined with definition (4.103) gives:

$$|I_4| \leq \frac{BM}{2} \int_0^T \left(\sum_{i=2}^n \left|\frac{v_{i-2}(t) - v_{i-1}(t)}{x_{i-2}(t) - x_{i-1}(t)} - \frac{v_{i-1}(t) - v_i(t)}{x_{i-1}(t) - x_i(t)}\right|(x_{i-1}(t) - x_i(t))\right) dt \qquad (4.118)$$

where $M = \max\left\{\mu(s): \frac{m}{b} \leq s \leq \frac{m}{a}\right\}$. Applying the Cauchy-Schwarz inequality in (4.119), using (3.23) and the fact that $\sum_{i=1}^n (x_{i-1}(t) - x_i(t)) = L$ as well as the inequality

$$\frac{1}{\sqrt{n}}\left(\sum_{i=2}^n \left|\frac{v_{i-2}(t) - v_{i-1}(t)}{x_{i-2}(t) - x_{i-1}(t)} - \frac{v_{i-1}(t) - v_i(t)}{x_{i-1}(t) - x_i(t)}\right|^2\right)^{1/2}$$

$$\leq \frac{1}{2n} + \frac{1}{2}\sum_{i=2}^n \left|\frac{v_{i-2}(t) - v_{i-1}(t)}{x_{i-2}(t) - x_{i-1}(t)} - \frac{v_{i-1}(t) - v_i(t)}{x_{i-1}(t) - x_i(t)}\right|^2$$

we get:

$$|I_4| \leq \frac{BM}{2}\int_0^T \left(\sum_{i=2}^n \left|\frac{v_{i-2}(t)-v_{i-1}(t)}{x_{i-2}(t)-x_{i-1}(t)} - \frac{v_{i-1}(t)-v_i(t)}{x_{i-1}(t)-x_i(t)}\right|^2\right)^{1/2} \left(\sum_{i=2}^n (x_{i-1}(t)-x_i(t))^2\right)^{1/2} dt$$

$$\leq \frac{BM\sqrt{bL}}{2}\int_0^T \frac{1}{\sqrt{n}} \left(\sum_{i=2}^n \left|\frac{v_{i-2}(t)-v_{i-1}(t)}{x_{i-2}(t)-x_{i-1}(t)} - \frac{v_{i-1}(t)-v_i(t)}{x_{i-1}(t)-x_i(t)}\right|^2\right)^{1/2} dt$$

$$\leq \frac{BM\sqrt{bL}}{4}\int_0^T \left(\frac{1}{n} + \sum_{i=2}^n \left|\frac{v_{i-2}(t)-v_{i-1}(t)}{x_{i-2}(t)-x_{i-1}(t)} - \frac{v_{i-1}(t)-v_i(t)}{x_{i-1}(t)-x_i(t)}\right|^2\right) dt$$

$$\leq \frac{BM\sqrt{bL}}{4}\left(\frac{T}{n} + \int_0^T \sum_{i=2}^n \left|\frac{v_{i-2}(t)-v_{i-1}(t)}{x_{i-2}(t)-x_{i-1}(t)} - \frac{v_{i-1}(t)-v_i(t)}{x_{i-1}(t)-x_i(t)}\right|^2 dt\right)$$

Combining the above estimate with (3.27) we obtain:

$$|I_4| \leq \frac{BM\sqrt{bL}}{4n}(R_5 + (1+R_6)T) \tag{4.119}$$

Finally, we estimate the term $I_5$. Using (4.106) and definition (4.104) we get:

$$|I_5| \leq \sum_{i=1}^n \int_{x_i(0)}^{x_{i-1}(0)} |\varphi(0,x)||\rho_0(x)v_0(x) - \rho_i(0)v_{i-1}(0)| dx$$

$$\leq \sum_{i=1}^n \int_{x_i(0)}^{x_{i-1}(0)} |\varphi(0,x)||v_0(x)||\rho_0(x) - \rho_i(0)| dx + \sum_{i=1}^n \rho_i(0) \int_{x_i(0)}^{x_{i-1}(0)} |\varphi(0,x)||v_0(x) - v_{i-1}(0)| dx \tag{4.120}$$

$$\leq B\|v_0\|_\infty \sum_{i=1}^n \int_{x_i(0)}^{x_{i-1}(0)} |\rho_0(x) - \rho_i(0)| dx + B\sum_{i=1}^n \rho_i(0) \int_{x_i(0)}^{x_{i-1}(0)} |v_0(x) - v_{i-1}(0)| dx$$

Using (3.31) and the fact that $\int_{x_i(0)}^{x_{i-1}(0)} \rho_0(x)dx = \frac{m}{n}$ for $i=1,...,n$, we conclude that $\rho_i(0) = \frac{1}{x_{i-1}(0)-x_i(0)} \int_{x_i(0)}^{x_{i-1}(0)} \rho_0(x)dx$ for $i=1,...,n$. Since $\rho_0$ is continuous, for each $i=1,...,n$ there exists $\xi_i \in [x_i(0), x_{i-1}(0)]$ such that $\rho_i(0) = \rho_0(\xi_i)$ and consequently, $|\rho_0(x) - \rho_i(0)| \leq \|\rho_0'\|_\infty (x_{i-1}(0) - x_i(0))$ for all $x \in [x_i(0), x_{i-1}(0)]$. Therefore, we get from (4.120):

$$|I_5| \leq B\|v_0\|_\infty \|\rho_0'\|_\infty \sum_{i=1}^n (x_{i-1}(0) - x_i(0))^2$$
$$+ \frac{mB}{n} \sum_{i=1}^n \frac{1}{x_{i-1}(0)-x_i(0)} \int_{x_i(0)}^{x_{i-1}(0)} |v_0(x) - v_{i-1}(0)| dx \tag{4.121}$$

Since $v_i(0) = v(x_i(0))$ for $i=0,1,...,n$, we get using Cauchy-Schwarz inequality for all $x \in [x_i(0), x_{i-1}(0)]$:

$$\left|v_{i-1}(0)-v_0(x)\right|=\left|v_0(x_{i-1}(0))-v_0(x)\right|\leq \int_{x}^{x_{i-1}(0)}\left|v_0'(s)\right|ds$$

$$\leq \left(x_{i-1}(0)-x_i(0)\right)^{1/2}\left(\int_{x_i(0)}^{x_{i-1}(0)}\left(v_0'(s)\right)^2 ds\right)^{1/2}$$

Thus, we get from (4.121) using the Cauchy-Schwarz inequality:

$$\begin{aligned}\left|I_5\right|&\leq B\left\|v_0\right\|_\infty \left\|\rho_0'\right\|_\infty \max_{i=1,\ldots,n}\left(x_{i-1}(0)-x_i(0)\right)\sum_{i=1}^{n}\left(x_{i-1}(0)-x_i(0)\right)\\ &+\frac{mB}{n}\sum_{i=1}^{n}\left(x_{i-1}(0)-x_i(0)\right)^{1/2}\left(\int_{x_i(0)}^{x_{i-1}(0)}\left(v_0'(s)\right)^2 ds\right)^{1/2}\\ &\leq B\left\|v_0\right\|_\infty \left\|\rho_0'\right\|_\infty \max_{i=1,\ldots,n}\left(x_{i-1}(0)-x_i(0)\right)\sum_{i=1}^{n}\left(x_{i-1}(0)-x_i(0)\right)\\ &+\frac{mB}{n}\left(\sum_{i=1}^{n}\left(x_{i-1}(0)-x_i(0)\right)\right)^{1/2}\left(\sum_{i=1}^{n}\int_{x_i(0)}^{x_{i-1}(0)}\left(v_0'(s)\right)^2 ds\right)^{1/2}\end{aligned} \quad (4.122)$$

Since $\sum_{i=1}^{n}\left(x_{i-1}(0)-x_i(0)\right)=L$ we get from (4.122) and (3.23):

$$\left|I_5\right|\leq \frac{B}{n}\left(b\left\|v_0\right\|_\infty \left\|\rho_0'\right\|_\infty L+m\sqrt{L}\left\|v_0'\right\|_2\right) \quad (4.123)$$

Estimate (4.90) is a direct consequence of (4.99), (4.110), (4.114), (4.117), (4.119) and (4.123). The proof is complete. ◁

**Lemma 4.5:** *Suppose that Assumption (A) holds. Let $(\rho_0,v_0)\in X$ that satisfies (3.30) be given. Consider the unique solution of (3.1), (3.2), (3.3) with initial condition $(x_1(0),\ldots,x_{n-1}(0),v_1(0),\ldots,v_{n-1}(0))=P_n(\rho_0,v_0)\in \Omega_n$ and define $v^{(n)},\rho^{(n)}:\mathbb{R}_+\times[0,L]\to\mathbb{R}$ by means of (3.31)-(3.34). Then there exists a constant $R_{12}>0$ (independent of $n\geq 2$ but dependent on $\rho_0,v_0$) such that the following estimate holds for all $n\geq 2$ and $t\geq 0$:*

$$\left|\int_0^L \rho^{(n)}(t,x)dx-m\right|\leq \frac{R_{12}}{n} \quad (4.124)$$

**Proof:** Definition (3.31) implies $m=\sum_{i=1}^{n}\int_{x_i(t)}^{x_{i-1}(t)}\rho_i(t)dx$. Using (3.33), (4.70), (3.23), the fact that $\sum_{i=1}^{n}\left(x_{i-1}(t)-x_i(t)\right)=L$ and the Cauchy-Schwarz inequality we get:

$$\left|\int_0^L \rho^{(n)}(t,x)dx - m\right| = \sum_{i=1}^n \int_{x_i(t)}^{x_{i-1}(t)} \left|\rho^{(n)}(t,x) - \rho_i(t)\right|dx$$

$$\leq \frac{1}{2}\sum_{i=1}^n |\rho_{i-1}(t) - \rho_i(t)|(x_{i-1}(t) - x_i(t))$$

$$\leq \frac{1}{2}\left(\sum_{i=1}^n |\rho_{i-1}(t) - \rho_i(t)|^2\right)^{1/2}\left(\sum_{i=1}^n (x_{i-1}(t) - x_i(t))^2\right)^{1/2}$$

$$\leq \frac{\sqrt{R_7}}{2\sqrt{n}}\left(\max_{i=1,\ldots,n}(x_{i-1}(t) - x_i(t))\sum_{i=1}^n (x_{i-1}(t) - x_i(t))\right)^{1/2} \leq \frac{\sqrt{LbR_7}}{2n}$$

The above inequality implies (4.124). The proof is complete. ◁

**Lemma 4.6:** *Suppose that Assumption (A) holds. Let $(\rho_0, v_0) \in X$ that satisfies (3.30) be given. Consider the unique solution of (3.1), (3.2), (3.3) with initial condition $(x_1(0),\ldots,x_{n-1}(0), v_1(0),\ldots,v_{n-1}(0)) = P_n(\rho_0, v_0) \in \Omega_n$ and define $v^{(n)}, \rho^{(n)} : \mathbb{R}_+ \times [0,L] \to \mathbb{R}$ by means of (3.31)-(3.34). Then there exists a constant $R_{13} > 0$ (independent of $n \geq 2$ but dependent on $\rho_0, v_0$) such that the following estimate holds for all $n \geq 2$ and $t \geq t_0 \geq 0$:*

$$E(\rho^{(n)}[t], v^{(n)}[t]) \leq E(\rho^{(n)}[t_0], v^{(n)}[t_0]) + \frac{R_{13}}{n} \tag{4.125}$$

*where the functional $E(\rho, v)$ is defined by (2.9).*

**Proof:** Since (4.1) holds (which implies that $E_n(t) \leq E_n(t_0)$ for all $n \geq 2$ and $t \geq t_0 \geq 0$), it suffices to show that there exists a constant $R_{14} > 0$ (independent of $n \geq 2$ but dependent on $\rho_0, v_0$) such that the following estimate holds for all $n \geq 2$ and $t \geq 0$:

$$\left|E(\rho^{(n)}[t], v^{(n)}[t]) - E_n(t)\right| \leq \frac{R_{14}}{n} \tag{4.126}$$

Using definitions (2.9), (2.11), (3.33) we get for all $t \geq 0$:

$$E(\rho^{(n)}[t], v^{(n)}[t]) = \frac{1}{2}\sum_{i=1}^n \int_{x_i(t)}^{x_{i-1}(t)} \rho^{(n)}(t,x)(v^{(n)}(t,x))^2 dx + \sum_{i=1}^n \int_{x_i(t)}^{x_{i-1}(t)} Q(\rho^{(n)}(t,x))dx$$

$$= \frac{1}{2}\sum_{i=1}^n \int_{x_i(t)}^{x_{i-1}(t)} (\rho^{(n)}(t,x) - \rho_i(t))(v^{(n)}(t,x))^2 dx + \sum_{i=1}^n \int_{x_i(t)}^{x_{i-1}(t)} (Q(\rho^{(n)}(t,x)) - Q(\rho_i(t)))dx$$

$$+ \frac{1}{2}\sum_{i=1}^n \int_{x_i(t)}^{x_{i-1}(t)} \rho_i(t)(v^{(n)}(t,x) - v_i(t))^2 dx + \sum_{i=1}^n \int_{x_i(t)}^{x_{i-1}(t)} \rho_i(t)v_i(t)(v^{(n)}(t,x) - v_i(t))dx$$

$$+ \frac{1}{2}\sum_{i=1}^n \rho_i(t)v_i^2(t)(x_{i-1}(t) - x_i(t)) + \sum_{i=1}^n Q(\rho_i(t))(x_{i-1}(t) - x_i(t))$$

$$\tag{4.127}$$

Using (3.8), (3.10), (3.31) and the facts that $\sum_{i=1}^{n}(x_{i-1}(t)-x_i(t))=L$, $\rho^*=m/L$, we also get for all $t\geq 0$:

$$\frac{1}{2}\sum_{i=1}^{n}\rho_i(t)v_i^2(t)(x_{i-1}(t)-x_i(t))+\sum_{i=1}^{n}Q(\rho_i(t))(x_{i-1}(t)-x_i(t))$$

$$=\frac{m}{2n}\sum_{i=1}^{n}v_i^2(t)+\sum_{i=1}^{n}Q\left(\frac{m}{n(x_{i-1}(t)-x_i(t))}\right)(x_{i-1}(t)-x_i(t))$$

$$=\frac{m}{2n}\sum_{i=1}^{n}v_i^2(t)+\frac{m}{n}\sum_{i=1}^{n}\Phi(n(x_{i-1}(t)-x_i(t)))$$

$$-\sum_{i=1}^{n}P(\rho^*)\left(\frac{m}{n\rho^*(x_{i-1}(t)-x_i(t))}-1\right)(x_{i-1}(t)-x_i(t)) \quad (4.128)$$

$$=\frac{m}{2n}\sum_{i=1}^{n}v_i^2(t)+\frac{m}{n}\sum_{i=1}^{n}\Phi(n(x_{i-1}(t)-x_i(t)))$$

$$-P(\rho^*)\sum_{i=1}^{n}\frac{m}{n\rho^*}+P(\rho^*)\sum_{i=1}^{n}(x_{i-1}(t)-x_i(t))$$

$$=\frac{m}{2n}\sum_{i=1}^{n}v_i^2(t)+\frac{m}{n}\sum_{i=1}^{n}\Phi(n(x_{i-1}(t)-x_i(t)))=E_n(t)$$

Combining (4.127) and (4.128) we get for all $t\geq 0$:

$$\left|E(\rho^{(n)}[t],v^{(n)}[t])-E_n(t)\right|$$

$$\leq \frac{1}{2}\sum_{i=1}^{n}\int_{x_i(t)}^{x_{i-1}(t)}\left|\rho^{(n)}(t,x)-\rho_i(t)\right|\left(v^{(n)}(t,x)\right)^2 dx+\sum_{i=1}^{n}\int_{x_i(t)}^{x_{i-1}(t)}\left|Q(\rho^{(n)}(t,x))-Q(\rho_i(t))\right|dx \quad (4.129)$$

$$+\frac{1}{2}\sum_{i=1}^{n}\int_{x_i(t)}^{x_{i-1}(t)}\rho_i(t)\left(v^{(n)}(t,x)-v_i(t)\right)^2 dx+\sum_{i=1}^{n}\int_{x_i(t)}^{x_{i-1}(t)}\rho_i(t)|v_i(t)|\left|v^{(n)}(t,x)-v_i(t)\right|dx$$

Using (3.23), (3.26), (3.33), (4.68) and (4.69) we obtain from (4.129):

$$\left|E(\rho^{(n)}[t],v^{(n)}[t])-E_n(t)\right|$$

$$\leq \frac{1}{2}\left(\frac{R_4^2}{2}+\max\left\{|Q'(s)|:\frac{m}{b}\leq s\leq \frac{m}{a}\right\}\right)\sum_{i=1}^{n}|\rho_{i-1}(t)-\rho_i(t)|(x_{i-1}(t)-x_i(t)) \quad (4.130)$$

$$+\frac{m}{6a}\sum_{i=1}^{n}|v_{i-1}(t)-v_i(t)|^2(x_{i-1}(t)-x_i(t))+\frac{mR_4}{2a}\sum_{i=1}^{n}|v_{i-1}(t)-v_i(t)|(x_{i-1}(t)-x_i(t))$$

Using (3.12), (3.24), (3.23), (4.70) and the Cauchy-Schwarz inequality we obtain from (4.130):

$$\left| E(\rho^{(n)}[t], v^{(n)}[t]) - E_n(t) \right|$$

$$\leq \frac{1}{2}\left( \frac{R_4^2}{2} + \max\left\{ |Q'(s)| : \frac{m}{b} \leq s \leq \frac{m}{a} \right\} \right) \left( \sum_{i=1}^{n} |\rho_{i-1}(t) - \rho_i(t)|^2 \right)^{1/2} \left( \sum_{i=1}^{n} (x_{i-1}(t) - x_i(t))^2 \right)^{1/2}$$

$$+ \frac{m}{6a} \sum_{i=1}^{n} \frac{|v_{i-1}(t) - v_i(t)|^2}{x_{i-1}(t) - x_i(t)} (x_{i-1}(t) - x_i(t))^2$$

$$+ \frac{mR_4}{2a} \left( \sum_{i=1}^{n} \frac{|v_{i-1}(t) - v_i(t)|^2}{x_{i-1}(t) - x_i(t)} \right)^{1/2} \left( \sum_{i=1}^{n} (x_{i-1}(t) - x_i(t))^3 \right)^{1/2}$$

$$\leq \frac{\sqrt{bLR_7}}{2n} \left( \frac{R_4^2}{2} + \max\left\{ |Q'(s)| : \frac{m}{b} \leq s \leq \frac{m}{a} \right\} \right) + \frac{mb^2 R_2}{3an^2} + \frac{mR_4 b\sqrt{2LR_2}}{2an}$$

The above estimate implies (4.126) with

$$R_{14} = \frac{\sqrt{bLR_7}}{2} \left( \frac{R_4^2}{2} + \max\left\{ |Q'(s)| : \frac{m}{b} \leq s \leq \frac{m}{a} \right\} \right) + \frac{mb^2 R_2}{3a} + \frac{mR_4 b\sqrt{2LR_2}}{2a}.$$

The proof is complete. ◁

**Lemma 4.7:** *Suppose that Assumption (A) holds. Let a constant $T > 0$ be given. Let $(\rho_0, v_0) \in X$ that satisfies (3.30) be given. Consider the unique solution of (3.1), (3.2), (3.3) with initial condition $(x_1(0), ..., x_{n-1}(0), v_1(0), ..., v_{n-1}(0)) = P_n(\rho_0, v_0) \in \Omega_n$ and define $v^{(n)}, \rho^{(n)} : \mathbb{R}_+ \times [0, L] \to \mathbb{R}$ by means of (3.31)-(3.34). Then there exists a constant $R_{15} > 0$ (independent of $n \geq 2$ but dependent on $\rho_0, v_0, T$) such that the following estimate holds for all $n \geq 2$ and $T \geq t \geq t_0 \geq 0$:*

$$W(\rho^{(n)}[t], v^{(n)}[t]) \leq W(\rho^{(n)}[t_0], v^{(n)}[t_0]) + \frac{R_{15}}{n} \tag{4.131}$$

*where the functional $W(\rho, v)$ is defined by (2.10).*

**Proof:** Since (4.2) holds (which implies that $W_n(t) \leq W_n(t_0)$ for all $n \geq 2$ and $t \geq t_0 \geq 0$), it suffices to show that there exists a constant $R_{16} > 0$ (independent of $n \geq 2$ but dependent on $\delta, \rho_0, v_0, T$) such that the following estimate holds for all $n \geq 2$ and $t \in [0, T]$:

$$\left| W(\rho^{(n)}[t], v^{(n)}[t]) - W_n(t) \right| \leq \frac{R_{16}}{n} \tag{4.132}$$

Using definitions (2.7), (2.10), (3.11), (3.14), (3.31) we get for all $t \geq 0$:

$$W(\rho^{(n)}[t], v^{(n)}[t]) - W_n(t) = \sum_{i=1}^{n} \int_{x_i(t)}^{x_{i-1}(t)} Q(\rho^{(n)}(t, x)) dx$$

$$+ \frac{1}{2} \sum_{i=1}^{n} \int_{x_i(t)}^{x_{i-1}(t)} \rho^{(n)}(t, x) \left( v^{(n)}(t, x) + \frac{\mu(\rho^{(n)}(t, x))}{(\rho^{(n)}(t, x))^2} \rho_x^{(n)}(t, x) \right)^2 dx \tag{4.133}$$

$$- \frac{m}{2n} \sum_{i=1}^{n-1} \left( v_i(t) - nK\left( \frac{m}{\rho_i(t)} \right) + nK\left( \frac{m}{\rho_{i+1}(t)} \right) \right)^2 - \frac{m}{n} \sum_{i=1}^{n} \Phi\left( \frac{m}{\rho_i(t)} \right)$$

Using (3.8), (3.33) and (4.133) we obtain for all $t \geq 0$:

$$W(\rho^{(n)}[t], v^{(n)}[t]) - W_n(t) = I_1 + I_2 + I_3 + I_4 \qquad (4.134)$$

where

$$I_1 = \sum_{i=1}^{n} Q(\rho_i(t))(x_{i-1}(t) - x_i(t)) - \frac{m}{n}\sum_{i=1}^{n} \frac{1}{\rho_i(t)} Q(\rho_i(t))$$
$$- \frac{m}{n} P(\rho^*) \sum_{i=1}^{n} \left( \frac{1}{\rho^*} - \frac{1}{\rho_i(t)} \right) + \sum_{i=1}^{n} \int_{x_i(t)}^{x_{i-1}(t)} \left( Q(\rho^{(n)}(t,x)) - Q(\rho_i(t)) \right) dx \qquad (4.135)$$

$$I_2 = \frac{1}{2} \sum_{i=1}^{n} \int_{x_i(t)}^{x_{i-1}(t)} (\rho^{(n)}(t,x) - \rho_i(t)) \left( v^{(n)}(t,x) + \frac{\mu(\rho^{(n)}(t,x))}{(\rho^{(n)}(t,x))^2} \rho_x^{(n)}(t,x) \right)^2 dx$$
$$+ \frac{1}{2} \sum_{i=1}^{n} \int_{x_i(t)}^{x_{i-1}(t)} \rho_i(t) \left( v^{(n)}(t,x) - v_{i-1}(t) \right)^2 dx \qquad (4.136)$$
$$+ \sum_{i=1}^{n} \int_{x_i(t)}^{x_{i-1}(t)} \rho_i(t)(v^{(n)}(t,x) - v_{i-1}(t)) \left( v_{i-1}(t) + \frac{\mu(\rho^{(n)}(t,x))}{(\rho^{(n)}(t,x))^2} \rho_x^{(n)}(t,x) \right) dx$$

$$I_3 = \frac{1}{2} \sum_{i=1}^{n} \rho_i(t) \int_{x_i(t)}^{x_{i-1}(t)} \left( \rho_x^{(n)}(t,x) \right)^2 \left( \frac{\mu(\rho^{(n)}(t,x))}{(\rho^{(n)}(t,x))^2} - \frac{\mu(\rho_i(t))}{\rho_i^2(t)} \right)^2 dx$$
$$+ \sum_{i=1}^{n} \rho_i(t) \int_{x_i(t)}^{x_{i-1}(t)} \rho_x^{(n)}(t,x) \left( v_{i-1}(t) + \frac{\mu(\rho_i(t))}{\rho_i^2(t)} \rho_x^{(n)}(t,x) \right) \left( \frac{\mu(\rho^{(n)}(t,x))}{(\rho^{(n)}(t,x))^2} - \frac{\mu(\rho_i(t))}{\rho_i^2(t)} \right) dx \qquad (4.137)$$

$$I_4 = \frac{1}{2} \sum_{i=1}^{n} \rho_i(t) \left( v_{i-1}(t) + \frac{\mu(\rho_i(t))}{\rho_i^2(t)} \frac{\rho_{i-1}(t) - \rho_i(t)}{x_{i-1}(t) - x_i(t)} \right)^2 (x_{i-1}(t) - x_i(t))$$
$$- \frac{m}{2n} \sum_{i=2}^{n} \left( v_{i-1}(t) + nK\left(\frac{m}{\rho_i(t)}\right) - nK\left(\frac{m}{\rho_{i-1}(t)}\right) \right)^2 \qquad (4.138)$$

Using (3.23), (3.31), (3.33), (4.68), (4.70), the facts that $\sum_{i=1}^{n}(x_{i-1}(t) - x_i(t)) = L$, $\rho^* = m/L$ and the Cauchy-Schwarz inequality, we get from (4.135):

$$|I_1| \leq \sum_{i=1}^{n} \int_{x_i(t)}^{x_{i-1}(t)} \left| Q(\rho^{(n)}(t,x)) - Q(\rho_i(t)) \right| dx$$
$$\leq \frac{1}{2} \max \left\{ |Q'(s)| : \frac{m}{b} \leq s \leq \frac{m}{a} \right\} \sum_{i=1}^{n} |\rho_{i-1}(t) - \rho_i(t)|(x_{i-1}(t) - x_i(t))$$
$$\leq \frac{1}{2} \max \left\{ |Q'(s)| : \frac{m}{b} \leq s \leq \frac{m}{a} \right\} \left( \sum_{i=1}^{n} |\rho_{i-1}(t) - \rho_i(t)|^2 \right)^{1/2} \left( \sum_{i=1}^{n} (x_{i-1}(t) - x_i(t))^2 \right)^{1/2} \qquad (4.139)$$
$$\leq \frac{\sqrt{bLR_7}}{2n} \max \left\{ |Q'(s)| : \frac{m}{b} \leq s \leq \frac{m}{a} \right\}$$

Using (3.12), (3.23), (3.24), (3.26), (3.31), (3.33), (4.68), (4.69), (4.70), (4.77) the facts that $\sum_{i=1}^{n}(x_{i-1}(t) - x_i(t)) = L$, $t \in [0,T]$ and the Cauchy-Schwarz inequality, we get from (4.136):

$$|I_2| \leq \frac{1}{2}\left(R_4 + \max\left\{\frac{\mu(s)}{s^2} : \frac{m}{b} \leq s \leq \frac{m}{a}\right\} R_9 \exp(\kappa T)\right)^2 \sum_{i=1}^{n}|\rho_{i-1}(t) - \rho_i(t)|(x_{i-1}(t) - x_i(t))$$

$$+ \frac{m}{2a}\sum_{i=1}^{n}|v_{i-1}(t) - v_i(t)|^2 (x_{i-1}(t) - x_i(t))$$

$$+ \frac{m}{a}\left(R_4 + \max\left\{\frac{\mu(s)}{s^2} : \frac{m}{b} \leq s \leq \frac{m}{a}\right\} R_9 \exp(\kappa T)\right) \sum_{i=1}^{n}|v_{i-1}(t) - v_i(t)|(x_{i-1}(t) - x_i(t)) \quad (4.140)$$

$$\leq \frac{\sqrt{bLR_7}}{2n}\left(R_4 + \max\left\{\frac{\mu(s)}{s^2} : \frac{m}{b} \leq s \leq \frac{m}{a}\right\} R_9 \exp(\kappa T)\right)^2 + \frac{mb^2 R_2}{an^2}$$

$$+ \frac{mb\sqrt{2LR_2}}{an}\left(R_4 + \max\left\{\frac{\mu(s)}{s^2} : \frac{m}{b} \leq s \leq \frac{m}{a}\right\} R_9 \exp(\kappa T)\right)$$

Using (3.23), (3.26), (3.31), (3.33), (4.68), (4.77) and the fact that $t \in [0,T]$, we get from (4.137):

$$|I_3| \leq \frac{m}{2a} R_9^2 \exp(2\kappa T) \bar{K}_1^2 \sum_{i=1}^{n}|\rho_{i-1}(t) - \rho_i(t)|^2 (x_{i-1}(t) - x_i(t))$$

$$+ \frac{m}{a} \bar{K}_1 R_9 \exp(\kappa T)(R_4 + \bar{K}_2 R_9 \exp(\kappa T)) \sum_{i=1}^{n}|\rho_{i-1}(t) - \rho_i(t)|(x_{i-1}(t) - x_i(t)) \quad (4.141)$$

where $\bar{K}_1 = \max\left\{\left|\frac{\mu'(s)}{s^2} - \frac{2\mu(s)}{s^3}\right| : \frac{m}{b} \leq s \leq \frac{m}{a}\right\}$, $\bar{K}_2 = \max\left\{\frac{\mu(s)}{s^2} : \frac{m}{b} \leq s \leq \frac{m}{a}\right\}$. Exploiting (3.23), (4.141) the Cauchy-Schwarz inequality and the fact that $\sum_{i=1}^{n}(x_{i-1}(t) - x_i(t)) = L$, we get:

$$|I_3| \leq \frac{m}{an}\left(\frac{bR_7}{2} R_9 \bar{K}_1 + \sqrt{bLR_7}(R_4 + \bar{K}_2 R_9)\right) \bar{K}_1 R_9 \exp(2\kappa T) \quad (4.142)$$

Finally, we deal with the term $I_4$. Equation (4.138) in conjunction with (3.31) gives:

$$I_4 = \frac{m}{2n}\sum_{i=2}^{n}\Delta_i\left(\frac{\mu(\rho_i(t))}{\rho_i^2(t)}\frac{\rho_{i-1}(t) - \rho_i(t)}{x_{i-1}(t) - x_i(t)} - nK\left(\frac{m}{\rho_i(t)}\right) + nK\left(\frac{m}{\rho_{i-1}(t)}\right)\right) \quad (4.143)$$

where $\Delta_i = 2v_{i-1}(t) + \frac{\mu(\rho_i(t))}{\rho_i^2(t)}\frac{\rho_{i-1}(t) - \rho_i(t)}{x_{i-1}(t) - x_i(t)} + nK\left(\frac{m}{\rho_i(t)}\right) - nK\left(\frac{m}{\rho_{i-1}(t)}\right)$. Clearly, we have by using (3.5):

$$K\left(\frac{m}{\rho_i(t)}\right) - K\left(\frac{m}{\rho_{i-1}(t)}\right) = \int_{\frac{m}{\rho_{i-1}(t)}}^{\frac{m}{\rho_i(t)}} K'(x)\,dx = \int_{\frac{m}{\rho_{i-1}(t)}}^{\frac{m}{\rho_i(t)}} \frac{1}{mx}\mu\left(\frac{m}{x}\right)dx = \frac{1}{m}\int_{\rho_i(t)}^{\rho_{i-1}(t)} \frac{\mu(u)}{u}\,du$$

The above equation implies the existence of $\lambda \in [0,1]$ such that

$$K\left(\frac{m}{\rho_i(t)}\right) - K\left(\frac{m}{\rho_{i-1}(t)}\right) = \frac{1}{m} \frac{\mu(\rho_i(t) + \lambda(\rho_{i-1}(t) - \rho_i(t)))}{\rho_i(t) + \lambda(\rho_{i-1}(t) - \rho_i(t))} (\rho_{i-1}(t) - \rho_i(t))$$

Consequently, we get by using (3.31):

$$\left| \frac{\mu(\rho_i(t))}{\rho_i^2(t)} \frac{\rho_{i-1}(t) - \rho_i(t)}{x_{i-1}(t) - x_i(t)} - nK\left(\frac{m}{\rho_i(t)}\right) + nK\left(\frac{m}{\rho_{i-1}(t)}\right) \right|$$

$$= \left| \frac{\mu(\rho_i(t))}{\rho_i^2(t)} \frac{\rho_{i-1}(t) - \rho_i(t)}{x_{i-1}(t) - x_i(t)} - \frac{n}{m} \frac{\mu(\rho_i(t) + \lambda(\rho_{i-1}(t) - \rho_i(t)))}{\rho_i(t) + \lambda(\rho_{i-1}(t) - \rho_i(t))} (\rho_{i-1}(t) - \rho_i(t)) \right|$$

$$= \frac{1}{\rho_i(t)} \left| \frac{\mu(\rho_i(t))}{\rho_i(t)} - \frac{\mu(\rho_i(t) + \lambda(\rho_{i-1}(t) - \rho_i(t)))}{\rho_i(t) + \lambda(\rho_{i-1}(t) - \rho_i(t))} \right| \left| \frac{\rho_{i-1}(t) - \rho_i(t)}{x_{i-1}(t) - x_i(t)} \right|$$

Using (3.23), (3.31), (4.71) and the facts that $\lambda \in [0,1]$, $t \in [0,T]$ we get:

$$\left| \frac{\mu(\rho_i(t))}{\rho_i^2(t)} \frac{\rho_{i-1}(t) - \rho_i(t)}{x_{i-1}(t) - x_i(t)} - nK\left(\frac{m}{\rho_i(t)}\right) + nK\left(\frac{m}{\rho_{i-1}(t)}\right) \right|$$
$$\leq \frac{b^2 R_8^2}{nm} \bar{K}_3 \exp(\kappa T) \tag{4.144}$$

where $\bar{K}_3 = \max\left\{ \left| \frac{\mu'(s)}{s} - \frac{\mu(s)}{s^2} \right| : \frac{m}{b} \leq s \leq \frac{m}{a} \right\}$. Similar arguments and the use of (3.26) give:

$$|\Delta_i| \leq 2R_4 + \bar{K}_2 R_8 b \exp(\kappa T) + \frac{R_8 b^2}{m} \exp(\kappa T) \bar{K}_4 \tag{4.145}$$

where $\bar{K}_4 = \max\left\{ \frac{\mu(s)}{s} : \frac{m}{b} \leq s \leq \frac{m}{a} \right\}$. Combining (4.143), (4.144) and (4.145) we get:

$$|I_4| \leq \frac{b^2 R_8^2}{2n} \left( 2R_4 + \bar{K}_2 R_8 b \exp(\kappa T) + \frac{R_8 b^2}{m} \exp(\kappa T) \bar{K}_4 \right) \bar{K}_3 \exp(\kappa T) \tag{4.146}$$

Estimate (4.131) is a direct consequence of (4.134), (4.139), (4.140), (4.142) and (4.146). The proof is complete. ◁

**Lemma 4.8:** *Suppose that Assumption (A) holds. Let a constant $T > 0$ be given. Let $(\rho_0, v_0) \in X$ that satisfies (3.30) be given. Consider the unique solution of (3.1), (3.2), (3.3) with initial condition $(x_1(0), \ldots, x_{n-1}(0), v_1(0), \ldots, v_{n-1}(0)) = P_n(\rho_0, v_0) \in \Omega_n$ and define $v^{(n)}, \rho^{(n)} : \mathbb{R}_+ \times [0, L] \to \mathbb{R}$ by means of (3.31)-(3.34). Then there exist constants $R_{17}, R_{18} > 0$ (independent of $n \geq 2$ but dependent on $\rho_0, v_0, T$) such that the following estimates hold for all $n \geq 2$ and $t, t_0 \in [0, T]$ with $t \geq t_0$:*

$$\left\| \dot{\rho}^{(n)}[t] \right\|_2 \leq R_{17} \tag{4.147}$$

$$\left| \int_{t_0}^{t} \int_{0}^{L} \varphi(s, x) \dot{v}^{(n)}(s, x) \, dx \, ds \right| \leq R_{18} \left( \int_{t_0}^{t} \left\| \varphi_x[s] \right\|_2^2 \, ds \right)^{1/2}, \text{ for all } \varphi \in L^2\left((0, T); H_0^1((0, L))\right) \tag{4.148}$$

**Remark:** Estimate (4.148) shows that $\dot{v}^{(n)} \in L^2\big((0,T); H^{-1}((0,L))\big)$.

**Proof:** Using definition (3.35) we get:

$$\begin{aligned}
\left\|\dot{\rho}^{(n)}[t]\right\|_2^2 &= \sum_{i=1}^n \int_{x_i(t)}^{x_{i-1}(t)} \left(\dot{\rho}^{(n)}(t,x)\right)^2 dx \\
&\leq 4\sum_{i=1}^n \int_{x_i(t)}^{x_{i-1}(t)} \left(\dot{\rho}_i(t)\right)^2 dx + 4\sum_{i=1}^n \int_{x_i(t)}^{x_{i-1}(t)} \left(\frac{\dot{\rho}_{i-1}(t)-\dot{\rho}_i(t)}{x_{i-1}(t)-x_i(t)}(x-x_i(t))\right)^2 dx \\
&+ 4\sum_{i=1}^n \int_{x_i(t)}^{x_{i-1}(t)} \left(\frac{\rho_{i-1}(t)-\rho_i(t)}{(x_{i-1}(t)-x_i(t))^2}(v_{i-1}(t)-v_i(t))(x-x_i(t))\right)^2 dx \\
&+ 4\sum_{i=1}^n \int_{x_i(t)}^{x_{i-1}(t)} \left(\frac{\rho_{i-1}(t)-\rho_i(t)}{x_{i-1}(t)-x_i(t)}v_i(t)\right)^2 dx \\
&= 4\sum_{i=1}^n \dot{\rho}_i^2(t)(x_{i-1}(t)-x_i(t)) + \frac{4}{3}\sum_{i=1}^n (\dot{\rho}_{i-1}(t)-\dot{\rho}_i(t))^2 (x_{i-1}(t)-x_i(t)) \\
&+ \frac{4}{3}\sum_{i=1}^n \frac{(\rho_{i-1}(t)-\rho_i(t))^2}{x_{i-1}(t)-x_i(t)}(v_{i-1}(t)-v_i(t))^2 + 4\sum_{i=1}^n \frac{(\rho_{i-1}(t)-\rho_i(t))^2}{x_{i-1}(t)-x_i(t)} v_i^2(t)
\end{aligned} \quad (4.149)$$

Using (3.26) we obtain from (4.149):

$$\left\|\dot{\rho}^{(n)}[t]\right\|_2^2 \leq \frac{8}{3}\sum_{i=1}^n \dot{\rho}_{i-1}^2(t)(x_{i-1}(t)-x_i(t)) + \frac{20}{3}\sum_{i=1}^n \dot{\rho}_i^2(t)(x_{i-1}(t)-x_i(t)) \\
+ \frac{20 R_4^2}{3}\sum_{i=1}^n \frac{(\rho_{i-1}(t)-\rho_i(t))^2}{x_{i-1}(t)-x_i(t)} \quad (4.150)$$

Using (3.31), (3.32), we get from (4.150) and (3.23):

$$\left\|\dot{\rho}^{(n)}[t]\right\|_2^2 \leq \frac{12b}{n}\sum_{i=1}^n \dot{\rho}_i^2(t) + \frac{20 R_4^2}{3a} n \sum_{i=1}^n (\rho_{i-1}(t)-\rho_i(t))^2 \quad (4.151)$$

Using (3.1), (3.2), (3.3) and (3.31) we get from (4.151) and (3.23):

$$\left\|\dot{\rho}^{(n)}[t]\right\|_2^2 \leq \frac{12b}{n}\sum_{i=1}^n \frac{(v_{i-1}(t)-v_i(t))^2}{(x_{i-1}(t)-x_i(t))^2}\rho_i^2(t) + \frac{20 R_4^2}{3a} R_7 \quad (4.152)$$

Using (4.152) in conjunction with (3.12), (3.23), (3.24) we obtain:

$$\left\|\dot{\rho}^{(n)}[t]\right\|_2 \leq R_{17} = \sqrt{\frac{24 m^2 b}{a^3} R_2 + \frac{20 R_4^2}{3a} R_7} \quad (4.153)$$

Equations (3.2), (3.31) imply the following equations for $i=1,...,n-1$:

$$\dot{v}_i(t) = n\Phi'\left(\frac{m}{\rho_i(t)}\right) - n\Phi'\left(\frac{m}{\rho_{i+1}(t)}\right)$$

$$+ n\left(G\left(\frac{m}{\rho_i(t)}\right) - G\left(\frac{m}{\rho_{i+1}(t)}\right)\right)\frac{v_{i-1}(t) - v_i(t)}{x_{i-1}(t) - x_i(t)} \qquad (4.154)$$

$$- nG\left(\frac{m}{\rho_{i+1}(t)}\right)\left(\frac{v_i(t) - v_{i+1}(t)}{x_i(t) - x_{i+1}(t)} - \frac{v_{i-1}(t) - v_i(t)}{x_{i-1}(t) - x_i(t)}\right)$$

where $G(x) = xK'(x)$, for $x > 0$. Equations (4.154) imply the following estimates for $i = 1, \ldots, n-1$:

$$|\dot{v}_i(t)| \leq nK_1 |\rho_i(t) - \rho_{i+1}(t)| + nK_2 |\rho_i(t) - \rho_{i+1}(t)| \frac{|v_{i-1}(t) - v_i(t)|}{x_{i-1}(t) - x_i(t)}$$

$$+ nK_3 \left|\frac{v_i(t) - v_{i+1}(t)}{x_i(t) - x_{i+1}(t)} - \frac{v_{i-1}(t) - v_i(t)}{x_{i-1}(t) - x_i(t)}\right| \qquad (4.155)$$

where $K_1 = \max\left\{\left|\Phi''\left(\frac{m}{s}\right)\right|\frac{m}{s^2} : \frac{m}{a} \leq s \leq \frac{m}{b}\right\}$, $K_2 = \max\left\{\left|G'\left(\frac{m}{s}\right)\right|\frac{m}{s^2} : \frac{m}{a} \leq s \leq \frac{m}{b}\right\}$,

$K_3 = \max\left\{\left|G\left(\frac{m}{s}\right)\right| : \frac{m}{a} \leq s \leq \frac{m}{b}\right\}$.

Let arbitrary $\varphi \in L^2\left((0,T); H_0^1((0,L))\right)$ be given. Definition (3.36) implies the following equations for all $t, t_0 \in [0,T]$ with $t \geq t_0$:

$$\int_{t_0}^{t}\int_0^L \varphi(s,x)\dot{v}^{(n)}(s,x)dxds = \sum_{i=1}^n \int_{t_0}^{t}\int_{x_i(s)}^{x_{i-1}(s)} \varphi(s,x)\left(\dot{v}_i(s) - \frac{v_{i-1}(s) - v_i(s)}{x_{i-1}(s) - x_i(s)}v_i(s)\right)dxds$$

$$+ \sum_{i=1}^n \int_{t_0}^{t}\int_{x_i(s)}^{x_{i-1}(s)} \varphi(s,x)\left(\frac{\dot{v}_{i-1}(s) - \dot{v}_i(s)}{x_{i-1}(s) - x_i(s)} - \frac{(v_{i-1}(s) - v_i(t))^2}{(x_{i-1}(s) - x_i(s))^2}\right)(x - x_i(s))dxds \qquad (4.156)$$

Using (3.1), (3.3), (3.12), (3.23), the fact that $\|\varphi[s]\|_\infty \leq \sqrt{L}\|\varphi_x[s]\|_2$ for $s \in [0,T]$ a.e. (a consequence of the fact that $\varphi[s] \in H_0^1((0,L))$ for $s \in [0,T]$ a.e.) we obtain from (4.156) for all $t, t_0 \in [0,T]$ with $t \geq t_0$:

$$\left| \int_{t_0}^{t} \int_{0}^{L} \varphi(s,x) \dot{v}^{(n)}(s,x) dx ds \right| \leq \sum_{i=1}^{n} \int_{t_0}^{t} \int_{x_i(s)}^{x_{i-1}(s)} |\varphi(s,x)| \left( |\dot{v}_i(s)| + \left| \frac{v_{i-1}(s) - v_i(s)}{x_{i-1}(s) - x_i(s)} \right| |v_i(s)| \right) dx ds$$

$$+ \sum_{i=1}^{n} \int_{t_0}^{t} \int_{x_i(s)}^{x_{i-1}(s)} |\varphi(s,x)| \left( \frac{|\dot{v}_{i-1}(s)| + |\dot{v}_i(s)|}{x_{i-1}(s) - x_i(s)} + \frac{(v_{i-1}(s) - v_i(s))^2}{(x_{i-1}(s) - x_i(s))^2} \right) (x - x_i(s)) dx ds$$

$$\leq \sum_{i=1}^{n} \int_{t_0}^{t} \|\varphi[s]\|_{\infty} |\dot{v}_i(s)| (x_{i-1}(s) - x_i(s)) ds + \sum_{i=1}^{n} \int_{t_0}^{t} \|\varphi[s]\|_{\infty} |v_{i-1}(s) - v_i(s)| |v_i(s)| ds$$

$$+ \frac{1}{2} \sum_{i=1}^{n} \int_{t_0}^{t} \|\varphi[s]\|_{\infty} (|\dot{v}_{i-1}(s)| + |\dot{v}_i(s)|)(x_{i-1}(s) - x_i(s)) ds + \frac{1}{2} \sum_{i=1}^{n} \int_{t_0}^{t} (v_{i-1}(s) - v_i(s))^2 \|\varphi[s]\|_{\infty} dt$$

$$\leq \frac{2b\sqrt{L}}{n} \int_{t_0}^{t} \|\varphi_x[s]\|_2 \sum_{i=1}^{n-1} |\dot{v}_i(s)| ds + R_4 \sqrt{L} \int_{t_0}^{t} \|\varphi_x[s]\|_2 \sum_{i=1}^{n} |v_{i-1}(s) - v_i(s)| ds$$

$$+ \frac{b\sqrt{L}}{n} \int_{t_0}^{t} \|\varphi_x[s]\|_2 Z_n(s) ds$$

(4.157)

Estimate (4.157) combined with estimate (4.155) and the use of the Cauchy-Schwarz inequality gives for all $t, t_0 \in [0, T]$ with $t \geq t_0$:

$$\left| \int_{t_0}^{t} \int_{0}^{L} \varphi(s,x) \dot{v}^{(n)}(s,x) dx ds \right| \leq 2bK_1 \sqrt{L} \int_{t_0}^{t} \|\varphi_x[s]\|_2 \sum_{i=1}^{n-1} |\rho_i(s) - \rho_{i+1}(s)| ds$$

$$+ 2bK_2 \sqrt{L} \int_{t_0}^{t} \|\varphi_x[s]\|_2 \sum_{i=1}^{n-1} |\rho_i(s) - \rho_{i+1}(s)| \frac{|v_{i-1}(s) - v_i(s)|}{x_{i-1}(s) - x_i(s)} ds$$

$$+ 2bK_3 \sqrt{L} \int_{t_0}^{t} \|\varphi_x[s]\|_2 \sum_{i=1}^{n-1} \left| \frac{v_i(s) - v_{i+1}(s)}{x_i(s) - x_{i+1}(s)} - \frac{v_{i-1}(s) - v_i(s)}{x_{i-1}(s) - x_i(s)} \right| ds + \frac{b\sqrt{L}}{n} \int_{t_0}^{t} \|\varphi_x[s]\|_2 Z_n(s) ds$$

$$+ R_4 \sqrt{L} \int_{t_0}^{t} \|\varphi_x[s]\|_2 \left( \sum_{i=1}^{n} \frac{|v_{i-1}(s) - v_i(s)|^2}{x_{i-1}(s) - x_i(s)} \right)^{1/2} \left( \sum_{i=1}^{n} (x_{i-1}(s) - x_i(s)) \right)^{1/2} ds$$

(4.158)

Using (4.71), (3.23), (3.12), the Cauchy-Schwarz inequality, the facts that $t, t_0 \in [0, T]$, $\sum_{i=1}^{n} (x_{i-1}(s) - x_i(s)) = L$, we get from (4.158) for all $t, t_0 \in [0, T]$ with $t \geq t_0$:

$$\left| \int_{t_0}^{t} \int_{0}^{L} \varphi(s,x) \dot{v}^{(n)}(s,x) dx ds \right| \leq 2bK_1 R_8 \exp(\kappa T) \sqrt{L} \int_{t_0}^{t} \|\varphi_x[s]\|_2 ds$$

$$+ 2bK_2 R_8 \exp(\kappa T) \frac{\sqrt{L}}{a} \int_{t_0}^{t} \|\varphi_x[s]\|_2 \sum_{i=1}^{n-1} |v_{i-1}(s) - v_i(s)| ds$$

$$+ 2bK_3 \sqrt{L} \int_{t_0}^{t} \|\varphi_x[s]\|_2 \left( n \sum_{i=1}^{n-1} \left| \frac{v_i(s) - v_{i+1}(s)}{x_i(s) - x_{i+1}(s)} - \frac{v_{i-1}(s) - v_i(s)}{x_{i-1}(s) - x_i(s)} \right|^2 \right)^{1/2} ds \quad (4.159)$$

$$+ \frac{b\sqrt{L}}{n} \int_{t_0}^{t} \|\varphi_x[s]\|_2 Z_n(s) ds + R_4 L \int_{t_0}^{t} \|\varphi_x[s]\|_2 \sqrt{2Z_n(s)} ds$$

Using (3.24), (3.12), the Cauchy-Schwarz inequality, we get from (4.159) for all $t, t_0 \in [0, T]$ with $t \geq t_0$:

$$\left| \int_{t_0}^{t} \int_{0}^{L} \varphi(s,x) \dot{v}^{(n)}(s,x) dx ds \right| \leq \left( 2bK_1 R_8 \exp(\kappa T) \sqrt{L} + \frac{bR_2 \sqrt{L}}{n} + R_4 L \sqrt{2R_2} \right) \int_{t_0}^{t} \|\varphi_x[s]\|_2 ds$$

$$+ 2bK_2 R_8 \exp(\kappa T) \frac{\sqrt{L}}{a} \int_{t_0}^{t} \|\varphi_x[s]\|_2 \left( \sum_{i=1}^{n-1} \frac{|v_{i-1}(s) - v_i(s)|^2}{x_{i-1}(s) - x_i(s)} \right)^{1/2} \left( \sum_{i=1}^{n-1} (x_{i-1}(s) - x_i(s)) \right)^{1/2} ds$$

$$+ 2bK_3 \sqrt{L} \left( n \int_{t_0}^{t} \sum_{i=1}^{n-1} \left| \frac{v_i(s) - v_{i+1}(s)}{x_i(s) - x_{i+1}(s)} - \frac{v_{i-1}(s) - v_i(s)}{x_{i-1}(s) - x_i(s)} \right|^2 ds \right)^{1/2} \left( \int_{t_0}^{t} \|\varphi_x[s]\|_2^2 ds \right)^{1/2}$$

$$(4.160)$$

Finally, using (3.12), (3.24), (3.27), the fact that $\sum_{i=1}^{n} (x_{i-1}(s) - x_i(s)) = L$ the Cauchy-Schwarz inequality, we get from (4.159) for all $t, t_0 \in [0, T]$ with $t \geq t_0$:

$$\left| \int_{t_0}^{t} \int_{0}^{L} \varphi(s,x) \dot{v}^{(n)}(s,x) dx ds \right| \leq 2bK_1 R_8 \exp(\kappa T) \sqrt{L} \sqrt{t - t_0} \left( \int_{t_0}^{t} \|\varphi_x[s]\|_2^2 ds \right)^{1/2}$$

$$+ \left( \frac{bR_2 \sqrt{L}}{n} + R_4 L \sqrt{2R_2} + 2bK_2 R_8 \exp(\kappa T) \sqrt{2R_2} \frac{L}{a} \right) \sqrt{t - t_0} \left( \int_{t_0}^{t} \|\varphi_x[s]\|_2^2 ds \right)^{1/2}$$

$$+ 2bK_3 \sqrt{L(R_5 + R_6 t)} \left( \int_{t_0}^{t} \|\varphi_x[s]\|_2^2 ds \right)^{1/2}$$

$$(4.161)$$

Estimate (4.161) implies estimate (4.148) with

$$R_{18} = 2bK_1 R_8 \exp(\kappa T) \sqrt{TL} + 2bK_3 \sqrt{L(R_5 + R_6 T)}$$

$$+ bR_2 \sqrt{TL} + R_4 L \sqrt{2TR_2} + 2bK_2 R_8 \exp(\kappa T) \sqrt{2TR_2} \frac{L}{a}$$

The proof is complete. ◁

**Lemma 4.9:** *Suppose that Assumption (A) holds. Let a constant $T>0$ be given. Let $(\rho_0, v_0) \in X$ that satisfies (3.30) be given. Consider the unique solution of (3.1), (3.2), (3.3) with initial condition $\left(x_1(0),...,x_{n-1}(0), v_1(0),...,v_{n-1}(0)\right) = P_n(\rho_0, v_0) \in \Omega_n$ and define $v^{(n)}, \rho^{(n)} : \mathbb{R}_+ \times [0,L] \to \mathbb{R}$ by means of (3.31)-(3.34). Then there exist constants $R_{19}, R_{20} > 0$ (independent of $n \geq 2$ but dependent on $\rho_0, v_0, T$) such that the following estimates hold for all $n \geq 2$ and $t, t_0 \in [0,T]$:*

$$\left\| \rho^{(n)}[t] - \rho^{(n)}[t_0] \right\|_2 \leq R_{19} |t - t_0|^{1/2} \quad (4.162)$$

$$\left\| v^{(n)}[t] - v^{(n)}[t_0] \right\|_2 \leq R_{20} |t - t_0|^{1/4} \quad (4.163)$$

**Proof:** Theorem 3 on page 287 in [5] in conjunction with (4.147), (4.148) and the triangle inequality implies that the following estimates hold for all $t, t_0 \in [0,T]$ with $t \geq t_0$:

$$\left\| \rho^{(n)}[t] - \rho^{(n)}[t_0] \right\|_2^2 \leq 2 \int_{t_0}^{t} \left\| \dot{\rho}^{(n)}[s] \right\|_2 \left\| \rho^{(n)}[s] - \rho^{(n)}[t_0] \right\|_2 ds$$
$$\leq 2 R_{17} \int_{t_0}^{t} \left( \left\| \rho^{(n)}[s] \right\|_2 + \left\| \rho^{(n)}[t_0] \right\|_2 \right) ds \quad (4.164)$$

$$\left\| v^{(n)}[t] - v^{(n)}[t_0] \right\|_2^2 \leq 2 \int_{t_0}^{t} \int_0^L \dot{v}^{(n)}(s,x) \left( v^{(n)}(s,x) - v^{(n)}(t_0,x) \right) dx\, ds$$
$$\leq 2 R_{18} \left( \int_{t_0}^{t} \left\| v_x^{(n)}[s] - v_x^{(n)}[t_0] \right\|_2^2 ds \right)^{1/2} \leq 2 R_{18} \left( 2 \int_{t_0}^{t} \left\| v_x^{(n)}[s] \right\|_2^2 ds + 2 \int_{t_0}^{t} \left\| v_x^{(n)}[t_0] \right\|_2^2 ds \right)^{1/2} \quad (4.165)$$

Using (4.164), (4.165) in conjunction with (4.68) and (4.75) we get for all $t, t_0 \in [0,T]$ with $t \geq t_0$:

$$\left\| \rho^{(n)}[t] - \rho^{(n)}[t_0] \right\|_2^2 \leq 4\sqrt{L} R_{17} \frac{m}{a} (t - t_0) \quad (4.166)$$

$$\left\| v^{(n)}[t] - v^{(n)}[t_0] \right\|_2^2 \leq 4 R_{18} \sqrt{R_9} |t - t_0|^{1/2} \quad (4.167)$$

Estimates (4.162), (4.163) are direct consequences of estimates (4.166), (4.167). The proof is complete. ◁

For the proof of Theorem 3.3 we use the following extension of the Arzela-Ascoli theorem. Its proof is exactly the same with the usual Arzela-Ascoli theorem for real-valued functions and exploits the properties of the spaces $H^1((0,L))$ and $L^2((0,L))$.

<u>Fact 6:</u> *Let $\{ f_n : [0,T] \to H^1((0,L)), n \geq 1 \}$ be a given sequence of functions for which the following two properties hold:*

*i) Uniform Boundedness: there exists a constant $M > 0$ such that the estimate $\left\| f_n[t] \right\|_2 + \left\| (f_n[t])_x \right\|_2 \leq M$ holds for all $t \in [0,T]$ and $n \geq 1$.*

*ii) Uniform Equicontinuity: for every $\varepsilon > 0$ there exists $\delta > 0$ such that the estimate $\|f_n[t] - f_n[t_0]\|_2 < \varepsilon$ holds for all $n \geq 1$, $t, t_0 \in [0,T]$ with $|t - t_0| < \delta$.*

*Then there exists a subsequence of $\{f_n : n \geq 1\}$ that converges in the topology of $C^0([0,T]; L^2((0,L)))$ to a function $f \in C^0([0,T]; L^2((0,L)))$.*

**Proof of Theorem 3.3:** Let $(\rho_0, v_0) \in X$ that satisfies (3.30) be given. Consider the unique solution of (3.1), (3.2), (3.3) with initial condition $(x_1(0), ..., x_{n-1}(0), v_1(0), ..., v_{n-1}(0)) = P_n(\rho_0, v_0) \in \Omega_n$ and define $v^{(n)}, \rho^{(n)} : \mathbb{R}_+ \times [0, L] \to \mathbb{R}$ by means of (3.31)-(3.34). It follows from (4.68), (4.69), (4.75) and (4.76) that for every $T > 0$ the sequence $\{(\rho^{(n)}, v^{(n)}) : n \geq 2\} \subseteq L^\infty((0,T); H^1((0,L))) \times L^\infty((0,T); H^1_0((0,L)))$ is bounded. Moreover, by virtue of (4.77), (4.147), (4.148) the sequences $\{\rho^{(n)}_x : n \geq 2\} \subseteq L^\infty((0,T) \times (0,L))$, $\{\dot{\rho}^{(n)} : n \geq 2\} \subseteq L^\infty((0,T); L^2((0,L)))$, $\{\dot{v}^{(n)} : n \geq 2\} \subseteq L^2((0,T); H^{-1}((0,L)))$ are bounded. Consequently, there exist $\rho \in L^\infty((0,T); H^1((0,L)))$, $v \in L^\infty((0,T); H^1_0((0,L)))$ and a subsequence of $\{(\rho^{(n)}, v^{(n)}) : n \geq 2\} \subseteq L^\infty((0,T); H^1((0,L))) \times L^\infty((0,T); H^1_0((0,L)))$ (still denoted by $\{(\rho^{(n)}, v^{(n)})\}$) for which the following convergence properties hold:

$$v^{(n)} \to v \text{ in } L^\infty((0,T); H^1_0((0,L))) \text{ weak star}$$

$$\dot{v}^{(n)} \to v_t \text{ in } L^2((0,T); H^{-1}((0,L))) \text{ weakly}$$

$$\rho^{(n)} \to \rho \text{ in } L^\infty((0,T); H^1((0,L))) \text{ weak star}$$

$$\dot{\rho}^{(n)} \to \rho_t \text{ in } L^\infty((0,T); L^2((0,L))) \text{ weak star}$$

$$\rho^{(n)}_x \to \rho_x \text{ in } L^\infty((0,T) \times (0,L)) \text{ weak star}$$

Lemma 4.9, Fact 6 and estimates (4.68), (4.75), (4.76) allow us to conclude that there exists a subsequence of $\{(\rho^{(n)}, v^{(n)}) : n \geq 2\} \subseteq L^\infty((0,T); H^1((0,L))) \times L^\infty((0,T); H^1_0((0,L)))$ for which

$$\rho^{(n)} \to \rho \text{ in } C^0([0,T]; L^2((0,L))) \text{ strongly}$$

$$v^{(n)} \to v \text{ in } C^0([0,T]; L^2((0,L))) \text{ strongly}$$

Moreover, by virtue of (4.162), (4.163) there exist constants $R_{19}, R_{20} > 0$ such that the following estimates hold for all $t, t_0 \in [0, T]$:

$$\|\rho[t] - \rho[t_0]\|_2 \leq R_{19} |t - t_0|^{1/2} \tag{4.168}$$

$$\|v[t]-v[t_0]\|_2 \leq R_{20}|t-t_0|^{1/4} \qquad (4.169)$$

Estimates (4.168), (4.169) imply that $\rho \in C^{0,1/2}\left([0,T];L^2\left((0,L)\right)\right)$ and $v \in C^{0,1/4}\left([0,T];L^2\left((0,L)\right)\right)$.

Define next the operator $\tilde{K}:L^2\left((0,L)\right) \to L^\infty\left((0,L)\right)$ by means of the following formula:

$$(\tilde{K}f)(x) = \begin{cases} m/a & \text{if } f(x) > m/a \\ f(x) & \text{if } m/b \leq f(x) \leq m/a, \\ m/b & \text{if } f(x) < m/b \end{cases} \text{ for all } f \in L^2\left((0,L)\right) \qquad (4.170)$$

where $0 < a \leq b$ are the constants involved in (4.68). Since $\left|(\tilde{K}f)(x)-(\tilde{K}g)(x)\right| \leq |f(x)-g(x)|$, it follows that $\|\tilde{K}f - \tilde{K}g\|_2 \leq \|f-g\|_2$ for all $f,g \in L^2\left((0,L)\right)$. Moreover, by virtue of (4.68) it holds that $\tilde{K}\rho^{(n)}[t] = \rho^{(n)}[t]$ and the triangle inequality gives:

$$\sup_{t\in[0,T]}\left(\|\rho[t]-\tilde{K}\rho[t]\|_2\right) \leq \sup_{t\in[0,T]}\left(\|\rho^{(n)}[t]-\tilde{K}\rho[t]\|_2 + \|\rho^{(n)}[t]-\rho[t]\|_2\right)$$
$$\leq \sup_{t\in[0,T]}\left(\|\rho^{(n)}[t]-\tilde{K}\rho[t]\|_2\right) + \sup_{t\in[0,T]}\left(\|\rho^{(n)}[t]-\rho[t]\|_2\right)$$
$$= \sup_{t\in[0,T]}\left(\|\tilde{K}\rho^{(n)}[t]-\tilde{K}\rho[t]\|_2\right) + \sup_{t\in[0,T]}\left(\|\rho^{(n)}[t]-\rho[t]\|_2\right)$$
$$\leq 2\sup_{t\in[0,T]}\left(\|\rho^{(n)}[t]-\rho[t]\|_2\right)$$

Since $\rho^{(n)} \to \rho$ in $C^0\left([0,T];L^2\left((0,L)\right)\right)$ strongly, it follows that $\tilde{K}\rho[t]=\rho[t]$ for all $t \in [0,T]$. It follows that for all $t \in [0,T]$ it holds that

$$\frac{m}{b} \leq \rho(t,x) \leq \frac{m}{a}, \text{ for } x \in (0,L) \text{ a.e.} \qquad (4.171)$$

Estimate (4.171) shows that for every $t \in [0,T]$ the function $\rho[t]$ is of class $L^\infty\left((0,L)\right)$. A similar argument can be provided for $v$ by using (4.69). It then follows that for all $t \in [0,T]$ it holds that

$$|v(t,x)| \leq R_4, \text{ for } x \in (0,L) \text{ a.e.} \qquad (4.172)$$

which also shows that for every $t \in [0,T]$ the function $v[t]$ is of class $L^\infty\left((0,L)\right)$. Since $\rho^{(n)} \to \rho$ and $v^{(n)} \to v$ in $C^0\left([0,T];L^2\left((0,L)\right)\right)$ strongly we have also have $\rho^{(n)} \to \rho$ and $v^{(n)} \to v$ in $L^2\left((0,T)\times(0,L)\right)$ strongly. Then, Theorem 4.9 on page 94 in [1] implies that there exists a further subsequence of $\left\{\left(\rho^{(n)},v^{(n)}\right):n\geq 2\right\}$ for which $\rho^{(n)}(t,x) \to \rho(t,x)$ and $v^{(n)}(t,x) \to v(t,x)$ for $(t,x) \in (0,T)\times(0,L)$ a.e.. Consequently, we also obtain from (4.68) and (4.69) that

$$\frac{m}{b} \leq \rho(t,x) \leq \frac{m}{a}, \text{ for } (t,x) \in (0,T)\times(0,L) \text{ a.e.} \qquad (4.173)$$

$$|v(t,x)| \leq R_4, \text{ for } (t,x) \in (0,T) \times (0,L) \text{ a.e.} \qquad (4.174)$$

Lemma 4.3 and Lemma 4.4 and the above convergence properties imply that (2.15), (2.16) hold. Moreover, inequalities (4.171), (4.172), (4.173), (4.174) show that (2.19), (2.20) hold with $\rho_{\max} = \dfrac{m}{a}$, $v_{\max} = R_4$ and $\rho_{\min} = \dfrac{m}{b}$.

Lemma 4.5 and Lemma 4.6 and the above convergence properties allow us to conclude that (2.4) and (2.17) hold.

The only thing left to prove is estimate (2.18). Define for each fixed $t \in [0,T)$ and $h > 0$ with $t + h \leq T$ the functional $F : L^2\big((t, t+h) \times (0, L)\big) \to \mathbb{R}$:

$$F(u) = \frac{1}{2} \int_t^{t+h} \int_0^L \rho(s,x)\left(v(s,x) + \frac{\mu(\rho(s,x))}{\rho^2(s,x)} u(s,x)\right)^2 dxds + \int_t^{t+h} \int_0^L Q(\rho(s,x)) dxds \qquad (4.175)$$

Definitions (2.7), (2.10), (2.11), (4.175) guarantee that

$$F(\rho_x) = \int_t^{t+h} W(\rho[s], v[s]) ds \qquad (4.176)$$

Moreover, $F$ is convex and strongly continuous. Consequently, Corollary 3.9 on page 61 in [1], (4.176) and the fact that $\rho_x^{(n)} \to \rho_x$ in $L^\infty\big((0,T) \times (0,L)\big)$ weak star (which implies that $\rho_x^{(n)} \to \rho_x$ in $L^2\big((t, t+h) \times (0, L)\big)$ weakly) gives that

$$\int_t^{t+h} W(\rho[s], v[s]) ds \leq \liminf_{n \to +\infty} \left(F\big(\rho_x^{(n)}\big)\right) \qquad (4.177)$$

Exploiting the above convergence properties and (4.77) we can show that

$$\lim_{n \to +\infty} \left(F\big(\rho_x^{(n)}\big) - \int_t^{t+h} W(\rho^{(n)}[s], v^{(n)}[s]) ds\right) = 0 \qquad (4.178)$$

Combining (4.177), (4.178) and (4.131) with $t_0 = 0$, we get:

$$\int_t^{t+h} W(\rho[s], v[s]) ds \leq \liminf_{n \to +\infty} \left(\int_t^{t+h} W(\rho^{(n)}[s], v^{(n)}[s]) ds\right)$$
$$\leq \liminf_{n \to +\infty} \left(hW(\rho^{(n)}[0], v^{(n)}[0]) + h\frac{R_{15}}{n}\right) = h \liminf_{n \to +\infty} \left(W(\rho^{(n)}[0], v^{(n)}[0])\right) \qquad (4.179)$$

Consequently, we obtain from (4.179) and (4.132):

$$\frac{1}{h}\int_t^{t+h} W(\rho[s], v[s]) ds \leq \liminf_{n \to +\infty} \left(W_n(0) + \frac{R_{16}}{n}\right) = \liminf_{n \to +\infty} \left(W_n(0)\right) \qquad (4.180)$$

Estimate (2.18) is a direct consequence of (4.180), (4.67) and definitions (3.11), (3.14). The proof is complete. ◁

## 5. Concluding Remarks

The present paper proposed a novel particle scheme that provides convergent approximations of a weak solution of the Navier-Stokes equations for the 1-D flow of a viscous compressible fluid. Moreover, it is shown that all differential inequalities that hold for the fluid model are preserved by the particle method: mass is conserved, mechanical energy is decaying and a modified mechanical energy functional is also decaying. The proposed particle method can be used both as a numerical method and as a method of proving existence of solutions for compressible fluid models. The method can be extended easily to the case where external forces (e.g., gravity) are acting on the fluid and to the case where the tank containing the fluid is moving (see [15, 16, 17, 18]). On the other hand, the extension of the method to two and three spatial dimensions is demanding and will require new ideas.

The extension of the particle scheme to macroscopic models describing the flow of automated vehicles (see [14, 32]) can also be studied. However, additional issues have to be resolved since the macroscopic models given in [14, 32] are similar but not identical to the Navier-Stokes equations.

**Acknowledgments:** The authors would like to thank Dr. Dionysis Theodosis for his careful reading of the manuscript. The research leading to these results has received funding from the European Research Council under the European Union's Horizon 2020 Research and Innovation programme/ ERC Grant Agreement n. [833915], project TrafficFluid.